\documentclass{amsart}
\usepackage{amssymb,euscript,amsmath, mathrsfs}
\usepackage{bm}
\usepackage{rotating}
\usepackage[table]{xcolor}

\makeindex

\newcounter{ENUM}
\newcommand{\itm}{\item}
\newenvironment{ilist}[1][0]{\renewcommand{\theENUM}{\roman{ENUM}}\renewcommand{\itm}{\addtocounter{ENUM}{1}\item[(\theENUM)]}\begin{itemize}\setcounter{ENUM}{#1}}{\end{itemize}}

\newcommand{\margh}[1]{}

\input xy
\xyoption{all}
\CompileMatrices

\def\ZZ{{\mathbb Z}}

\def\PP{{\mathbb P}}

\def\cD{{\mathcal D}}

\def\cM{{\mathcal M}}

\def\sL{{\mathscr L}}

\def\sO{{\mathscr O}}

\def\fg{{\mathfrak g}}

\def\vp{\varphi}

\def\dv{\operatorname{div}}

\def\Pic{\operatorname{Pic}}

\def\ord{\operatorname{ord}}

\def\Sym{\operatorname{Sym}}

\def\md{\operatorname{md}}
\def\df{\operatorname{def}}

\newtheorem{thm}{Theorem}[section]
\newtheorem{prop}[thm]{Proposition}
\newtheorem{lem}[thm]{Lemma}
\newtheorem{cor}[thm]{Corollary}
\newtheorem{conj}[thm]{Conjecture}

\theoremstyle{definition}
\newtheorem{defn}[thm]{Definition}

\newtheorem{ex}[thm]{Example}
\newtheorem{sit}[thm]{Situation}
\theoremstyle{remark}
\newtheorem{notn}[thm]{Notation}
\newtheorem{rem}[thm]{Remark}

\numberwithin{equation}{section}

\begin{document}
\title[The strong maximal rank conjecture]{The strong maximal rank conjecture and moduli spaces of curves}
\author[F. Liu]{Fu Liu}
\author[B. Osserman]{Brian Osserman}
\author[M. Teixidor i Bigas]{Montserrat Teixidor i Bigas}
\author[N. Zhang]{Naizhen Zhang}

\begin{abstract} Building on recent work of the authors, we use degenerations
to chains of elliptic curves to prove two cases of the Aprodu-Farkas strong
maximal rank conjecture, in genus $22$ and $23$. This constitutes a major
step forward in Farkas' program to prove that the moduli spaces 
of curves of genus $22$ and $23$ are of general type. Our techniques involve 
a combination of the Eisenbud-Harris theory of limit linear series, and the 
notion of linked linear series developed by the second author.
\end{abstract}
\thanks{Fu Liu is partially supported by a grant from the
Simons Foundation \#426756. Brian
Osserman is partially supported by a grant from the Simons Foundation
\#279151. Naizhen Zhang is supported by the Methusalem Project Pure 
Mathematics at KU Leuven.}
\address[Fu Liu]{Department of Mathematics, University of California, 
One Shields Ave., Davis, CA 95616 USA}
\email{fuliu@math.ucdavis.edu}
\address[Brian Osserman]{Department of Mathematics, University of California, 
One Shields Ave., Davis, CA 95616 USA}
\email{osserman@math.ucdavis.edu}
\address[Montserrat Teixidor i Bigas]{Department of Mathematics, Tufts 
University, 503 Boston Ave., Medford, MA 02155 USA}
\email{montserrat.teixidoribigas@tufts.edu}
\address[Naizhen Zhang]{Department of Mathematics, KU Leuven, 
Celestijnenlaan 200B, Leuven, 3001, Belgium}
\email{naizhen.zhang@kuleuven.be}
\maketitle
\tableofcontents

\section{Introduction}

The study of the moduli space $\cM_g$ of curves of fixed genus $g$ is one of 
the most classical in algebraic geometry. Going back to Severi, based on
examples in low genus there was a general expectation that these moduli
spaces ought to be unirational. However, groundbreaking work of Harris
and Mumford and Eisenbud \cite{h-m2} \cite{ha6} \cite{e-h6} in the 1980's 
showed that not only is $\cM_g$ not unirational for large $g$, but it is in
fact of general type for $g \geq 24$. 
Their fundamental technique was to compute the classes of certain explicit
effective divisors on $\overline{\cM}_g$ arising from Brill-Noether theory, 
and use this to show that the canonical class of $\overline{\cM}_g$ can
be written as the sum of an ample and an effective divisor. The particular
families of divisors they considered were computable in all applicable
genera, but did not suffice to prove that $\cM_g$ is of general type for
$g \leq 23$. For the last $30$ years, no new cases have been proved of
$\cM_g$ being of general type. Roughly a decade ago, Farkas proposed (see 
\cite{fa2}, \cite[\S 4]{fa5}, \cite[\S 7]{fa4})
new families of expected divisors on $\cM_g$ as an approach to showing
that $\cM_{22}$ and $\cM_{23}$ of are general type. He computed `virtual
classes' for these expected divisors in \cite{fa4} for genus $22$ and 
in the new paper \cite{fa3} for genus $23$, and in both cases found that
the classes satisfy the necessary inequalities to conclude that $\cM_{22}$
and $\cM_{23}$ are of general type, provided that they are indeed 
represented by effective divisors.\footnote{\cite{fa4} also includes an 
announcement that $\cM_{22}$ is of general type, but in \cite{fa3} Farkas 
states that the intended proof ``has not materialized.''}

The remaining steps of Farkas' program can be described in terms of the 
following definition.

\begin{defn}\label{defn:divisors} Given $g \geq 21$, let
$d=3+g$. Let $\cD_{g} \subseteq \cM_g$
consist of curves $X$ of genus $g$ which admit a $\fg^6_d$ such that the
resulting image of $X$ in $\PP^6$ lies on a quadric hypersurface.
\end{defn}

For $g=22$ or $23$, in order to conclude that $\cM_g$ is of general type,
one has to check two statements: first, that $\cD_g$ yields an effective 
divisor, or equivalently, that $\cD_g \subsetneq \cM_g$; and second, that the 
class induced by $\cD_g$ agrees with the class computed by Farkas, or 
equivalently, that the subset of $\cD_g$ consisting of curves carrying 
\emph{infinitely} many 
$\fg^6_d$s whose image lie on a quadric occurs in codimension strictly 
higher than $1$.
In this paper, we prove the first of these two statements, for both
$g=22$ and $g=23$. An independent proof of this result has been obtained
by Jensen and Payne in \cite{j-p4}, using their tropical approach.
Our main theorem is thus the following:

\begin{thm}\label{thm:main} In characteristic $0$, the loci $\cD_{22}$ and 
$\cD_{23}$ are proper subsets of $\cM_{22}$ and $\cM_{23}$ respectively.
\end{thm}

In fact, we show further that the closure of $\cD_{22}$ does not contain the
locus of chains of genus-$1$ curves; see Theorem \ref{thm:rho-1} below.
In addition, our proof goes through unmodified for characteristic 
$p \geq 29$, and our techniques can in principle be applied to lower 
characteristics as well, but due to characteristic restrictions on the
application to the geometry of $\cM_g$, we have not pursued this. See
Remark \ref{rem:char} below.

The genus $21$ case of Theorem \ref{thm:main} is a special case of the 
classical ``maximal rank
conjecture,'' and was proved by Farkas in \cite{fa2} as part of an
infinite family of divisors of small slopes. However, in this case the
divisor in question does not quite satisfy the necessary inequality to
obtain that $\cM_{21}$ is of general type. With applications to moduli
spaces of curves in mind, Aprodu and Farkas proposed in
Conjecture 5.4 of \cite{a-f1} a ``strong maximal rank conjecture,'' of
which Theorem \ref{thm:main} constitutes two cases. These conjectures
study ranks of multiplication maps. Specifically, given a linear series
$(\sL,V)$ on a curve $X$, we have the multiplication map
\begin{equation}\label{eq:mult-map}
\Sym^2 V \to \Gamma(X,\sL^{\otimes 2}).
\end{equation}
Note that the source has dimension $\binom{r+2}{2}$, and assuming $X$
is Petri-general, the target has dimension $2d+1-g$.
The image of $X$ under the linear series lies on a quadric if and only
if \eqref{eq:mult-map} has a nonzero kernel. The classical maximal rank
conjecture asserts that if $r \geq 3$, for a general $X$ and a general 
$\fg^r_d$ on $X$, the map
\eqref{eq:mult-map} should always be injective or surjective (and 
similarly for the higher-order multiplication maps). Many special cases
of this were proved by various people; we omit discussion of most of 
these, but mention that the case of quadrics was first proved by
Ballico \cite{ba3}, and subsequent proofs were given by Jensen and Payne 
using a tropical approach \cite{j-p2}, and by the present authors 
\cite{o-l-t-z} using
a degeneration to a chain of genus-$1$ curves. Very recently, Larson
has proved the full classical maximal rank conjecture \cite{la5}.

In contrast, the strong maximal rank conjecture remains wide open, even
in the case of quadrics. Since the failure of \eqref{eq:mult-map} to have
maximal rank is a determinantal condition, the strong maximal rank
conjecture of Aprodu and Farkas (Conjecture 5.4 of \cite{a-f1}) is the 
following:

\begin{conj}[Aprodu-Farkas]\label{conj:strong-mrc} 
Set $\rho:=g-(r+1)(g+r-d)$.

On a general curve of genus $g$, if $\rho<r-2$,
the locus of $\fg^r_d$s for which \eqref{eq:mult-map} fails to have
maximal rank is equal to the expected determinantal codimension, which
is $1+\left|\binom{r+2}{2}-(2d+1-g)\right|$. In particular, when this
expected codimension exceeds $\rho$, every linear series
on $X$ should have maximal rank.\footnote{In fact, Aprodu and Farkas also
include higher-degree multiplication maps in their conjecture. Farkas
and Ortega \cite{f-o1} subsequently relax the $\rho<r-2$ hypothesis in cases 
such as ours, where $\rho$ is less than the expected codimension.} 
\end{conj}

For the family of cases considered in Definition \ref{defn:divisors},
we compute that $\rho=g-21=(2d+1-g)-\binom{r+2}{2}$, so in
this case Conjecture \ref{conj:strong-mrc} predicts that every
linear series should yield \eqref{eq:mult-map} of maximal rank, and more 
specifically, should have injective multiplication map, just as we prove
in Theorem \ref{thm:main} for the cases $g=22,23$.

Our proof builds on the ideas introduced in \cite{o-l-t-z}, which
combine the Eisenbud-Harris theory of limit linear series with ideas
from the theory of linked linear series introduced by the second author
in \cite{os8} and \cite{os20}. The idea is to start with a limit linear
series on a chain $X_0$ of genus-$1$ curves, and describe a collection of 
global sections living in different multidegrees on $X_0$. We then take
tensors of these sections and consider their image in a carefully chosen
multidegree, showing that they have the correct-dimensional span.
The first major difficulty in moving from the classical maximal rank 
conjecture to the strong maximal rank conjecture is that instead of being
able to work with a single limit linear series in each case, we have to
consider all possible limit linear series. We are able to overcome this
for the full infinite family of cases described in Definition 
\ref{defn:divisors}, and we expect that these ideas should extend to
cover an infinite sequences of similarly constructed (infinite) families.

The more serious difficulty is that for certain degenerate limit linear
series (which occur already in codimension $1$), we do not have very good
control over global sections occurring in the expected multidegrees when
we have a family of linear series on smooth curves specializing to the
given limit linear series. This can be expressed in terms of trying to
understand the possible linked linear series lying over the given limit
linear series. Even for these limit linear series, a nonempty open subset 
of the possible linked linear series will always behave well, but there
will be some cases which are more slippery. To overcome this, we 
systematically use ideas from linked linear series to prove that when
$\rho=1$ or $\rho=2$ we can always produce global sections of certain 
prescribed forms which must lie in the specialization of the family of
linear series. For the case $\rho=1$ (i.e., genus $22$), this leads to
a relatively brief proof of Theorem \ref{thm:main}. However, already 
for $\rho=2$ the situation is quite a bit more complicated. To handle
the degenerate cases, we consider variant multidegrees which depend 
more tightly on the limit linear series in question, and (partially
inspired by the earlier work \cite{j-p2} of Jensen and Payne on a
tropical approach to the classical maximal rank conjecture) we also consider 
families of curves 
with highly specialized directions of approach, which gives us further
control over the behavior of the global sections in different multidegrees.

We expect that the tools we develop here will lead to proofs of 
infinite families of the strong maximal rank conjecture, and have written
the different parts of the argument to be independent of $r$ and/or $\rho$
wherever this does not lead to unnecessary complication. The nature of our 
approach also allows for proving cases 
of the maximal rank conjecture where the expected codimension does not 
exceed $\rho$, so that the locus of linear series which do not have maximal 
rank is nonempty. Our approach should also
be useful in other questions involving multiplication maps for linear
series, such as the conjecture of Bakker and Farkas (Remark 14 of 
\cite{b-f3}),
which was motivated by connections to higher-rank Brill-Noether
theory. Their conjecture treats a certain specific family of cases, but 
with products of distinct linear series in place of symmetric squares of a 
fixed one. 
In addition, our work in \S \ref{sec:nondegen} on nondegeneracy
of certain morphisms from genus-$1$ curves to projective spaces and in
\S \ref{sec:degen-link} on the structure of exact linked linear series is
likely to be useful in other settings as well.

The structure of the paper is as follows. In \S \ref{sec:nondegen} we
analyze certain maps from genus-$1$ curves to projective spaces
which arise naturally from tensor squares of linear series, and show that
these are nondegenerate morphisms in cases of interest. In \S \ref{sec:bg},
we review the Eisenbud-Harris theory of limit linear series, and the 
related theory of linked linear series introduced by the second author.
In \S \ref{sec:degen-link}, we analyze the possible structures of linked
linear series lying over a given limit linear series in the cases that can
arise when $\rho \leq 2$. In \S \ref{sec:setup} and \S \ref{sec:crit},
we analyze a certain collection of sections which arise from taking the
tensor square of a limit linear series, and give an elementary criterion 
for them to be linearly independent. In \S \ref{sec:r6}, we apply this
criterion to a family of examples with $r=6$, which include the genus-$22$
and genus-$23$ cases of interest for the proof of Theorem \ref{thm:main}.
Finally, in \S \ref{sec:proofs}, we put together the analysis of the 
structure of linked linear series with the independence results of
\S \ref{sec:r6} to complete the proof of Theorem \ref{thm:main}.

\begin{rem}\label{rem:char}
We mention that although we impose characteristic-$0$ hypotheses in our
main theorem, these do not appear to be essential. Nearly everything we
do is characteristic-independent, but we use a characteristic-dependent
result (Theorem \ref{thm:refined-specialize} below)
of Eisenbud and Harris to simplify the situation slightly by restricting
our attention to ``refined'' limit linear series (Definition 
\ref{def:eh-lls} below). In fact, the only characteristic dependence in
Theorem \ref{thm:refined-specialize} is the use of the Pl\"ucker inequality,
which still holds in characteristic $p$ and degree $d$ when $p>d$; see
for instance Proposition 2.4 and Corollary 2.5 of \cite{os8}. Thus, our 
proof of Theorem \ref{thm:main} extends as written to characteristic 
$p >25$ for $g=22$ and $p>26$ for $g=23$.

Moreover,
since our key specialization result (Proposition \ref{prop:exact-limit} 
below) on linked linear series applies in arbitrary characteristic, there
is no visible obstruction to extending our proof to lower characteristics as 
well.
However, key portions of the argument for the implications for the geometry 
of $\cM_g$ were written using characteristic $0$, and as far as we are
aware no one has carefully analyzed which positive characteristics they may
apply to, so for the present paper it seems preferable
to work in the simpler setting.
\end{rem}

\subsection*{Acknowledgements} We would like to thank Gavril Farkas, Dave 
Jensen and Sam Payne for helpful contextual conversations.

\section{Nondegeneracy calculations}\label{sec:nondegen}

In this section, we study maps from elliptic curves to projective space
determined by comparing values of tensor products of certain tuples of 
sections at two points
$P$ and $Q$. We will need two distinct results in this direction:
first, we consider the situation that we let the point $Q$ vary. This
is already considered in \cite{o-l-t-z}, where we showed that these maps 
are morphisms, described them explicitly, and gave partial criteria for 
nondegeneracy. Here we extend the nondegeneracy criterion to a sharp 
statement for the case of tensor pairs. This is used to show that if we vary 
the location of the nodes on individual components, we can get possible
linear dependencies to vary sufficiently nontrivially. Next, we will 
consider a new case, where $Q$ is fixed, but we have a separate varying
parameter. This situation was not considered in \cite{o-l-t-z}, but will
be important to us in dealing with situations where the discrete data
of the limit linear series does not fix the underlying line bundle in some
components.

First, given a genus-$1$ curve $C$ and distinct $P,Q$ on $C$, and 
$c,d \geq 0$, let
$\sL=\sO_C(cP+(d-c)Q)$. Then for any $a,b \geq 0$ with $a+b=d-1$, there is
a unique section (up to scaling by $k^\times$) of $\sL$ vanishing to order
at least $a$ at $P$ and at least $b$ at $Q$. Thus, we have a uniquely
determined point $R$ such that the divisor of the aforementioned
section is $aP+bQ+R$; explicitly, $R$ is determined by
$aP+bQ+R \sim cP+(d-c)Q$, or
\begin{multline}\label{eq:r-formula}
R \sim (c-a)P+(d-c-b)Q = (c-a)P+(1+a-c)Q\\
=P + (a+1-c)(Q-P)=Q+(a-c)(Q-P).
\end{multline}
Thinking of $C$ as a torsor over $\Pic^0(C)$, we see that $R=P$ if and
only if $Q-P$ is $|a+1-c|$-torsion, and $R=Q$ if and only if $Q-P$ is
$|a-c|$-torsion. Note that \eqref{eq:r-formula} makes sense even when
$Q=P$ (in which case $R=Q=P$), so we will use the formula for all $P,Q$,
understanding that it has the initial interpretation as long as $Q \neq P$.
To avoid trivial cases, we will assume that $a \neq c-1$,
and $b \neq d-c-1$, or equivalently, $a+1-c \neq 0$, and $a-c \neq 0$.

\begin{sit}\label{sit:nondegen-1} Fix $\ell \geq 1$, and for
$j=0,\dots,\ell$, set numbers $a_1^j,a_2^j$, $b_1^j,b_2^j$
satisfying:
\begin{itemize}
\item $a_i^j+b_i^j=d-1$ for all $i,j$;
\item $a_i^j - c \neq 0,-1$ for all $i,j$;
\item $a_1^j+a_2^j$ is independent of $j$.
\end{itemize}
\end{sit}

We now have sections
$s_i^j$ with divisors $a_i^jP+b_i^j Q +R_i^j$, and
forming tensor products yields sections
$s^j = s_1^j \otimes s_2^j \in \Gamma(C,\sL^{\otimes 2})$,
with divisors
$$\left(a_1^j+a_2^j\right) P + \left(b_1^j+b_2^j\right) Q + R_1^j + R_2^j,$$
having the property that any two $R_1^j+R_2^j$ are linearly equivalent.
Now, if $Q-P$ is not $|a_i^j+1-c|$-torsion for any $i,j$, we can normalize
the $s^j$, uniquely up to simultaneous scalar, so that their values at $P$
are all the same. Then provided that there is some $j$ such that
$Q-P$ is not $|a_i^j-c|$-torsion for any $i$, considering
$(s^0(Q),\dots,s^{\ell}(Q))$ gives a well-defined point of $\PP^{\ell}$.

\begin{notn} With discrete data as in Situation \ref{sit:nondegen-1},
suppose $P$ is fixed. For a given $Q\in C$, denote by $R_i^{j,Q}$ the point
determined as above by $P$ and $Q$, and by $f_Q$ the point of $\PP^{\ell}$
determined by $(s^0(Q),\dots,s^{\ell}(Q))$.
Let $U$ be the open subset of $C$ consisting of all $Q$ such that $Q-P$ is
not $|a_i^j-c|$-torsion or $|a_i^j+1-c|$-torsion for any
$i,j$.
\end{notn}

In \cite{o-l-t-z} we showed that 
the map $U \to \PP^{\ell}$ given by
$Q \mapsto f_Q$ extends to a morphism $f:C \to \PP^{\ell}$.

Our main result is then the following, extending Corollary 2.7 of
\cite{o-l-t-z}.

\begin{prop}\label{prop:nondegenerate} 
If $C$ is not supersingular, all the $a_i^j$ are distinct, and
$a_1^j+a_2^j \neq 2c-1$, then $f$ is nondegenerate.
\end{prop}

The proof relies on reduction to a good understanding of the $\ell=1$
case. Indeed, we can view our map as being given by
$(1,f_1,\dots,f_{\ell})$, where $f_j$ is the rational function constructed
from the sections $s^j,s^0$ (note that we are switching the order of $0$
and $j$, because we are dividing through all terms by $s^0$).
Thus, nondegeneracy is equivalent to linear
independence of the rational functions $1,f_1,\dots,f_{\ell}$, whose zeroes
and poles are described explicitly by the following result, which combines
Lemma 2.3 and Corollary 2.4 of \cite{o-l-t-z}.

\begin{lem}\label{lem:nondegen} In the $\ell=1$ case of 
Situation \ref{sit:nondegen-1}, the
function $f:U \to k^\times$ given by $Q \mapsto (s^0/s^1)(Q)$
determines a rational function on $C$. We then have
\begin{multline*} \dv f  =
\sum_{i=1}^2 ((P+\Pic^0(C)[|a_i^0-c|])-(P+\Pic^0(C)[|a_i^1-c|])
\\-(P+\Pic^0(C)[|a_i^0+1-c|])+(P+\Pic^0(C)[|a_i^1+1-c|])),
\end{multline*}
where for a divisor $D=\sum_j c_j P_j $ on $\Pic^0(C)$, the notation
$P+D$ indicates the divisor $\sum_j c_j (P+P_j)$ on $C$, using the
$\Pic^0(C)$-torsor structure on $C$.

Moreover, if $C$ is not supersingular, $f$ is nonconstant if and only if 
$a_1^j+a_2^j \neq 2c-1$.
\end{lem}

In the above, the torsion subgroups $\Pic^0(C)[n]$ should be equipped with
the multiplicities arising from the inseparable degree of the appropriate
multiplication map. Thus, if $k$ has characteristic $0$ or characteristic
$p$ not dividing $n$, then $\Pic^0(C)[n]$ is a reduced divisor, but
otherwise all the points of $\Pic^0(C)[n]$ have coefficients given by the
appropriate power of $p$.

\begin{proof}[Proof of Proposition \ref{prop:nondegenerate}]
By Lemma \ref{lem:nondegen}, $f_1,\dots,f_{\ell}$ are all non-constant.
By re-indexing the pairs we may further assume 
$$a^{\ell}_1<a^{\ell-1}_1<\dots<a^0_1\leq 
a^0_2<a^1_2<\dots<a^{\ell}_2.$$

Let $n^{j}_i:=|a^j_i-c+1|$, $m^j_i:=|a^j_i-c|$, and
$n^j=\max\{n^j_1,n^j_2\}$. 

A first observation is that $n^j>n^{j-1}$ for all $j$: if $a^{j-1}_1<c$ 
(respectively, $a^{j-1}_1>c$), then $n^{j-1}_1<n^j_1$ (respectively, 
$n^{j-1}_1<n^j_2$), and thus $n^{j-1}_1<n^j$; by a similar calculation, 
$n^{j-1}_2<n^j$; thus, $n^{j-1}<n^j$.

A second observation is that $n^j\geq\max\{m^0_1,m^0_2\}$ for all $j\geq 1$,
and if equality is attained, $j$ must be $1$. Indeed, when $c<a^0_2$, we have 
$m^0_2<m^j_2<n^j_2$ for all $j\geq 1$; meanwhile, either 
$m^0_1\leq m^0_2$ (if $c<a^0_1$) or $m^0_1\leq n^1_1<n^j_1$ (if $c>a^0_1$) 
for all $j>1$; thus, $\max\{m^0_1,m^0_2\}\leq n^j$ for all $j\geq 1$. 
When $c>a^0_2$, $m^0_2\leq m^0_1\leq n^1_1<n^j_1$ for all $j>1$, and hence the
same conclusion holds. 

Now, we claim that $f_j$ has poles at the strict $n_j$-torsion points. 
Recalling from Lemma \ref{lem:nondegen} that the poles of $f_j$ are
supported among the $m_i^0$- and $n_i^j$-torsion points for $i=1,2$, the
above two observations show that $1,\dots,f_{j-1}$ cannot have any poles
at strict $n^j$-torsion points, which immediately implies that
$1,f_1,\dots,f_{\ell}$ are $k$-linearly independent. Thus, it suffices
to prove the claim. Since the potential zeroes of $f_j$ are supported among
the $m_i^j$- and $n_i^0$-torsion points, we just need to show that
$n^j$ does not divide $m_i^j$ or $n_i^0$ for $i=1,2$ and any $j \geq 1$.
Moreover, we already know that $n^j>n^0 \geq n_i^0$, so it is enough to
consider the $m_i^j$. We consider two cases.

\textbf{Case 1:} $c<a^0_2$, so that also $c<a^j_2$ for all $j$. In this 
case, 
$$m^0_2<n^0_2\leq m^1_2<n^1_2\leq\dots \leq m^{\ell}_2<n^{\ell}_2.$$
In particular, we have $n^j>m^j_2$, so
it remains to compare $n^j$ against $m^j_1$. If 
$n^j=n^j_1$, since $n^j_1$ is always coprime to $m^j_1$, the claim follows 
instantly. If $n^j=n^j_2>n^j_1$, since $|n^j_1-m^j_1|=1$, we have 
$n^j\geq m^j_1$. But equality cannot hold as it would imply that 
$a^j_1+a^j_2=2c-1$, which is ruled out by our assumption. So we conclude
the claim in this case. 

\textbf{Case 2:} $c>a^0_1$, so that $c>a^j_1$ for all $j$. 
If $a^j_2>c$, $n^j_2=m^j_2+1$ and hence $n^j>m^j_2$. Meanwhile, 
$n^j_1=m^j_1-1$. Similarly to the previous case, either $n^j>m^j_1$ or 
$n^j$ is coprime to $m^j_1$, and the claim follows. If $a^j_2<c$, 
$n^j_2=m^j_2-1$. Under our assumption, $n^j=n^j_1$ so is coprime 
to $m^j_1$. But because $j\geq 1$, we have $n^j_1 \geq n^j_2+2 = m^j_2+1$,
so $n^j_1>m^j_2$ and the claim follows.
\end{proof}

We now move on to the new situation, where our point $Q$ will be fixed,
but our line bundle $\sL$ is allowed to vary. If we have $a,b$ with 
$a+b=d-1$, then the isomorphism class of a line bundle $\sL$ of degree $d$ 
can be uniquely determined by a point $R$ by setting 
$\sL \cong \sO_C(aP+bQ+R)$.  If we have 
$a',b'$ with $a'+b'=d-1$ also, and $\sO_C(aP+bQ+R)\cong \sO_C(a'P+b'Q+R')$,
then we find that 
$$R'=R+(a-a')(P-Q).$$

We fix discrete data as before, except that since $\sL$ will vary, we
do not have any $c$.

\begin{sit}\label{sit:nondegen-2} Fix $\ell \geq 1$, and for
$j=0,\dots,\ell$, set nonnegative integers $a_1^j,a_2^j$, $b_1^j,b_2^j$
satisfying:
\begin{itemize}
\item $a_i^j+b_i^j=d-1$ for all $i,j$;
\item $a_1^j+a_2^j$ is independent of $j$.
\end{itemize}
\end{sit}

First suppose we fix $\sL$. As before, we have sections $s_i^j$ of $\sL$ 
with divisors $a_i^j P+b_i^j Q +R_i^j$, and we can take tensor products
to obtain $s^j=s_1^j \otimes s_2^j$ having divisors 
$(a_1^j+a_2^j) P+(b_1^j+b_2^j) Q +R_1^j+R_2^j$. Note that the divisors
$R_1^j+R_2^j$ will all be linearly equivalent to one another by construction.
If we assume that none
of the $R_i^j$ are equal to $P$ (which will be the case provided that
$\sL$ and $Q$ are general relative to $P$), we can normalize the $s_j$
to have the same value at $P$, and then we obtain a well-defined point
$(s^0(Q),\dots,s^{\ell}(Q)) \in \PP^{\ell}$. But because we have said that
$\sL$ is uniquely determined by $R_1^0$, we can view this procedure as
giving a rational map from $C$ to $\PP^{\ell}$, which we will now study.
The argument will be similar to that of Lemma 2.3 and Corollary 2.7 of
\cite{o-l-t-z}, but a bit simpler.

\begin{prop}\label{prop:nondegen-again} Suppose that $P-Q$ is not
$m$-torsion for any $m \leq d$, and let $U \subseteq C$ consist of the
open subset of points not differing from $P$ by $m$-torsion for any 
$m \leq d$. Let $\vp:U \to \PP^{\ell}$ be the map which uses the above 
construction with $R_1^0=R$ to send $R\in U$ to $(s^0(Q),\dots,s^{\ell}(Q))$.
Then $\vp$ extends to a nondegenerate morphism
$C \to \PP^{\ell}$.
\end{prop}

\begin{proof} We first consider the case $\ell=1$, proving that we 
obtain a nonconstant rational function, and showing further that the
divisor of this function is equal to 
\begin{multline*} Q+(Q-(a_1^0-a_2^0)(P-Q)) +(P-(a_1^0-a_1^{1})(P-Q))
+(P-(a_1^0-a_2^{1})(P-Q))
\\ -(Q-(a_1^0-a_1^1)(P-Q))-(Q-(a_1^0-a_2^1)(P-Q))
-P -(P-(a_1^0-a_2^0)(P-Q)).\end{multline*}

Consider $R_i^j$ as a divisor on $C \times C$ by setting 
$R_i^j$ to be the graph of the morphism 
$$R \mapsto R+(a_1^0-a_i^j)(P-Q)$$
(so that $R_1^0$ is simply the diagonal, and in each fiber over $R$ we 
obtain our original point $R_i^j$ for the case that $\sL$ is determined by
setting $R_1^0=R$).
Set 
$$D^j=R_1^j+R_2^j+(P-(a_1^0-a_1^{1-j})(P-Q)) \times C 
+(P-(a_1^0-a_2^{1-j})(P-Q)) \times C$$
for $j=0,1$. Then we claim that $D^0$ and $D^1$ are linearly equivalent.
By construction, if we restrict to $\{R\} \times C$ for any $R$ not among
the $P-(a_1^0-a_i^j)(P-Q)$, we get that $D^0$ and $D^1$ are linearly
equivalent,
so $D^0 - D^1 \sim D \times C$ for some divisor $D$ on $C$. But if we
restrict to $C \times \{P\}$, we see that $R_1^j+R_2^j$ restricts to
$(P-(a_1^0-a_1^j)(P-Q))+(P-(a_1^0-a_2^j)(P-Q))$, so the restrictions of
$D^0$ and $D^1$ are linearly equivalent on $C \times \{P\}$, and hence on
$C \times C$, as desired. Moreover, this shows that if $t_0$ and $t_1$
are sections of the resulting line bundle having $D^0$ and $D^1$ as
divisors, then $t_0|_{C \times \{P\}}$ has the same divisor as
$t_1|_{C \times \{P\}}$, so we can scale so that $t_0$ and $t_1$ are
equal on $C \times \{P\}$.\footnote{In fact, we see that the divisors 
$R_1^j+R_2^j$ are already linearly equivalent, but we need the given 
definition of the $D^j$ precisely so that $t_0$ and $t_1$ can be normalized 
as desired along $P$.} We then see that our map $U \to \PP^1$ is
given by composing $R \mapsto (R,Q)$ with the rational function induced
by our normalized choice of $(t_0,t_1)$. Thus, it is a rational function,
as desired. We compute its divisor simply by looking at the restrictions
of $D^0$ and $D^1$ to $C \times \{Q\}$, which gives the claimed formula.

Now, for the case of arbitrary $\ell$, we can consider the map to 
$\PP^{\ell}$ to be given by a tuple of rational functions induced from
the $\ell=1$ case, specifically by $(f_0,\dots,f_{\ell-1},1)$, where
$f_j$ comes from looking at $s^j$ and $s^{\ell}$. To show nondegeneracy, 
it suffices to show that the $f_j$ are linearly independent, which we do by 
showing that each of them (other than $f_{\ell}=1$) has a pole which none 
of the others have. If we order so that 
$$a_1^0 < a_1^1<\dots<a_1^{\ell} \leq a_2^{\ell} < a_2^{\ell-1} 
< \dots < a_2^0,$$ 
we see that $P-(a_1^j-a_2^j)(P-Q)$ occurs among the poles of $f_j$: indeed,
given our non-torsion hypothesis on $P-Q$, the only positive term in the
divisor which could possibly cancel it is $Q$, which would require
$a_1^j-a_2^j=1$, which is not possible with our above ordering.
But again using our nontorsion hypothesis, and the fact that 
$a_2^j-a_1^j$ strictly decreases as $j$ increases, we see that we obtain
the desired distinct poles.
\end{proof}

\section{Background on limit linear series and linked linear series}\label{sec:bg}

In this section we review background on limit linear series, as introduced
by Eisenbud and Harris in \cite{e-h1}, and on linked linear series, 
introduced by the second author in \cite{os8} for two-component curves
and generalized to arbitrary curves of compact type in 
\cite{os20}.\footnote{In \cite{os8}, linked linear series were called
`limit linear series,' but the name was changed subsequently to reduce
confusion.} Recall that a curve of \textbf{compact type} is a projective
nodal curve such that every node is disconnecting, or equivalently, the
dual graph is a tree. To streamline our presentation, we will largely 
restrict our attention to the situation of curves of compact type together 
with one-parameter smoothings.

\begin{defn}\label{def:eh-lls} Let $X_0$ be a curve of compact type,
with dual graph $\Gamma$.
Given $r,d \geq 0$, a \textbf{limit linear series} on $X_0$ of dimension
$r$ and degree $d$ is a tuple $(\sL^v,V^v)_{v \in V(\Gamma)}$, where
each $(\sL^v,V^v)$ is a linear series of dimension $r$ and degree $d$
on the component $Z_v$ of $X_0$ corresponding to $v$. This tuple is
further required to satisfy the following condition: if $Z_v$ and $Z_{v'}$
meet at a node $P_e$, and $a^{(v,e)}_{\bullet}$ and $a^{(v',e)}_{\bullet}$
are the vanishing sequences at $P_e$ of $(\sL^v,V^v)$ and 
$(\sL^{v'},V^{v'})$ respectively, then 
$$a^{(v,e)}_j+a^{(v',e)}_{r-j} \geq d \quad \text{ for }j=0,\dots,r.$$

A limit linear series is said to be \textbf{refined} if the above
inequalities are equalities for all $e$ and $j$.
\end{defn}

We now consider a one-parameter smoothing of $X_0$, as follows.

\begin{sit}\label{sit:smoothing} Suppose $B$ is the spectrum of a 
discrete valuation ring with algebraically closed residue field, and 
$\pi:X \to B$ is flat and proper, with
special fiber $X_0$ a curve of compact type, and smooth generic fiber
$X_{\eta}$. Suppose further that the total space $X$ is regular, 
that $\pi$ admits a section. 
\end{sit}

Now, suppose we have a line bundle $\sL_{\eta}$ generically -- more
precisely, we allow for the possibility that $\sL_{\eta}$ is only
defined after a finite extension of the base field of $X_{\eta}$. 
We can then take a finite base change $B' \to B$ so that $\sL_{\eta}$
is defined over $X'_{\eta}$, and then $X'$ may not be regular, but
the line bundle $\sL_{\eta}$ will still extend over $X_0$ because
$X_0$ is of compact type. Moreover, there is a unique extension of
$\sL_{\eta}$ having any specified \textbf{multidegree} (i.e., tuple of 
degrees on each component) adding up to $d$: because $X$ was regular each
component $Z_v$ of $X_0$ is a Cartier divisor in $X$, and twisting by
the $\sO_X(Z_v)$ (or more precisely, their pullbacks to $X'$) will
increase the degree by $1$ on each component meeting $Z_v$, and 
decrease the degree on $Z_v$ correspondingly. For a multidegree $\omega$,
we denote this unique extension by $\widetilde{\sL}_{\omega}$. In particular, 
for each $Z_v$, we can consider the multidegree $\omega^v$ which concentrates 
degree $d$ on $Z_v$, and has degree $0$ elsewhere. 

Eisenbud and Harris (Proposition 2.1 of \cite{e-h1}) show the following 
specialization result:

\begin{prop}[Eisenbud-Harris]\label{prop:eh-specialize}
Given a linear series $(\sL_{\eta},V^v)$
on $X'_{\eta}$ of dimension $r$ and degree $d$, if we set 
$\sL^v:= (\widetilde{\sL}_{\omega^v})|_{Z_v}$, and 
$V^v:=(V_{\eta} \cap \Gamma(X',\widetilde{\sL}_{\omega^v}))|_{Z_v})$,
then the resulting tuple $(\sL^v,V^v)_v$ is a limit linear series on $X_0$.
\end{prop}

They also show (Theorem 2.6 of \cite{e-h1}) the following:

\begin{thm}[Eisenbud-Harris]\label{thm:refined-specialize} In characteristic
$0$, after finite base change and blowing up nodes in the special fiber, 
we may assume that the specialized limit linear series constructed by 
Proposition \ref{prop:eh-specialize} is refined.
\end{thm}

Note that the only effect on $X_0$ of the base change and blowup is that
chains of genus-$0$ curves are introduced at the nodes. Assuming we blow
up to fully resolve the singularities resulting from the base change,
these chain of
curves have length equal to one less than the ramification index of the
base change, so in particular they are the same at every node.

We now move on to linked linear series. The first observation is that if
we have two multidegrees $\omega$ and $\omega'$, then there is a unique 
collection of nonnegative coefficients $c_v \in \ZZ$, not all positive, such 
that 
$\widetilde{\sL}_{\omega} \cong \widetilde{\sL}_{\omega'}(-\sum_v c_v Z_v)$. 
In this way, we obtain an inclusion 
$\widetilde{\sL}_{\omega} \hookrightarrow \widetilde{\sL}_{\omega'}$ which is 
defined uniquely up to scaling. If we
define $\sL_{\omega}:=\widetilde{\sL}_{\omega}|_{X_0}$, we get induced maps
$\sL_{\omega} \to \sL_{\omega'}$ which are no longer injective, as they vanish
identically on the components $Z_v$ with $c_v>0$. However, they are
injective on the remaining components. Passing to global sections we
obtain maps
$$f_{\omega,\omega'}:\Gamma(X_0,\sL_\omega) \to \Gamma(X_0,\sL_{\omega'}).$$
From the construction we see that 
$f_{\omega,\omega'} \circ f_{\omega',\omega}$ always
vanishes identically. Although the twisted line bundles $\sL_{\omega}$ can be
described intrinsically on the special fiber, the maps 
$f_{\omega,\omega'}$ depend
on the smoothing of $X_0$ whenever the locus on which they are nonvanishing
is disconnected.

To minimize notation, we will define linked linear series only in the
above specialization context.

\begin{defn}\label{defn:linked-ls} Given $\sL_{\eta}$ of degree $d$ and the 
induced tuple $(\sL_{\omega})_{\omega}$ of line bundles, a \textbf{linked 
linear series} of dimension $r$ (and degree $d$) on the $\sL_{\omega}$ is a 
tuple $(V_{\omega})_{\omega}$ for all multidegrees of total degree $d$ where 
each $V_{\omega} \subseteq \Gamma(X_0,\sL_{\omega})$ is an 
$(r+1)$-dimensional space of global sections, and for every $\omega,\omega'$,
we have
$$f_{\omega,\omega'} (V_{\omega}) \subseteq V_{{\omega}'}.$$
\end{defn} 

We then see easily from the definitions that we have:

\begin{prop}\label{prop:linked-specialize} 
If we have $(\sL_{\eta},V_{\eta})$ generically, and for all $\omega$ we set
$V_{\omega}=(V_{\eta} \cap \Gamma(X',\widetilde{\sL}_{\omega}))|_{X_0}$, we
obtain a linked linear series. 
\end{prop}

Moreover, this process is visibly compatible with the Eisenbud-Harris
specialization process, and we have a forgetful map which visibly commutes 
with specialization:

\begin{thm}\label{thm:linked-forget} If $(V_{\omega})_{\omega}$ is a linked linear
series on $\sL_{\omega}$, and we set $\sL^v=\sL_{\omega^v}|_{Z_v}$ and 
$V^v=V_{\omega^v}|_{Z_v}$ for all $v \in V(\Gamma)$, then $(\sL^v, V^v)$
is a limit linear series.
\end{thm}

This is explicitly stated (in the generality of higher-rank vector bundles)
as part of Theorem 4.3.4 of \cite{os20}, but
is primarily a consequence of Lemma 4.1.6 of \textit{loc.\ cit.}

In \cite{os20}, the following notion is introduced:

\begin{defn} A linked linear series is \textbf{simple} if there exist
multidegrees $\omega_0,\dots,\omega_r$ and sections 
$s_j \in \Gamma(X_0,\sL_{\omega_j})$ such that for every $\omega$, the 
$f_{\omega_j,\omega} (s_j)$ form a basis of $V_{\omega}$.
\end{defn}

The simple linked linear series form an open subset, and are particularly 
easy to understand (hence the name). However, we will be forced to consider
more general linked linear series arising under specialization. We
therefore introduce the following open subset, originally introduced in
\cite{os8} in the two-component case.

\begin{defn} A linked linear series is \textbf{exact} if for every
multidegree $\omega$, and every proper subset $S \subseteq V(\Gamma)$,
if $\sL_{\omega'} \cong \sL_{\omega}(-\sum_{v \in S} Z_v)$, then
$$f_{\omega,\omega'} (V_{\omega})=V_{\omega'} \cap \ker f_{\omega',\omega}.$$
\end{defn}

An important special case in the definition, and the only one which we will
use in the present paper, is that $\omega'$ is obtained from $\omega$ by 
decreasing the degree by $1$ on a single component and increasing it 
correspondingly on an adjacent component.

While we cannot always ensure our linked linear series are simple, we
can ensure they are exact:

\begin{prop}\label{prop:exact-limit} If $(\sL_{\eta},V_{\eta})$ is
defined over $X_{\eta}$ itself, then the resulting linked linear series
is exact.
\end{prop}

The proof is exactly the same as in the two-component case, which is
explained immediately before the statement of Theorem 5.2 of \cite{o-e1}.
Thus, even if $(\sL_{\eta},V_{\eta})$ is not defined over $X_{\eta}$,
we can take a finite base change to make it defined, and blow up the
resulting singularities of the total space to put ourselves into position
to apply Proposition \ref{prop:exact-limit}.

\section{Degenerate linked linear series}\label{sec:degen-link}

The purpose of this section is to analyze the structures of the possible
exact linked linear series lying over limit linear series in the situations
that can arise when $\rho \leq 2$. We will henceforth restrict our attention
to the case that our reducible curve $X_0$ is a chain, although for the
moment we don't have to place any restrictions on the genus of the 
components.

\begin{sit}\label{sit:chain} Suppose that $X_0$ is obtained by starting
with smooth curves $Z_1,\dots,Z_N$, with each $Z_i$ having distinct
marked points $P_i, Q_i$, and gluing $Q_i$ to $P_{i+1}$ for each 
$i=1,\dots,N-1$.
\end{sit}

In this situation, we can encode a multidegree $w$ as follows:

\begin{notn}\label{notn:multidegree} If we have fixed a total degree 
$d$, if we write $w=(c_2,\dots,c_N)$, we let $\md_d(w)$ be the multidegree
which has degree equal to $c_2$ on $Z_1$, to $c_{i+1}-c_i$ on $Z_i$
for $i=2,\dots,N-1$, and to $d-C_N$ on $Z_N$. 

$w$ is \textbf{bounded} if $0 \leq c_i \leq d$ for all $i$.
\end{notn}

To avoid notational clutter, we will frequently write simply $\md(w)$
when the total degree is clear, and we will write abbreviate 
$\sL_{\md(w)}$ by $\sL_w$, $f_{\md(w),\md(w')}$ by $f_{w,w'}$, and so forth.
Note that $\md$ is invertible: in any fixed total degree, any multidegree
$\omega$ has a unique $w$ such that $\omega=\md(w)$.
The total degrees will always be equal to $d$ for the remainder of this
section.

We will assume without further comment that all $w$ are bounded.
The point of this is that if $w$ is bounded, then for all $i$, the map
$f_{w,w^i}$ will be injective on the component $Z_i$ (see the second part 
of Proposition 4.6 of \cite{o-l-t-z}), so we can understand sections in
multidegree $w$ as being glued from the $Z_i$-parts of sections in the
multidegrees $w^i$.

It will also be convenient to use the convention that $c_1=0$ and $c_N=d$ 
always, so that the various conditions we will describe below do not have
to treat endpoints as special cases.

The idea behind Notation \ref{notn:multidegree} is that for $1<i<N$, 
the line bundle 
$\sL_w|_{Z_i}$ is obtained from $\sL^i$ by twisting down by $c_i P_i$ and 
by $(d-c_{i+1})Q_i$, leaving degree $d-c_i-(d-c_{i+1})=c_{i+1}-c_i$.
This notation is very helpful in connection with the way in which we encode
the combinatorial data of a limit linear series. We first describe the
behavior of the maps $f_{w,w'}$ under the above encoding.

\begin{prop}\label{prop:vanishing-comps} Given 
$w=(c_2,\dots,c_N),w'=(c_2',\dots,c_N')$ and total degree $d$, the
map $\sL_{w'} \to \sL_w$ vanishes identically on the component $Z_i$ if and 
only if 
$$\sum_{j=i+1}^N (c'_j-c_j)>
\min_{1 \leq i' \leq N} \sum_{j=i'+1}^N(c'_j-c_j).$$
In particular, if $c'_i<c_i$ or $c'_{i+1}>c_{i+1}$ then the map vanishes 
identically on $Z_i$, and if $c'_i=c_i$ for $i>1$, then the map vanishes 
identically on $Z_i$ if and only if it vanishes identically on $Z_{i-1}$.
\end{prop}

See Proposition 4.6 and Remark 3.14 of \cite{o-l-t-z}.

We now move onto how we encode the discrete data of limit linear series.
First, the following is easy to check via an elimination argument.

\begin{prop}\label{prop:PQ-basis} Let $Z$ be a smooth projective curve,
and $P,Q \in Z$ distinct. Let $(\sL,V)$ be a $\fg^r_d$ on $Z$. 
Then there
is a unique (unordered) set of pairs $(a_0,b_0),\dots,(a_r,b_r)$ 
with all $a_j$ distinct and all $b_j$ distinct such that there
exists a basis $s_0,\dots,s_r$ of $V$ with $\ord_{P} s_j=a_j$ and
$\ord_Q s_j=b_{j}$ for $j=0,\dots,r$. 
\end{prop}

Note that the $s_j$ themselves are not unique, although a given $s_j$
can be modified only by adding multiples of $s_{j'}$ which simultaneously 
satisfy $\ord_P s_{j'}>\ord_P s_j$ and $\ord_Q s_{j'} > \ord_Q s_j$.
We then can introduce a table of numbers to a refined limit linear series as
follows.

\begin{notn}\label{notn:lls-table} Let $(\sL^i,V^i)$ be a refined limit 
$\fg^r_d$ on $X_0$, and for each $i$ let $(a^i_j,b^i_j)_j$ be the set
of pairs given by Proposition \ref{prop:PQ-basis}. 

Construct the 
$(r+1) \times N$ table $T'$ from left to right, with the $i$th column of $T'$ 
consisting of the pairs $(a^i_j,b^i_j)$ for $j=0,\dots,r$, and the ordering
of each column determined as follows: $a^1_j$ should be strictly increasing,
and for $i>1$ and each $j$, we require $a^{i}_j=d-b^{i-1}_j$. For fixed
$i$, we refer to the $a^i_j$ and the $b^i_j$ as making up the 
\textbf{subcolumns} of the $i$th column of $T'$.

For each $j$, let $w_j=(a^2_j,\dots,a^N_j)$, and set
$\omega_j=\md_d(w_j)$. 
\end{notn}

Note that the set of pairs of Proposition \ref{prop:PQ-basis} is giving
a relative ordering of the vanishing sequences at $P$ and $Q$, so the
condition that the limit linear series is refined means that we can
always impose that $a^i_j=d-b^{i-1}_j$. The reason for arranging our table
ordering in this way is that we can always choose sections $s^i_j \in V^i$ 
such that $\ord_{P_i} s^i_j=a^i_j$ and $\ord_{Q_i} s^i_j=b^i_j$, and then
in multidegree $\omega_j$ there is a unique section $s_j$ obtained from gluing 
together the $s^i_j$ (although as noted above, the choices of $s^i_j$ are
not unique in general).

\begin{defn}\label{defn:swaps} We say that a \textbf{swap} occurs in
column $i$ between rows $j,j'$ if $a^i_j<a^i_{j'}$ and $b^i_j<b^i_{j'}$
or if $a^i_j>a^i_{j'}$ and $b^i_{j}>b^i_{j'}$. A swap is \textbf{minimal}
if further $|a^i_j-a^i_{j'}|=|b^i_j-b^i_{j'}|=1$ and either $a^i_j+b^i_j=d$
or $a^i_{j'}+b^i_{j'}=d$.
\end{defn}

Now, given a limit linear series on our chain of curves, there may be
more than one linked linear series lying over it. If the limit linear 
series is ``chain-adaptable'' in the sense of \cite{os20} (i.e., if there are 
no swaps in the table $T'$), the linked linear series is unique, and 
simple, generated by $s_j$ described above. However, in the non-chain-adaptable
case it is not unique. A nonempty open subset of the set of possible linked 
linear series will always be simple, generated by sections similar to the
$s_j$ described above.\footnote{For instance, in the case of a single swap 
they may differ from the $s_j$ by adding multiples of
certain other sections in the first and/or last columns. This results in 
distinct possibilities for the simple linked linear series; see
Example 4.3.5 of \cite{os20}.}
From the point of view of proving Theorem \ref{thm:main}, these simple
cases behave essentially as if they contained all the $s_j$, and are much more
straightforward to handle.
However, even 
among the exact linked linear series, not all of them are necessarily
simple. We can nonetheless use exactness to obtain fairly good control
over what these linked linear series look like. We address all the cases
that can arise for $\rho \leq 2$ below. 

For the rest of this section, we suppose we have fixed a refined limit
linear series along with the resulting table $T'$ as described above, 
as well as a choice of all the $s^i_j$.

The starting point of our analysis is that for any 
$w=(c_2,\dots,c_N)$ (always implicitly assumed bounded), the linkage 
condition implies that the
$(r+1)$-dimensional space $V_w$ in our linked linear series must consist
of sections which are obtained by linear combinations of sections obtained
by gluing, for a fixed $j$, the sections $s^i_j$ to one another as $i$
varies, where each $s^i_j$ that appears must satisfy $a^i_j \geq c_i$ and
$b^i_j \geq d-c_{i+1}$,
and if the first (respectively, second) 
inequality is an equality we must also have $s^{i-1}_j$ (respectively,
$s^{i+1}_j$) included in the gluing. Indeed, a section in $V_w$ must be a 
linear combination of such $s^i_j$, and since the $a^i_j$ and $b^i_j$
are all distinct for fixed $i$, at most one can have equality on each
side, leading to the desired form for the gluing.

\begin{prop}\label{prop:simple} Suppose that 
the $j_0$th row of $T'$ has the property that for all $j<j_0$ we have 
$b^i_j>b^i_{j_0}$ for $i=1,\dots,N-1$, and for all $j>j_0$ we have
$a^i_j>a^i_{j_0}$ for $i=2,\dots,N$. 
Then any linked linear series lying over the given limit linear 
series contains the expected section $s_{j_0}$.
\end{prop}

\begin{proof} We just have to see that the space of global sections in 
multidegree $\omega_{j_0}$ obtained from all possible gluings of the $s^i_j$
has dimension exactly $r+1$, so that any linked linear series must 
contain the whole space, including $s_{j_0}$. But for $j<j_0$ since
$b^i_j>b^i_{j_0}$ for $i<N$, we have $a^{i+1}_j<a^{i+1}_{j_0}$, so
$s^{i+1}_j$ cannot appear at all in multidegree $\omega_{j_0}$. Thus,
only $s^1_j$ can appear, glued to the zero section on every other 
component. Similarly, for $j>j_0$ only $s^N_j$ can appear. And since
each $s^i_{j_0}$ has precisely the desired vanishing at the nodes, 
$s_{j_0}$ is the unique way to glue them together, so we obtain an
$(r+1)$-dimensional space in total, as desired.
\end{proof}

When the hypotheses of Proposition \ref{prop:simple} are not satisfied
for every $j_0$, then we can have linked linear series -- even exact 
ones -- which do not contain all of the $s_{j_0}$, and are not even
simple. Moreover, we expect
that these actually occur as specializations of linear series on the
generic fiber. This leads us to introduce the following notion:

\begin{defn} For $\ell>1$, let $\vec{S}=(S_1,\dots,S_{\ell})$ be a tuple of
subsets of $\{1,\dots,N\}$ such that for all pairs $i<i'$,
every element of $S_i$ is less than or equal to every element of $S_{i'}$,
and such that every element of $\{1,\dots,g\}$ is contained in some
$S_i$.
Let $\vec{j}=(j_1,\dots,j_{\ell})$ be a tuple of elements of 
$\{0,\dots,r\}$, possibly with repetitions.

Then given a fixed limit linear series and corresponding
choices of the $s^i_j$, a \textbf{mixed section} of type $(\vec{S},\vec{j})$
is a $w$ and a section $s$ in multidegree $\md(w)$ which is a 
sum from $i=1$ to $\ell$ of sections obtained by gluing $s^{i'}_{j_i}$
for all $i' \in S_i$ to the zero section on other components.
\end{defn}

In the above definition, it is convenient to allow the possibility that
some of the $S_i$ are empty. The choice of $w$ is not always determined 
uniquely by the type of a mixed section when there are sufficiently large 
gaps between the relevant values of the $a^i_j$, but in our arguments the 
particular value of $w$ will never arise. In cases where the $s^i_j$ are
not uniquely determined, the type of a mixed section may depend on these
choices. 
However, this dependence will be irrelevant to our independence
arguments.

We will show that in the cases of interest, even when a given $s_j$
is not in our linked linear series, we can ensure that there are mixed
sections of rather precise forms, which can in some sense take the place
of the missing $s_j$. It will be convenient to carry out these 
constructions in two steps.

The following single swap between a pair of rows is the only form of 
degeneracy that can occur in the $\rho=1$ case; in the below lemma, we
also allow for the possibility that there could be other swaps occurring
in other rows.

\begin{prop}\label{prop:single-swap-1st}
Suppose that our limit linear series has a single swap between the 
$j_0$th and $(j_0-1)$st rows, occurring in the $i_0$th column,
and for all $j<j_0-1$ we have 
$b^i_j>b^i_{j_0-1},b^i_{j_0}$ for $i=1,\dots,N-1$, and for all $j>j_0$ 
we have $a^i_j>a^i_{j_0},a^i_{j_0-1}$ for $i=2,\dots,N$. 

Then any linked linear series lying over the given limit linear 
series contains the expected section $s_{j_0-1}$, and the multidegrees
associated to
$(a^2_{j_0-1},\dots,a^{i_0}_{j_0-1},a^{i_0+1}_{j_0},\dots,a^N_{j_0})$
and 
$(a^2_{j_0},\dots,a^{i_0}_{j_0},a^{i_0+1}_{j_0-1},\dots,a^N_{j_0-1})$
must contain the respective images of the section $s_{j_0}$. These
images consist respectively of $0$ on the first $i_0-1$ components and
$s^i_{j_0}$ for $i=i_0,\dots,N$, and of $s^i_{j_0}$ for $i=1,\dots,i_0$,
and $0$ on the last $N-i_0$ components.

Let $w=(c_2,\dots,c_N)$. If 
$c_i<a^i_{j_0-1},a^i_{j_0}$ for all $i$, the linked linear series contains
$s^1_{j_0}$ in multidegree $\md(w)$, and if $c_i>a^i_{j_0-1},a^i_{j_0}$ for 
all $i$, the linked linear series contains $s^N_{j_0}$ in multidegree $\md(w)$
(in both cases, glued to $0$ on the other components).
\end{prop}

\begin{proof} 
First, in the multidegree $\omega_{j_0-1}$, as in the proof of Proposition
\ref{prop:simple}, the $s^i_j$ for $j \neq j_0-1,j_0$
can only contribute for $i=1$ (if $j<j_0-1$) or $i=N$ (if $j>j_0$),
and the $s^i_{j_0-1}$ glue uniquely to give $s_{j_0-1}$. Finally, the
$s^i_{j_0}$ can only contribute at $i=i_0$, so we find that the space
obtained from all the $s^i_j$ is $(r+1)$-dimensional, and $s_{j_0-1}$
must be in the linked linear series, as desired.

Next, consider 
$w'= (a^2_{j_0-1},\dots,a^{i_0}_{j_0-1},a^{i_0+1}_{j_0},\dots,a^N_{j_0})$.
Note that $f_{w_{j_0},w'} (s_{j_0})$
is equal to $s_{j_0}$ from $i_0$ to $N$ (inclusive), and $0$ strictly before 
$i_0$. We claim that the 
space of possible sections from the $s^i_j$ in multidegree $\md(w')$ is
precisely $(r+1)$-dimensional, so the linked linear series is uniquely
determined in this multidegree. By hypothesis, the $s^i_j$ for $j<j_0-1$ can 
only contribute for $i=1$, and the $s^i_j$ for $j>j_0$ can only 
contribute for $i=N$. The $s^i_{j_0-1}$ could in principle contribute for
$i<i_0$ and $i=N$, but if the $s^i_{j_0-1}$ appeared for $i<i_0$, they
all would be nonvanishing at the relevant nodes, and they would have to
glue to something nonvanishing in the $i_0$th column. But this would have
to be $s^{i_0}_{j_0-1}$, which does not have enough vanishing on the right
to appear in multidegree $\md(w')$. Thus, we conclude that the $s^i_{j_0-1}$
can only appear for $i=N$ (where it is glued to the zero section on all
other columns). Finally, the $s^i_{j_0}$ can only appear for $i \geq i_0$,
where they are nonzero at all interior nodes, and therefore have a unique 
gluing, which must yield $f_{w_{j_0},w'} (s_{j_0})$. Thus we get the claimed
dimension $r+1$, and conclude that
$f_{w_{j_0},w'} (s_{j_0})$ is contained in the linked linear series.

Similarly, if 
$w''= (a^2_{j_0},\dots,a^{i_0}_{j_0},a^{i_0+1}_{j_0-1},\dots,a^N_{j_0-1})$,
we find that space of possible sections is $(r+1)$-dimensional, 
and contains $f_{w_{j_0},w''} (s_{j_0})$.

Now, suppose we are given $w$ with $c_i<a^i_{j_0-1},a^i_{j_0}$ for all
$i$. Then Proposition \ref{prop:vanishing-comps} implies that $f_{w'',w}$
is nonzero precisely on
the $1$st component, so $f_{w'', w} (f_{w_{j_0},w''} (s_{j_0}))$ is equal
to $s^1_{j_0}$ glued to $0$, as desired. The situation with
$c_i>a^i_{j_0-1},a^i_{j_0}$ is similar, but with $w'$ in place of $w''$.
\end{proof}

Now we will systematically consider cases where the limit linear series
has only one or two swaps, with no swaps in any other rows. We start with
the case of a single swap.

\begin{prop}\label{prop:single-swap}
Suppose that our limit linear series has precisely one swap, between the 
$j_0$th and $(j_0-1)$st rows, and occurring in the $i_0$th column.

Then any linked linear series lying over the given limit linear 
series contains the expected sections $s_{j}$ for all $j \neq j_0$.
If the linked linear series is exact, then 
it must contain mixed sections 
$s_{j_0}'$ and $s_{j_0}''$ of type $((S_1',S_2'),(j_0-1,j_0))$ and
$((S_1'',S_2''),(j_0,j_0-1))$ respectively, with 
$S_1'$ supported strictly left of $i_0$, and $S_2''$ supported strictly
right of $i_0$.
\end{prop}

Note that the possibility that the linked linear series contains the section
$s_{j_0}$ itself is contained in the proposition by allowing $S_1'$ and
$S_2''$ both to be empty.

\begin{proof} Start with the $w'$ from the proof of Proposition
\ref{prop:single-swap-1st}. Note that if $i_0=1$, then the proposition says
that $s_{j_0}$ itself is in our linked linear series, consistent with the
stated form for $s_{j_0}'$. Otherwise, if $i_0>1$
we consider iteratively changing $w'$ by increasing the twists by $1$ 
for $i'\leq i_0$ (starting at $i_0$) until they each agree with 
$a^{i'}_{j_0}$. We note that
every such modified $w'$ has an $(r+2)$-dimensional space of global sections
obtained from the $s^i_j$, described explicitly as follows: $s^1_j$ for 
$j<j_0-1$;
$s^N_j$ for $j>j_0$; $s^N_{j_0-1}$; a section obtained by gluing 
the $s^i_{j_0-1}$ for $i$ from $1$ to $i'-1$ (which is the last column in 
which $w'$ agrees with $a^i_{j_0-1}$); and a section 
obtained by gluing the $s^i_{j_0}$ from either $i'-1$ or $i'$ to $N$,
beginning with the last column in which $w'$ has coefficient 
strictly less than $a^i_{j_0}$. For each 
$j \neq j_0-1$, since there is a unique section constructed from the 
$s^i_j$, it is necessarily equal to $f_{w_j,w'} (s_j)$. In addition, since
we know $s_j$ is in our linked linear series for $j\neq j_0$, we
have that $f_{w_j,w'} (s_j)$ is necessarily contained in our linked linear
series for $j \neq j_0-1,j_0$.

Now, suppose that our linked linear series contained 
$f_{w_{j_0},w'} (s_{j_0})$ for the old $w'$; we claim that it either also
contains it for the new $w'$, or contains a section of the form desired 
for $s_{j_0}'$. Indeed, increasing the twist in the $i$th column 
corresponding to twisting once by every component from $i$ to $N$. We 
observe that $f_{w_{j_0},w'} (s_{j_0})$ is in the kernel of the map from the 
old $w'$ to the new one, so by the definition of exactness, the linked 
linear series must contain some $s$ in the new multidegree mapping to 
$f_{w_{j_0},w'} (s_{j_0})$ in the old one. Using the above description of the
space of global sections, this is necessarily a combination of the
$f_{w_j,w'} (s_j)$ for $j<j_0-1$ and $j=j_0$, together with the section from 
the $s^i_{j_0-1}$ for $i=1$ to $i'-1$. Moreover, since we observed above
that $f_{w_j,w'} (s_j)$ is contained in our linked linear series for $j<j_0-1$,
we can subtract these off to obtain a combination of the sections from the
$j_0-1$ and $j_0$ rows. If the $j_0-1$ term vanishes, we have that
$f_{w_{j_0},w'} (s_{j_0})$ is contained in our linked linear series for
the new $w'$, and if the $j_0-1$ term is nonzero, we have something of the 
desired form
for $s_{j_0}'$ (with the minimal element of $S_2'$ being either $i'$ or 
$i'-1$ according to where 
$f_{w_{j_0},w'} (s_{j_0})$ begins), as claimed.
Iterating this process, we either obtain the desired $s_{j_0}'$, or we
eventually reach $w'=w_{j_0}$ and find that the linked linear series
actually contains $s_{j_0}$ itself. 

As the situation is completely symmetric, the construction of $s_{j_0}''$ is
similar, starting from the multidegree $w''$ from the proof of 
Proposition \ref{prop:single-swap-1st}.
\end{proof}

\begin{ex}\label{ex:g22-grd} We will use the below as a running example,
showing a possible table $T'$
associated to a limit linear series in the case $r=6$, $g=22$, $d=25$ (here
we assume every component is of genus $1$).
\begin{center}
\resizebox{\textwidth}{!}{
\begin{tabular}{lr|lr|lr|lr|lr|lr|lr|lr|lr|lr|lr|lr|lr|lr|lr|lr|lr|lr|lr|lr|lr|lr}
$0$ & $25$ & $0$ & $24$ & $1$ & $23$ & $2$ & $22$ & $3$ & $21$ & $4$ & $20$ & $5 $ & $19$ & $6$ & $19$ & $6$ & $18$ & $7$ & $17$ & $8$ & $16$ & $9$ & $15$ & $10$ & $14$ & $11$ & $13$ & $12$ & $12$ & $13$ & $12$ & $13$ & $11$ & $14$ & $10$ & $15$ & $9$ & $16$ & $8$ & $17$ & $7$ & $18$ & $6$ \\
$1$ & $23$ & $2$ & $23$ & $2$ & $22$ & $3$ & $21$ & $4$ & $20$ & $5$ & $19$ & $6 $ & $18$ & $7$ & $17$ & $8$ & $16$ & $9$ & $16$ & $9$ & $15$ & $10$ & $14$ & $11 $ & $13$ & $12$ & $12$ & $13$ & $11$ & $14$ & $10$ & $15$ & $10$ & $15$ & $9$ & $16$ & $8$ & $17$ & $7$ & $18$ & $6$ & $19$ & $5$ \\
$2$ & $22$ & $3$ & $21$ & $4$ & $21$ & $4$ & $20$ & $5$ & $19$ & $6$ & $18$ & $7 $ & $17$ & $8$ & $16$ & $9$ & $14$ & $11$ & $13$ & $12$ & $12$ & $13$ & $12$ & $ 13$ & $11$ & $14$ & $10$ & $15$ & $9$ & $16$ & $8$ & $17$ & $7$ & $18$ & $6$ & $ 19$ & $6$ & $19$ & $5$ & $20$ & $4$ & $21$ & $3$ \\
$3$ & $21$ & $4$ & $20$ & $5$ & $19$ & $6$ & $19$ & $6$ & $18$ & $7$ & $17$ & $8 $ & $16$ & $9$ & $15$ & $10$ & $15$ & $10$ & $14$ & $11$ & $14$ & $11$ & $13$ & $12$ & $12$ & $13$ & $11$ & $14$ & $10$ & $15$ & $9$ & $16$ & $8$ & $17$ & $8$ & $17$ & $7$ & $18$ & $6$ & $19$ & $5$ & $20$ & $4$ \\
$4$ & $20$ & $5$ & $19$ & $6$ & $18$ & $7$ & $17$ & $8$ & $17$ & $8$ & $16$ & $9 $ & $15$ & $10$ & $14$ & $11$ & $13$ & $12$ & $12$ & $13$ & $11$ & $14$ & $10$ & $15$ & $10$ & $15$ & $9$ & $16$ & $8$ & $17$ & $7$ & $18$ & $6$ & $19$ & $5$ & $20$ & $4$ & $21$ & $4$ & $21$ & $3$ & $22$ & $2$ \\
$5$ & $19$ & $6$ & $18$ & $7$ & $17$ & $8$ & $16$ & $9$ & $15$ & $10$ & $15$ & $ 10$ & $14$ & $11$ & $13$ & $12$ & $12$ & $13$ & $11$ & $14$ & $10$ & $15$ & $9$ & $16$ & $8$ & $17$ & $8$ & $17$ & $7$ & $18$ & $6$ & $19$ & $5$ & $20$ & $4$ & $21$ & $3$ & $22$ & $2$ & $23$ & $2$ & $23$ & $1$ \\
$6$ & $18$ & $7$ & $17$ & $8$ & $16$ & $9$ & $15$ & $10$ & $14$ & $11$ & $13$ & $12$ & $13$ & $12$ & $12$ & $13$ & $11$ & $14$ & $10$ & $15$ & $9$ & $16$ & $8$ & $17$ & $7$ & $18$ & $6$ & $19$ & $6$ & $19$ & $5$ & $20$ & $4$ & $21$ & $3$ & $22$ & $2$ & $23$ & $1$ & $24$ & $0$ & $25$ & $0$ \\
\end{tabular}
}
\end{center}

Since there is no ramification at $P_1$, the first entries of the table
agree with the row labels, so we have not shown the labels separately.

Note that we have a single swap, occurring in the $9$th column between 
the $j=2$ and $j=3$ rows. This leads to having an extra dimension of
possibilities in the multidegree obtained from the $j=3$ row, as the $j=2$
row can appear either in the first or last columns. Consequently, it is
possible that an exact linked linear series lying over this limit linear
series might not contain $s_3$, but might only contain mixed sections
$s_3'$ and $s_3''$ as in Proposition \ref{prop:single-swap}, with 
$s_3'$ agreeing with $s_3$ for $i \geq 9$, but switching to $s_2$
at some $i<9$, and $s_3''$ agreeing with $s_3$ for $i \leq 9$, but 
switching to $s_2$ at some $i>9$. In both cases, the switch occurs in a column
mixing $s^i_2$ and $s^i_3$ unless, the column in question has a gap of at
least $2$ between the $j=2$ and $j=3$ rows. Since this doesn't occur for 
$i <9$, we see that $s_3'$ always has a mixed column, while $s_3''$ may
not. 
\end{ex}

When $\rho=2$, there are four additional forms of degeneracy that 
can occur, which we consider one by one. They all involve having
exactly two swaps, occurring in distinct columns. The first case is
when the swaps occur in disjoint pairs of rows.

\begin{prop}\label{prop:disjoint-swap} (``Disjoint swap'')
Suppose that our limit linear series
contains precisely two swaps, and these occur in disjoint pairs of rows, 
say $j_0-1,j_0$ and $j_1-1,j_1$. 
Then any linked linear series lying over the given limit linear 
series contains the expected sections $s_{j}$ for all $j \neq j_0, j_1$.
If the linked linear series is exact, then for $\ell=0,1$ 
it must contain mixed sections 
$s_{j_{\ell}}'$ and $s_{j_{\ell}}''$ of type 
$((S_{1+2\ell}',S_{2+2\ell}'),(j_{\ell}-1,j_{\ell}))$ and
$((S_{1+2\ell}'',S_{2+2\ell}''),(j_{\ell},j_{\ell}-1))$ respectively, with 
$S_{1+2\ell}'$ supported strictly left of $i_{\ell}$ and 
$S_{2+2\ell}''$ supported strictly right of $i_{\ell}$.
\end{prop}

\begin{proof} This is essentially identical to the proof of Proposition
\ref{prop:single-swap}. The only new point which arises is that in
constructing the sections $s_{j_0}',s_{j_0}''$, we need to know that
we can always subtract off any $s_{j_1}$ part which arises in the
iterative procedure, and similarly with $j_0$ and $j_1$ switched. But
this follows from the last assertion of Proposition 
\ref{prop:single-swap-1st}.
\end{proof}

The next case is that
a single pair of rows can undergo two swaps in different columns.

\begin{table}\label{tab:repeat-swap}
\begin{center}
\resizebox{\textwidth}{!}{
\begin{tabular}{cc|lr|lr|lc|lr|lr|lc}
\multicolumn{4}{c}{} & \multicolumn{2}{c}{$i_0$} & \multicolumn{4}{c}{} 
& \multicolumn{2}{c}{$i_1$} & \multicolumn{2}{c}{} \\
\vdots & & & & & & & & & & & & & \\
$(j_0-1)$ & & $a-1$ & $d-a$ & $a$ & $d-a-2$ & $a+2$ & \dots
& $a'$ & $d-a'-1$ & $a'+1$ & $d-a'-1$ & $a'+1$ & \\
$(j_0)$ & \dots & $a$ & $d-a-1$ & $a+1$ & $d-a-1$ & $a+1$ & \dots
& $a'-1$ & $d-a'$ & $a'$ & $d-a'-2$ & $a'+2$ & \dots \\
\vdots & & & & & & & & & & & & & \\
\end{tabular}
}
\end{center}
\caption{A typical example of the ``repeated swap'' situation. In general,
there may be larger gaps between the rows, although when $\rho=2$ the
gaps at the $i_0$ and $i_1$ columns must both be equal to $1$.}
\end{table}

\begin{prop}\label{prop:repeat-swap} (``Repeated swap'')
Suppose that our limit linear series has precisely two swaps, both
between the 
$j_0$th and $(j_0-1)$st rows, with the first occurring in the $i_0$th column,
and the second in the $i_1$st column for some $i_1>i_0$. 

Then any exact linked linear series lying over the given limit linear 
series 
contains mixed sections 
$s_{j_0-1}'$ and $s_{j_0-1}''$ of type $((S_1',S_2',S_3'),(j_0-1,j_0,j_0-1))$
and $((S_1'',S_2''),(j_0-1,j_0)$ respectively, with 
$S_2'$ supported strictly left of $i_1$ 
and $S_2''$ supported strictly right of $i_1$,
and 
it contains mixed sections 
$s_{j_0}'$ and $s_{j_0}''$ of type $((S_4',S_5'),(j_0-1,j_0))$
and $((S_3'',S_4'',S_5''),(j_0,j_0-1,j_0)$ respectively, with 
$S_4'$ supported strictly left of $i_0$ 
and 
$S_4''$ supported strictly right of $i_0$.
\end{prop}

\begin{proof} The proof is similar to the proof of Proposition 
\ref{prop:single-swap}. For $s_{j_0-1}'$, we first consider 
$w'=(a^2_{j_0-1},\dots,a^{i_0}_{j_0-1},a^{i_0+1}_{j_0},\dots,a^{i_1}_{j_0},
a^{i_1+1}_{j_0-1},\dots,a^N_{j_0-1})$. 
Note that $f_{w_{j_0-1},w'} (s_{j_0-1})$
is equal to $s_{j_0-1}$ from $i_1$ to $N$ (inclusive), and $0$ elsewhere.
Indeed, these are the only columns in which the $s^i_{j_0-1}$ can be
supported, since they do not satisfy the correct inequalities from $i_0$
to $i_1-1$, and for $i<i_0$ they satisfy them with equality, so would have
to be glued to a nonzero element in the $i_0$th column.
As in the proof of Proposition \ref{prop:simple}, we check that we have
dimension exactly $r+1$ in multidegree $\md(w')$, with the unique contribution
from the $j_0$ row coming from $s^N_{j_0}$.
Thus, we find that $f_{w_{j_0-1},w'} (s_{j_0-1})$ is necessarily contained
in multidegree $\md(w')$. 

We then iterate changing $w'$ by $1$, increasing the
twist by $1$ in the $i'$th column for $i' \leq i_1$ to change them from 
$a^{i'}_{j_0}$ to $a^{i'}_{j_0-1}$.
Using exactness, at each stage we either 
find the linked linear series still contains $f_{w_{j_0-1},w'} (s_{j_0-1})$ 
for the new value of $w'$, or it contains the sum of 
$f_{w_{j_0-1},w'} (s_{j_0-1})$ with a section obtained by gluing the 
$s^{i}_{j_0}$ for $i=i_0,\dots,i'-1$. In the first case, we continue to 
iterate the process of changing $w'$, and if we do not ever get the second 
case, we end up with $s_{j_0-1}$ itself in our linked linear series. On
the other hand, once the second case occurs, we begin to iteratively 
change $w'$ by increasing the twist by $1$ in the $i'$th column for
$i' \leq i_0$ to change them from $a^{i'}_{j_0-1}$ to $a^{i'}_{j_0}$. 
Each time the twist increases above $a^{i'}_{j_0-1}$, we could obtain
a contribution obtained from gluing $s^i_{j_0-1}$ from $i=1$ to $i'-1$,
and if this occurs, we get our desired $s_{j_0-1}'$. Otherwise, we
keep iterating, and each time the twist at $i'$ reaches $a^{i'}_{j_0}$,
the portion of the section obtained from the $s^i_{j_0}$ extends to
include $i'-1$. Again, if we never get a contribution from the 
$s^i_{j_0-1}$ for $i \leq i'$, we will end up with a section as required
for $s_{j_0-1}'$, having $S_1'=\emptyset$. 

The construction of $s_{j_0-1}''$ is similar, but simpler: we set our 
initial 
$w''=(a^2_{j_0-1},\dots,a^{i_1}_{j_0-1},a^{i_1+1}_{j_0},\dots,a^N_{j_0})$, 
and then we 
iteratively decrease the twists
for $i'>i_1$ by $1$ to change them from $a^{i'}_{j_0}$ to $a^{i'}_{j_0-1}$,
until we obtain the desired result.

The construction of $s_{j_0}'$ and $s_{j_0}''$ follows the same process.
For $s_{j_0}'$, we start with 
$w'=(a^2_{j_0-1},\dots,a^{i_0}_{j_0-1},a^{i_0+1}_{j_0},\dots,a^N_{j_0})$,
and we iteratively increase the twists
for $i' \leq i_0$ by $1$ to change them from $a^{i'}_{j_0-1}$ to 
$a^{i'}_{j_0}$.
Finally, for $s_{j_0}''$, we start with 
$w''=(a^2_{j_0},\dots,a^{i_0}_{j_0},a^{i_0+1}_{j_0-1},\dots,a^{i_1}_{j_0-1},
a^{i_1+1}_{j_0},\dots,a^N_{j_0})$,
obtaining a section glued from
the $s^i_{j_0}$ for $i \leq i_0$. We iteratively decrease the twists
for $i' > i_0$ by $1$ to change them from $a^{i'}_{j_0-1}$ to 
$a^{i'}_{j_0}$, until we obtain a contribution from the $a^i_{j_0-1}$
(necessarily ending at $i_1$), and then we iteratively decrease the twists
for $i' > i_1$ by $1$ to change them from $a^{i'}_{j_0}$ to 
$a^{i'}_{j_0-1}$, eventually obtaining either $s_{j_0}$ itself, or the
desired $s_{j_0}''$.
\end{proof}

The last cases involve three consecutive rows undergoing two swaps. There
are only two different ways this can occur for $\rho=2$, but it turns
out that these two ways behave quite differently. The two cases can
be understood as either having one row which sums to $d$ in both of the
relevant columns, or one row which sums to $d-2$ in both of the relevant
columns. The latter turns out to be more degenerate.

\begin{table}\label{tab:3-cycle-1}
\begin{center}
\resizebox{\textwidth}{!}{
\begin{tabular}{cc|lr|lr|lc|lr|lr|lc}
\multicolumn{4}{c}{} & \multicolumn{2}{c}{$i_0$} & \multicolumn{4}{c}{} 
& \multicolumn{2}{c}{$i_1$} & \multicolumn{2}{c}{} \\
\vdots & & & & & & & & & & & & & \\
$(j_0-1)$ & & $a-1$ & $d-a$ & $a$ & $d-a-1$ & $a+1$ & 
& $a'-1$ & $d-a'$ & $a'$ & $d-a'-2$ & $a'+2$ & \\
$(j_0)$ & \dots & $a$ & $d-a-1$ & $a+1$ & $d-a-3$ & $a+3$ & \dots
& $a'+1$ & $d-a'-2$ & $a'+2$ & $d-a'-3$ & $a'+3$ & \dots \\
$(j_0+1)$ & & $a+1$ & $d-a-2$ & $a+2$ & $d-a-2$ & $a+2$ &
& $a'$ & $d-a'-1$ & $a'+1$ & $d-a'-1$ & $a'+1$ & \\
\vdots & & & & & & & & & & & & & \\
\end{tabular}
}
\end{center}
\caption{A typical example of the ``first $3$-cycle'' situation. In general,
there may be larger gaps between the rows, although when $\rho=2$ the
gaps between the $j_0$ and $j_0+1$ rows at $i_0$ and between the $j_0-1$ and
$j_0+1$ rows at $i_1$ must both be equal to $1$.}
\end{table}

\begin{prop}\label{prop:3-cycle-1} (``First $3$-cycle'')
Suppose that our limit linear series has one swap between the 
$j_0$th and $(j_0+1)$st rows occurring in the $i_0$th column,
and a second swap between the $(j_0-1)$st and $(j_0+1)$st rows in the 
$i_1$st column for some $i_1>i_0$, and no other swaps.

Then any linked linear series lying over the given limit linear series
contains $s_{j_0-1}$ and $s_{j_0}$. If further the linked linear series is
exact, then it 
contains mixed sections 
$s_{j_0+1}'$, $s_{j_0+1}''$ and $s_{j_0+1}''$
of type $((S_1',S_2',S_3'),(j_0-1,j_0,j_0+1))$, 
$((S_1'',S_2'',S_3''),(j_0+1,j_0-1,j_0))$ 
and $((S_1''',S_2''',S_3'''),(j_0,j_0+1,j_0-1))$, 
respectively, with 
$S_1'$ supported strictly left of $i_1$, $S_2'$ supported strictly left of
$i_0$, $S_2''$ supported strictly right of $i_1$, $S_3''$ supported strictly
right of $i_0$, $S_1'''$ supported strictly left of $i_0$, and $S_3'''$
supported strictly right of $i_1$.
\end{prop}

Note that if $S_2'=\emptyset$, then $S_1'$ may contain elements greater
than $i_0$, and similarly if $S_2''=\emptyset$, then $S_3''$ may contain
elements less than $i_1$.

\begin{proof} First, it is routine to check that the multidegrees 
$\omega_{j_0-1}$
and $\omega_{j_0}$ both have only $(r+1)$-dimensional spaces of
possible sections, so that $s_{j_0-1}$ and $s_{j_0}$ must both lie in any
linked linear series. Indeed, for the former, the $s^i_{j_0}$ can contribute
only for $i=N$, while the $s^i_{j_0+1}$ can contribute only for $i=i_1$,
while for the latter, the $s^i_{j_0-1}$ can contribute only for $i=1$, and
the $s^i_{j_0+1}$ can contribute only for $i=i_0$.

Now, to contruct the sections $s_{j_0+1}'$, $s_{j_0+1}''$ and $s_{j_0+1}'''$
we proceed as in the previous propositions. For $s_{j_0+1}'$, we start with
$w'=(a^2_{j_0-1},\dots,a^{i_1}_{j_0-1},a^{i_1+1}_{j_0+1},\dots,a^N_{j_0+1})$,
and then iteratively increase the twist by $1$ at a time for $i' \leq i_1$,
initially increasing it from $a^{i'}_{j_0-1}$ to $a^{i'}_{j_0+1}$. For
$i'>i_0$, this process behaves as before, either extending the contribution
from the $a^i_{j_0+1}$ iteratively to the left without introducing any
other nonzero contributions, or producing a section $s_{j_0+1}'$ as
desired, having $S_2=\emptyset$. Once $i' \leq i_0$, we still iteratively
increase the twist from $a^{i'}_{j_0-1}$ to $a^{i'}_{j_0+1}$, but we
are required to pass $a^{i'}_{j_0}$ in the process. This introduces a 
third possibility: once the twist at $i'$ is strictly greater than
$a^{i'}_{j_0}$, we could obtain a contribution from $s^{i'-1}_{j_0}$.
Also, for $i'<i_0$, once the twist at $i'$ is equal to $a^{i'}_{j_0}$,
we could obtain a contribution from both $s^{i'-1}_{j_0}$ and $s^{i'}_{j_0}$.
If either of these occurs, we move to the next $i'$, and for the remaining $i'$,
instead of increasing the twist from $a^{i'}_{j_0-1}$ to $a^{i'}_{j_0+1}$,
we only increase to $a^{i'}_{j_0}$. Note that we may obtain contributions from
the $s^i_{j_0}$ (for $i=i'-1$ or $i=i'-1,i'$) and $s^i_{j_0-1}$ (for 
$i=1,\dots,i'-1$) simultaneously at some point, which still gives an
$s_{j_0+1}'$ of the desired form. On the other hand, if we never obtain
a contribution from the $s^i_{j_0}$, then the resulting $s_{j_0+1}'$ 
simply has $S_2'=\emptyset$.

For $s_{j_0+1}''$, we start with
$w''=(a^2_{j_0+1},\dots,a^{i_0}_{j_0+1},a^{i_0+1}_{j_0},\dots,a^N_{j_0})$,
and then follow the same procedure 
as for $s_{j_0+1}'$, iteratively decreasing the twist at $i'>i_0$ from
$a^{i'}_{j_0}$ to $a^{i'}_{j_0+1}$, with the possibility of a contribution
from the $s^i_{j_0-1}$ once $i'$ passes $i_1$.

Finally, for $s_{j_0+1}'''$ set
$w'''=(a^2_{j_0},\dots,a^{i_0}_{j_0},a^{i_0+1}_{j_0+1},\dots,a^{i_1}_{j_0+1},
a^{i_1+1}_{j_0-1},\dots,a^N_{j_0-1})$ initially.
We then iteratively increase the twist at $i' \leq i_0$ from $a^{i'}_{j_0}$
to $a^{i'}_{j_0+1}$, and iteratively decrease the twist at $i'>i_1$ from
$a^{i'}_{j_0-1}$ to $a^{i'}_{j_0+1}$ to construct $s_{j_0+1}'''$.
\end{proof}

\begin{table}\label{tab:3-cycle-2}
\begin{center}
\resizebox{\textwidth}{!}{
\begin{tabular}{cc|lr|lr|lc|lr|lr|lc}
\multicolumn{4}{c}{} & \multicolumn{2}{c}{$i_0$} & \multicolumn{4}{c}{} 
& \multicolumn{2}{c}{$i_1$} & \multicolumn{2}{c}{} \\
\vdots & & & & & & & & & & & & & \\
$(j_0-1)$ & & $a-1$ & $d-a$ & $a$ & $d-a-2$ & $a+2$ & 
& $a'$ & $d-a'-1$ & $a'+1$ & $d-a'-3$ & $a'+3$ & \\
$(j_0)$ & \dots & $a$ & $d-a-1$ & $a+1$ & $d-a-1$ & $a+1$ & \dots
& $a'-1$ & $d-a'$ & $a'$ & $d-a'-1$ & $a'+1$ & \dots \\
$(j_0+1)$ & & $a+1$ & $d-a-2$ & $a+2$ & $d-a-3$ & $a+3$ &
& $a'+1$ & $d-a'-2$ & $a'+2$ & $d-a'-2$ & $a'+2$ & \\
\vdots & & & & & & & & & & & & & \\
\end{tabular}
}
\end{center}
\caption{A typical example of the ``second $3$-cycle'' situation. In general,
there may be larger gaps between the rows, although when $\rho=2$ the
gaps between the $j_0-1$ and $j_0$ rows at $i_0$ and between the $j_0-1$ and
$j_0+1$ rows at $i_1$ must both be equal to $1$.}
\end{table}

\begin{prop}\label{prop:3-cycle-2} (``Second $3$-cycle'')
Suppose that our limit linear series has one swap between the 
$(j_0-1)$st and $j_0$th rows occurring in the $i_0$th column,
and a second swap between the $(j_0-1)$st and $(j_0+1)$st rows in the 
$i_1$st column for some $i_1>i_0$, and no other swaps.

Then any linked linear series lying over the given limit linear series
contains $s_{j_0-1}$. If further the linked linear series is
exact, then it 
contains mixed sections 
$s_{j_0}'$ and $s_{j_0}''$ of type $((S_1',S_2'),(j_0-1,j_0))$
and $((S_1'',S_2'',S_3''),(j_0,j_0+1,j_0-1)$ respectively, with 
$S_1'$ supported strictly left of $i_0$, $S_2''$ supported at or right of
$i_1$, and $S_3''$ supported strictly right of $i_0$.
Similarly, it 
contains mixed sections 
$s_{j_0+1}'$ and $s_{j_0+1}''$ of type $((S_3',S_4',S_5'),(j_0-1,j_0,j_0+1))$
and $((S_4'',S_5''),(j_0+1,j_0-1)$ respectively, with 
$S_3'$ supported strictly left of $i_1$, $S_4'$ supported at or left of 
$i_0$, and $S_5''$ supported strictly right of $i_1$. Moreover, if
$i_1 \in S_2''$ then also $i_1 \in S_1''$, and if $i_0 \in S_4'$, then also
$i_0 \in S_5'$.
Finally, either we can have $S_2'=S_4''=\{1,\dots,N\}$, or
it also contains a mixed section $s'''$ of type 
$((S_1''',S_2''',S_3'''),(j_0,j_0-1,j_0+1))$, where every element of $S_2'''$
is strictly between $i_0$ and $i_1$.
\end{prop}

\begin{proof} For the most part, this is straightforward and similar
to the previous propositions, 
but there is one new subtlety to
address, and the idea for the construction of $s'''$ is new. We first 
construct $s_{j_0}'$, starting with 
$w'=(a^2_{j_0-1},\dots,a^{i_0}_{j_0-1},a^{i_0+1}_{j_0},\dots,a^N_{j_0})$.
We then iteratively increase the twist
from $a^{i'}_{j_0-1}$ to $a^{i'}_{j_0}$ for $i' \leq i_0$, and
obtain our $s_{j_0}'$ as usual. We then do the same procedure for 
$s_{j_0+1}''$, starting with
$w''=(a^2_{j_0+1},\dots,a^{i_1}_{j_0+1},a^{i_1+1}_{j_0-1},\dots,a^N_{j_0-1})$.

Next, we construct $s_{j_0}''$, starting with
$w''=(a^2_{j_0},\dots,a^{i_0}_{j_0},a^{i_0+1}_{j_0-1},\dots,a^N_{j_0-1})$.
We then iteratively decrease the twist at $i' >i_0$ from
$a^{i'}_{j_0-1}$ to $a^{i'}_{j_0}$. For $i' \leq i_1$, this behaves as
in the previous propositions, with one new subtlety: for each intermediate
value of $w'$, the $s^i_{j_0+1}$ can contribute only in the $i_1$
column, but because we do not know that $s_{j_0+1}$ is contained in our
linked linear series, we also do not know \emph{a priori} that this 
contribution from $s^{i_1}_{j_0+1}$ in multidegree $\md(w')$ is contained in 
our linked linear series. However, since we have already constructed 
$s_{j_0+1}''$, we can use its image in $\md(w')$. One checks that its only
possible support in $\md(w')$ is in the $i_1$ column, 
so that in fact the multidegree-$\md(w')$ part of our linked linear series 
necessarily contains the section given by $s^{i_1}_{j_0+1}$, and we
can subtract it off as necessary from the section we are constructing.
Thus, for $i' \leq i_1$, we can iterate as before, and will either
obtain an $s_{j_0}''$ as desired (with $S_2''=\emptyset$), or we will
obtain a section made up of the $s^i_{j_0}$ for $i \leq i_1$, and
vanishing identically for $i>i_1$. In the latter case, we continue to
iteratively decrease the twists defining $w'$ for $i>i_1$, but as in
the construction of $s_{j_0+1}'$ in the proof of Proposition 
\ref{prop:3-cycle-1}, to get from $a^{i'}_{j_0-1}$ to $a^{i'}_{j_0}$ we
need to pass $a^{i'}_{j_0+1}$, which is where the possible contribution
from the $j_0+1$ may occur.

The construction of $s_{j_0+1}'$ follows the same pattern as
that of $s_{j_0}''$, but starting with
$w'=(a^2_{j_0-1},\dots,a^{i_1}_{j_0-1},a^{i_1+1}_{j_0+1},\dots,a^N_{j_0+1})$.
Here we use the image of $s_{j_0}'$ in order to subtract 
off any contributions of $s^{i_0}_{j_0}$ which occur.

Finally, for $s'''$, we start with $w'=w_{j_0}$.
We observe that there is an $(r+2)$-dimensional space of potential sections
in multidegree $\omega_{j_0}$, with the $s^i_j$ for $j<j_0-1$ contributing 
only for
$i=1$, the $s^i_j$ for $j \geq j_0+1$ contributing only for $i=N$, the
$s^i_{j_0}$ contributing only with $s_{j_0}$ itself, and the $s^i_{j_0-1}$
contributing separately for $i=1$ and $i=N$. 
We must therefore have a three-dimensional space of combinations of the
four sections $s^1_{j_0-1}$, $s^N_{j_0-1}$, $s^N_{j_0+1}$, and $s_{j_0}$.
It follows by elimination that this space must contain (at least) one of 
the following:
$s_{j_0}$ plus a (possibly zero) multiple of $s^1_{j_0-1}$; 
$s_{j_0}$ plus a (possibly zero) multiple of $s^N_{j_0+1}$; 
$s^1_{j_0-1}$ and $s^N_{j_0+1}$.
The first case means that we can take $S_2'=\{1,\dots,N\}$, while in the
second we get a valid choice of $s'''$.
In the third case, we begin with $s^N_{j_0+1}$, and iteratively twist the
multidegree as before. For $i'>i_1$, we change
$w'$ from twisting down by $a^{i'}_{j_0}$ to $a^{i'}_{j_0+1}$, and at each 
stage, we must either obtain the desired $s'''$, or a section made up purely 
of the $s^i_{j_0+1}$, in which case we continue to iterate. Note that in
these multidegrees, we continue to have that the only possible contributions
of the $s^i_j$ (for $j \neq j_0$) supported strictly left of $i'$ come for 
$j \leq j_0-1$, and we can take the image of $s^1_{j_0-1}$ from multidegree
$\omega_{j_0}$, 
so all these can be subtracted off as necessary.
When $i'\leq i_1$, we will have $a^{i'}_{j_0-1}$ between $a^{i'}_{j_0}$ and
$a^{i'}_{j_0+1}$; we still iteratively increase the twist, but a new
possibility occurs: once we are twisting down by strictly more than 
$a^{i'}_{j_0-1}$, we could obtain a contribution from $a^{i'-1}_{j_0-1}$.
If this occurs, we will continue to iterate, but stopping after increasing
the twist from $a^{i'}_{j_0}$ to $a^{i'}_{j_0-1}$ for each smaller $i'$.

If we have continued with contributions from $s^{i'}_{j_0+1}$ for each $i'$,
then once we reach $i_0$, we will again have no other $a^i_j$ between
$a^i_{j_0}$ and $a^i_{j_0+1}$, so we will ultimately obtain an $s'''$ of
the desired form, with $S_2'''=\emptyset$. On the other hand, if we have
switched from the $s^{i'}_{j_0+1}$ to the $s^{i'}_{j_0-1}$, then we
see that this must terminate (necessarily with an $s'''$ of the desired
form) before we reach $i'=i_0$, because there is no section in column $i_0$
which can glue to $s^{i_0+1}_{j_0-1}$.

Now, if the above construction did not give $s'''$ because we had 
$S_2'=\{1,\dots,N\}$, we apply precisely the same process starting in
multidegree $\omega_{j_0+1}$, and we find that unless we also have 
$S_4''=\{1,\dots,N\}$, we end up with the desired $s'''$.
\end{proof}

Up until now, everything we have done has been insensitive to insertion
of genus-$0$ components. However, to handle the genus-$23$ case, we will
need to impose restrictions on direction of approach; more precisely,
we will require that the genus-$1$ components be separated by exponentially
increasing numbers of genus-$0$ components (going from right to left). 
The reason for doing this is that, if our limit linear series has all 
changes to the $\lambda_i$ occurring in the genus-$1$ components, the 
pattern of the genus-$0$ components will force the support of every $s_j$ 
in every multidegree to be precisely the leftmost segment of potential
support (see Proposition \ref{prop:left-weighted} below), so we obtain better 
control over the situation when the potential support is disconnected. 
However, when dealing with mixed sections, if
the transition from one row to another happens inside the chain of genus-$0$
components, we may lose control over the support. Thus, it turns out to
be important to also analyze what extra control we obtain on our mixed
sections by restricting the direction of approach.
That this sort of restriction could potentially be useful is already pointed 
out in Remark 4.12 of \cite{o-l-t-z}, but our approach is also influenced
by the earlier work of Jensen and Payne \cite{j-p2} on their tropical
approach to the classical maximal rank conjecture.

It turns out that it is convenient to count not the number of genus-$0$
components between a pair of genus-$1$ components, but rather the number
of nodes. When we consider restricted directions, we will assume that the
first and last components have genus $1$, and we will denote by
$\ell_1,\dots,\ell_{N-1}$ the number of nodes between each consecutive
pair of genus-$1$ components, so that $\ell_i-1$ is the number of genus-$0$
components. The reason why this is more convenient is that if we take
a ramified base change with ramification index $e$, and then blow up to
resolve the resulting singularities, we will insert $e-1$ new genus-$0$
components at every node, which has the effect of multiplying all the
$\ell_i$ by $e$. Thus, the ratios of the $\ell_i$ are invariant under
this operation.

\begin{defn}\label{defn:special-dir} We say that $X_0$ is 
\textbf{left-weighted} if we have
$$\ell_i \geq 4d \sum_{i'=i+1}^{N-1} \ell_{i'}.$$
\end{defn}

\begin{defn}\label{defn:mixed-controlled} 
Given $\vec{S}=(S_1,\dots,S_{\ell})$ and
$\vec{j}=(j_1,\dots,j_{\ell})$, a mixed section of type $(\vec{S},\vec{j})$
is said to be \textbf{controlled} if for every $i=2,\dots,\ell$ with
$S_{i} \neq \emptyset$, the
minimal element of $S_{i}$ is 
either a genus-$1$ component or strictly closer to the next genus-$1$ 
component to the right than to the previous one on the left.
\end{defn}

\begin{prop}\label{prop:mixed-controlled} Suppose that $X_0$ is 
left-weighted. Then:

\begin{enumerate}
\item In the situation of Proposition \ref{prop:disjoint-swap}, if we assume
further that 
$i_0$ and $i_1$ have genus $1$, then we may require that
$s_{j_0}'$ and $s_{j_1}'$ are controlled, that $S_2'$ does not contain any 
$i<i_0$ which has genus $1$, and that $S_4'$ does not contain any $i<i_1$
which has genus $1$.
\item In the situation of Proposition \ref{prop:3-cycle-2}, if we assume
further that $i_0$ and $i_1$ have genus $1$, then we may require
that $s_{j_0}'$ 
is controlled, and that $S_2'$ does not 
contain any $i<i_0$ which has genus $1$. 
\end{enumerate}
\end{prop}

\begin{proof} (1) 
For $s=0,1$, in the proof of Proposition \ref{prop:disjoint-swap}, 
every value of $w'$ arising in the iterative procedure will have
that the potential support of the $s^i_{j_s-1}$ in multidegree $\md(w')$ has
two connected components: one extending from $i=1$ to $i=i'-1$, and the
other supported at $i=N$. The reason that we cannot continue our iterative
procedure indefinitely is that we may have that 
$f_{w_{j_s-1},w'} (s_{j_s-1})$ is supported (partially or entirely) at $i=N$.
If we write $w'=(c_2',\dots,c_N')$, we will have that $a^i_{j_s-1}-c_i'>0$
for $i>i_s$, and $a^i_{j_s-1}-c_i'<0$ for $i' \leq i \leq i_s$, and
$a^i_{j_s-1}-c_i'=0$ for $i<i'$.
Suppose $i_s$ is the $m_s$th genus-$1$ component. Then in the notation used
in Definition \ref{defn:special-dir} above, 
we can say (extremely conservatively) that
\begin{equation}\label{eq:mixed-controlled}
\sum_{i=i_s+1}^N (a^i_{j_s-1}-c_i') \leq d \sum_{i=m_s}^{N-1} \ell_i
\leq \frac{\ell_{m_s-1}}{4}.
\end{equation}
Thus, if $i_s-i' > \frac{\ell_{m_s-1}}{4}$, then 
$\sum_{i=i'}^N (a^i_{j_s-1}-c_i') <0$, so we
certainly have that
$f_{w_{j_s-1},w'} (s_{j_s-1})$ is supported entirely on the left, so that
we can subtract off any contribution from the $s^i_{j_s-1}$ and continue
our iterative procedure. The desired conditions on $s_{j_s}'$ follow.

(2) This is essentially the same as (1) (the analogous statement for
$s_{j_0+1}'$ is a bit more complicated, but we don't need it). 
\end{proof}

\section{General setup}\label{sec:setup}

We now describe the basic situation for taking the tensor square of a
limit linear series, and considering images in a fixed multidegree of
total degree $2d$. We will specialize to the case that the components $Z_i$
all have genus at most $1$, but we begin by extending our discrete data 
from the base limit linear series to its tensor square.

\begin{notn}\label{notn:tensor-table} In the situation of Notation 
\ref{notn:lls-table}, let $T$ be the $\binom{r+2}{2} \times N$ table
with rows indexed by unordered pairs $(j,j')$ with $j,j' \in \{0,\dots,r\}$,
having entries $(a^i_{(j,j')}, b^i_{(j,j')})$ defined by
$$a^i_{(j,j')}=a^i_j+a^i_{j'}, \quad \text{ and } \quad 
b^i_{(j,j')}=b^i_j+b^i_{j'}.$$
\end{notn}

The following definition controls when sections can have nonzero image in
a given multidegree. Note that Proposition \ref{prop:vanishing-comps} 
involved only
multidegrees and not limit linear series, so we can apply it equally well
with $2d$ in place of $d$. This then motivates the following definition.

\begin{defn}\label{defn:potential-rows} In the situation of Notation
\ref{notn:tensor-table}, fix total degree $2d$, and $w=(c_2,\dots,c_N)$.
We say that the $(j,j')$ row is \textbf{potentially appearing}
(respectively \textbf{potentially starting}, respectively 
\textbf{potentially ending}) in column $i$ and multidegree $\md_{2d}(w)$ if
$a^i_{(j,j')} \geq c_i$ and $b^i_{(j,j')} \geq 2d-c_{i+1}$
(respectively
$a^i_{(j,j')} > c_i$ and $b^i_{(j,j')} \geq 2d-c_{i+1}$,
respectively
$a^i_{(j,j')} \geq c_i$ and $b^i_{(j,j')} > 2d-c_{i+1}$).
\end{defn}

Now, we will specialize to the genus-$1$ case:

\begin{sit}\label{sit:chain-2} In Situation \ref{sit:chain}, suppose
further that every $Z_i$ has genus at most $1$. Fix $d \geq 0$, and
suppose also that all the $P_i$ and $Q_i$ are general (and in particular
each $P_i-Q_i$ is not $\ell$-torsion for any $\ell \leq d$).
\end{sit}

Note that our generality condition cannot be imposed component by component,
but also involves interaction between components; this arises in the
proof of Lemma \ref{lem:independence-crit} below.

In any genus, we always have $a^i_j+b^i_j \leq d$ for all $j$; if
$Z_i$ has genus $0$, there are no further restrictions, but under the
genericity hypothesis of Situation \ref{sit:chain-2} we have that if the
genus of $Z_i$ is equal to $1$, we can have $a^i_j+b^i_j=d$ for exactly 
one value of $j$, and in this case the underlying line bundle is uniquely
determined as $\sO(a^i_j P_i + b^i_j Q_i)$. The generic situation is
that $a^i_j+b^i_j=d-1$ for all other $j$, but in positive codimension we
can have strictly smaller sums as well -- see the proof of Proposition 2.1
of \cite{os26} for an analysis of these codimensions.
As compared to \cite{o-l-t-z},
we have to consider arbitrary refined limit linear series, 
allowing columns to sum to $d-2$ or less, and to switch orders.
Summing to $d-2$ or less means that where sections actually appear in a 
given multidegree can be more complicated to analyze,
even with well-behaved (i.e., `unimaginative') multidegrees. However, the 
key point is that
the natural necessary and sufficient conditions for a row to appear in
a given column in a given multidegree in the case that all sums are at
least $d-1$ still gives a necessary condition in full generality.
When we have swaps, the limit linear
series are not chain-adaptable, so the linked linear series living over
them can be non-simple, and in fact quite degenerate. 

\begin{defn}\label{defn:exceptional} For a given limit linear series, we 
say that the $j$th row is \textbf{exceptional} in column $i$ if
$a^i_j+b^i_j<d-1$ when $Z_i$ has genus $1$, or if $a^i_j+b^i_j<d$ when $Z_i$
has genus $0$.
\end{defn}

While we imagine starting from a chain of genus-$1$ curves, we allow 
for inserting any number of rational components at nodes, 
so that we will have $N$ components, of which exactly $g$ will have
genus $1$ (including the first and last components) and $N-g$ will 
have genus $0$. Given $i$ between $1$ and $N$, we will denote by $g(i)$
the number of genus-$1$ components between $1$ and $i$, inclusive.

It will be convenient to package our 
sequences $a^i_j$ of vanishing sequences in a slightly different form,
as follows.

\begin{notn} For $j=0,\dots,r$ and $i=0,\dots,N$, we can determine an
integer (possibly negative) $\lambda_{i,j}$ by
$$a^{i+1}_j=g(i)+j-\lambda_{i,j}.$$ 
If $\lambda_{i,j}>\lambda_{i-1,j}$, 
we say $\delta_{i}=j$; otherwise, we say there is no $\delta_i$.
\end{notn}

Here $g(0)=0$ by convention, and we also set the convention that
$a^{N+1}_j=b^N_{r-j}$ for all $j$.

Intuitively, we think of the $\lambda_{i,j}$ as forming
a sequence
$\lambda_i$ of generalized `shapes' (not necessarily skew, or connected), 
behaving as follows: $\lambda_{0,j} \leq 0$ for all $j$; any number of 
`squares' can be removed from
the righthand side in each step, and at most one `square' is added at each 
stage, with the possibility of adding a `square' only in the genus-$1$
components. It is however possible for the $\lambda_{i,j}$ to be negative,
either because we are starting with nontrivial ramification sequences, or
because we remove too many squares early on.
Because the $a^i_j$ are always distinct for a fixed $i$, we see that we 
will always have 
the set $j-\lambda_{i,j}$ consisting of $r+1$ distinct integers.

\begin{rem}\label{rem:rho-meaning}
Before discussing tensor squares, we briefly recall the significance of
$\rho$ in this setup. We need to have $b^N_j$ nonnegative (and distinct)
for all $j$, or equivalently $a^{N+1}_j$ bounded by $d-r,d-r+1,\dots,d$.
In particular, $\sum_{j=0}^r a^{N+1}_j \leq (r+1)d-\binom{r+1}{2}$, so
$$\sum_{j=0}^r \lambda_{N,j} \geq(r+1)g+\binom{r+1}{2}-(r+1)d+\binom{r+1}{2}
=(r+1)(g+r-d)=g-\rho.$$ 
Since $\sum_j \lambda_{i,j}$ can increase by at most $1$ as $i$ increases
(and only on genus-$1$ components), and $\lambda_{0,j} \leq 0$ for all $j$,
we see that for $\rho=0$, we must have $\lambda_{0,j}=0$ for all $j$ (i.e.,
minimal initial vanishing sequence at $P_1$), no
places where $\lambda_{i,j}$ decreases (i.e., no exceptional columns for
any row), and a $\delta_i$ for every genus-$1$ column $i$. When $\rho>0$,
the total amount of initial ramification, exceptional columns, and genus-$1$
columns without $\delta_i$ is bounded by $\rho$. A swap is necessarily
a case of an exceptional column, and can contribute exactly $1$ to $\rho$
precisely when it is minimal and occurs in a genus-$1$ column. Note also
that a minimal swap
can occur at most once in a given (genus-$1$) column.
\end{rem}

Moving on to tensor squares,
we now recall the following definition from \cite{o-l-t-z}, updated to
allow for genus-$0$ components:

\begin{defn}\label{defn:unimaginative} We say a multidegree of total degree
$2d$ is
\textbf{unimaginative} if it assigns degree $0$ to every genus-$0$
component, and degree $2$ or $3$ to every genus-$1$ component.
By extension, we will say that $w$ is unimaginative if $\md_{2d}(w)$ is.
Given a fixed unimaginative multidegree, we will
let $\gamma_i$ be the number of $3$s in the first $i$ columns.
\end{defn}

We will work throughout only with unimaginative multidegrees. 
Thus, the multidegree is encoded by twisting down by $2d-2g(i)-\gamma_i$ on
the righthand of the $i$th column, and by twisting down by 
$2g(i)+\gamma_i$ on the lefthand side of the $(i+1)$st column, for all $i<N$.

Then the following is straightforward from the definitions:

\begin{prop}\label{prop:when-appears}
A row $(j_1,j_2)$ is potentially appearing in the $i$th and $(i+1)$st columns
only if 
$$j_1+j_2-\lambda_{i,j_1}-\lambda_{i,j_2}=\gamma_{i}.$$

A row $(j_1,j_2)$ is potentially appearing in the $i$th column only if
$$j_1+j_2-\lambda_{i-1,j_1}-\lambda_{i-1,j_2}\geq \gamma_{i-1}
\quad\text{ and }\quad
j_1+j_2-\lambda_{i,j_1}-\lambda_{i,j_2}\leq \gamma_{i}.$$

A row $(j_1,j_2)$ is potentially starting in the $i$th column only if
$$j_1+j_2-\lambda_{i-1,j_1}-\lambda_{i-1,j_2}> \gamma_{i-1}
\quad\text{ and }\quad
j_1+j_2-\lambda_{i,j_1}-\lambda_{i,j_2}\leq \gamma_{i}.$$

A row $(j_1,j_2)$ is potentially ending in the $i$th column only if
$$j_1+j_2-\lambda_{i-1,j_1}-\lambda_{i-1,j_2}\geq  \gamma_{i-1}
\quad\text{ and }\quad
j_1+j_2-\lambda_{i,j_1}-\lambda_{i,j_2}< \gamma_{i}.$$

A row $(j_1,j_2)$ is potentially starting and ending in the $i$th column 
only if
$$j_1+j_2-\lambda_{i-1,j_1}-\lambda_{i-1,j_2}> \gamma_{i-1}
\quad\text{ and }\quad
j_1+j_2-\lambda_{i,j_1}-\lambda_{i,j_2}< \gamma_{i}.$$
\end{prop}

Note that the sequence $j_1+j_2-\lambda_{i,j_1}-\lambda_{i,j_2}$ 
decreases by at most $1$ each time $i$ increases,
unless $j_1=j_2$, when it can decrease by $2$.
Similarly, $\gamma_i$ is nondecreasing, and increases by at most $1$ each
time $i$ increases. As mentioned previously, for fixed $i$, the
$j-\lambda_{i,j}$ consists of $r+1$ distinct values.

\begin{cor}\label{cor:2-case} If the multidegree has a $2$ in the $i$th 
column, then there can be at most one row potentially starting in it, and at 
most one row potentially ending in it. 

There can be a row potentially starting in the $i$th column only if 
$\delta_{i}$ exists and either there exists $j$ such that
$$\delta_{i}+j-\lambda_{i,\delta_{i}}-\lambda_{i,j}
=\gamma_{i}$$
or
$$2\delta_{i}-2\lambda_{i,\delta_{i}}=\gamma_{i}-1.$$
In these cases, the potentially starting rows are $(\delta_{i},j)$ or
$(\delta_{i},\delta_{i})$, respectively.

There can be a row potentially ending in the $i$th column only if 
$\delta_{i}$ exists and either there exists $j$ such that
$$\delta_{i}+j-\lambda_{i,\delta_{i}}-\lambda_{i,j}
=\gamma_{i}-1$$
or
$$2\delta_{i}-2\lambda_{i,\delta_{i}}=\gamma_{i}-2.$$
In these cases, the potentially ending rows are $(\delta_{i},j)$ or
$(\delta_{i},\delta_{i})$, respectively.
\end{cor}

\begin{proof} Since in this case $\gamma_{i}=\gamma_{i-1}$, Proposition 
\ref{prop:when-appears} implies that the $(j_1,j_2)$ row can be potentially
starting in the 
$i$th column only if $\lambda_{i,j_1}>\lambda_{i-1,j_1}$ or
$\lambda_{i,j_2}>\lambda_{i-1,j_2}$, which is to say if $\delta_{i}$
exists and $j_1$ or $j_2$ is equal to $\delta_{i}$. Moreover, in this
case $\lambda_{i,\delta_{i}}=\lambda_{i-1,\delta_{i}}+1$, so 
we conclude that the two stated cases are the only possibilities for
having 
$$j_1+j_2-\lambda_{i-1,j_1}-\lambda_{i-1,j_2}> \gamma_{i-1}=\gamma_{i}
\geq 
j_1+j_2-\lambda_{i,j_1}-\lambda_{i,j_2},$$
and that moreover in the first case we must also have 
$\lambda_{i-1,j}=\lambda_{i,j}$ unless $j=\delta_{i}$.

Now, there is at most one $j$ satisfying the first identity of the corollary, 
since the $j-\lambda_{i,j}$ are all distinct. Moreover, if there is
some $j$ satisfying the first, then the second one cannot hold, since
this would force
$$\delta_{i}-\lambda_{i-1,\delta_{i}}=
\delta_{i}-\lambda_{i,\delta_{i}}+1=j-\lambda_{i,j}=
j-\lambda_{i-1,j},$$
which is not allowed. This completes the proof of the assertions on rows
potentially starting in the $i$th column, and the assertion on rows 
potentially ending in the $i$th column is proved similarly.
\end{proof}

The following corollary has a similar proof, which we omit.

\begin{cor}\label{cor:3-case} If the multidegree has a $3$ in the $i$th 
column, then there can be at most one row potentially starting and ending in 
the $i$th column, and this occurs
only if $\delta_{i}$ exists and either there exists $j$ such that
$$\delta_{i}+j-\lambda_{i,\delta_{i}}-\lambda_{i,j}
=\gamma_{i}-1$$
or
$$2\delta_{i}-2\lambda_{i,\delta_{i}}=\gamma_{i}-2.$$

In addition, for a fixed $j \neq \delta_{i}$, there is at most one value
of $j'$ such that the $(j,j')$ row is potentially starting in column $i$,
and at most one value of $j'$ such that the $(j,j')$ row is potentially 
ending in column $i$.
\end{cor}

Now, given a refined limit linear series, we can also construct a second
table $\bar{T}$ of vanishing numbers which is obtained from the first simply 
by reordering each subcolumn into strict increasing (respectively, 
decreasing) order. Put differently, $\bar{T}$ is obtained from the limit 
linear series simply by taking vanishing sequences at each point, and 
ignoring the interplay between the pair of points. We will denote the 
$\lambda$ sequence obtained from $\bar{T}$ by $\bar{\lambda}_i$, and
the entries of the table $\bar{T}$ by $(\bar{a}^i_j, \bar{b}^i_j)$.
Here, if we picture skewing the rows of the $\bar{\lambda}_i$ according
to the initial ramification sequence $a^1_j-j$, the sequence 
$\bar{\lambda}_i$ will give a genuine sequence of skew shapes, terminating
with a skew shape containing the one obtained by starting from the usual 
$(r+1)\times (r+g-d)$ center rectangle, and adding squares
on the left determined by the initial ramification sequence.

For $\ell \geq 1$, we denote by
$\bar{\lambda}_i^{\ell}$ the number of $j$ such that 
$\bar{\lambda}_{i,j} \geq \ell$, which we can visualize as the number of 
squares in the $\ell$th column of $\bar{\lambda}_i$,
numbered so that the ``first column'' is the first column of the main 
$(r+1) \times (r+g-d)$ rectangle.\footnote{Although this definition could
in principle be applied also to $\lambda_i$, it does not seem to have any
particular significance when $\lambda_i \neq \bar{\lambda}_i$.}

The following lemma is the key to our analysis, showing in particular
that if we place multidegree $3$ in the correct places, we can obtain fine 
control over what happens with the rows involving $\delta_{i+1}$.

\begin{lem}\label{lem:impossible}
Given $1 \leq \ell_1 < \ell_2$ and $n>0$, suppose that 
$\bar{\lambda}_{i}^{\ell_1}+\bar{\lambda}_{i}^{\ell_2}=n$, but
$\bar{\lambda}_{i-1}^{\ell_1}+\bar{\lambda}_{i-1}^{\ell_2}<n$.
Then there does not exist a $j$ such that 
\begin{equation}\label{eq:impossible-1}
\delta_{i}+j-\lambda_{i,\delta_{i}}-\lambda_{i,j}
=n-1-\ell_1-\ell_2,
\end{equation}
and we do not have
\begin{equation}\label{eq:impossible-2}
2\delta_{i}-2\lambda_{i,\delta_{i}}=n-2-\ell_1-\ell_2.
\end{equation}

Moreover, if we have some $j$ with $\lambda_{i,j}<\lambda_{i-1,j}$,
then we cannot have
\begin{equation}\label{eq:semipossible-1}
\delta_{i}+j-\lambda_{i,\delta_{i}}-\lambda_{i,j}
=n-\ell_1-\ell_2
\end{equation}
or 
\begin{equation}\label{eq:semipossible-2}
\delta_{i}+j-\lambda_{i-1,\delta_{i}}-\lambda_{i-1,j}
=n-1-\ell_1-\ell_2.
\end{equation}
\end{lem}

\begin{proof} We first prove the case that $\bar{\lambda}_{i'}=\lambda_{i'}$ 
for all $i'$. Note that we necessarily have a $\delta_{i}$, and it
must be the row of the lowest square in either the $\ell_1$th or
$\ell_2$th column of $\lambda_{i}$. 
Note also that if
$\lambda_{i,j}<\lambda_{i-1,j}$ for some $j$, then since we assumed that
$\lambda_{i-1}^{\ell_1}+\lambda_{i-1}^{\ell_2}<n$, we cannot have
$\lambda_{i,j}$ equal to $\ell_1-1$ or $\ell_2-1$. Thus,
we will prove the desired 
statement on \eqref{eq:semipossible-1} 
by proving that if \eqref{eq:semipossible-1} is satisfied for any $j$, then 
we must have $\lambda_{i,j}$ equal to $\ell_1-1$ or $\ell_2-1$.

First consider the case that $\lambda_{i}^{\ell_1},\lambda_{i}^{\ell_2}$ 
are distinct and positive, and set
$j_s=\lambda_{i}^{\ell_s}-1$ for $s=1,2$.
Thus, we necessarily have $\delta_{i}=j_1$ or $j_2$.
Note that $(j_1+1)+(j_2+1)=n$ by hypothesis, and
for $s=1,2$ write $m_s=\lambda_{i,j_s}-\ell_s$, so that
necessarily $m_s \geq 0$ for $s=1,2$, with equality for (at least) one $s$.
It follows that we have
\begin{multline*}j_1+j_2-\lambda_{i,j_1}-\lambda_{i,j_2}
=(j_1+1)+(j_2+1)-2-(m_1+\ell_1)
-(m_2+\ell_2)\\
=n-2-\ell_1-\ell_2-m_1-m_2<n-1-\ell_1-\ell_2.\end{multline*}
Thus, the only way to get \eqref{eq:impossible-1} would be to set $j$
to be strictly greater than whichever $j_s$ is not equal to $\delta_{i}$.
Now, because $j_s$ was determined as
the lowest row with a square in the $\ell_s$th column, we have
$$\lambda_{i,j_s+1}<\ell_s
=\lambda_{i,j_s}-m_s,$$
so if we use $j>j_s$ in place of $j_s$,
the value of the above expression jumps by at least $2+m_s$. Moreover,
we can only use $j$ in place of $j_1$ if $\delta_{i}=j_2$, in which
case we must have $m_2=0$, and similarly if we use $j$ in place of
$j_2$, so we conclude that \eqref{eq:impossible-1} is not possible. 
We also see that if we have \eqref{eq:semipossible-1}, then necessarily
$j=j_s+1$ and $\lambda_{i,j}=\ell_s-1$, as asserted.
By the same reasoning, if $\delta_{i}=j_1>j_2$, then 
\eqref{eq:impossible-2} is also impossible, because $m_1=0$ replacing $j_2$ by
$\delta_{i}$ increases the lefthand side by at least $2+m_2$. 
On the other hand, if $\delta_{i}=j_2$ then replacing $j_1$ by 
$\delta_{i}$ decreases the lefthand side, making it too small to satisfy
\eqref{eq:impossible-2}.

Finally, suppose we have some $j$ such that 
$\lambda_{i,j}<\lambda_{i-1,j}$; say 
$\lambda_{i,j}=\lambda_{i-1,j}-p$ for some $p>0$. Then
\eqref{eq:semipossible-2} is equivalent to 
$$\delta_{i}+j-\lambda_{i,\delta_{i}}-\lambda_{i,j}
=n-2-\ell_1-\ell_2+p,$$
so if as above we have $j_s \neq \delta_{i}$, then necessarily
$j>j_s$, so that by definition of $j_s$ we must have
$\lambda_{i,j} < \ell_s$.
On the other hand, since we have assumed that
$\lambda_{i-1}^{\ell_1}+\lambda_{i-1}^{\ell_2}<n$, we must have that
$\lambda_{i,j},\dots,\lambda_{i,j}+p$ does not
contain $\ell_s$, so it follows that
$\lambda_{i,j}+p<\ell_s=\lambda_{i,j_s}-m_s$.
We conclude that
$j-\lambda_{i,j}>1+j_s-\lambda_{i,j_s}+p+m_s$, so
$$\delta_{i}+j-\lambda_{i,\delta_{i}}-\lambda_{i,j}
> (n-2-\ell_1-\ell_2-m_s)+1+p+m_s
=n-1-\ell_1-\ell_2+p,$$
proving the desired impossibility of \eqref{eq:semipossible-2}.

The next case is that $\lambda_{i}$ has no entries in the $\ell_2$th
column, so that $\delta_{i}+1=n$, and 
$\lambda_{i,\delta_{i}}=\ell_1$.
In this case, we have
$$\delta_{i}-\lambda_{i,\delta_{i}}
=(\delta_{i}+1)-1-\ell_1
=n-1-\ell_1.$$
But since the $\ell_2$th column is empty, for all $j$ we have
$\lambda_{i,j}<\ell_2$, so we find that
$$\delta_{i}+j-\lambda_{i,\delta_{i}}-\lambda_{i,j}>
n-1-\ell_1+j-\ell_2\geq n-1-\ell_1-\ell_2,$$
showing that \eqref{eq:impossible-1} cannot hold, and that 
\eqref{eq:semipossible-1} can hold only if 
$\lambda_{i,j}=\ell_2-1$.
We also see that
$$2\delta_{i}-2\lambda_{i,\delta_{i}}
=2n-2-2\ell_1>2n-2-\ell_1-\ell_2>n-2-\ell_1-\ell_2,$$
so \eqref{eq:impossible-2} does not hold either.
Finally, because $\lambda_{i-1}^{\ell_1}+\lambda_{i-1}^{\ell_2}<n$ 
we necessarily have also $\lambda_{i-1,j}<\ell_2$, so 
$$\delta_{i}+j-\lambda_{i-1,\delta_{i}}-\lambda_{i-1,j}
= \delta_{i}+j-\lambda_{i,\delta_{i}}-\lambda_{i-1,j}+1>
n-1-\ell_1-\ell_2+1=n-\ell_1-\ell_2,$$
proving that \eqref{eq:semipossible-2} also cannot hold.

The final case is that $\lambda_{i}$ has the same number of entries in
the $\ell_1$th and $\ell_2$th columns, so that we must have $n$ even,
with $\delta_{i}+1=n/2$, and also
$\lambda_{i,\delta_{i}}=\ell_2$.
In this case, we have
$$\delta_{i}-\lambda_{i,\delta_{i}}
=(\delta_{i}+1)-1-\ell_2
=n/2-1-\ell_2.$$
Thus if we set $j_1=j_2=\delta_{i}$, we find that
$$j_1+j_2-\lambda_{i,j_1}-\lambda_{i,j_2}
=n-2-2\ell_2<n-2-\ell_1-\ell_2,$$
so \eqref{eq:impossible-2} does not hold.
But because the $\ell_1$th column has exactly $\delta_{i}+1$ entries,
leaving $j_2=\delta_{i}$ and using $j_1>\delta_{i}$ results in an 
increase of at least $2+\ell_2-\ell_1$, 
yielding
$$j_1+j_2-\lambda_{i,j_1}-\lambda_{i,j_2}
\geq n-\ell_2-\ell_1,$$
so we see that \eqref{eq:impossible-1} also cannot be satisfied.
Moreover, we can have \eqref{eq:semipossible-1} only if $j=\delta_{i}+1$
and $\lambda_{i,j}=\ell_1-1$. Finally, as in the first case
considered, if $\lambda_{i,j}=\lambda_{i-1,j}-p$ for $p>0$, then in order to 
have \eqref{eq:semipossible-2} we would need to have $j>j_2$, which then
implies that 
$\lambda_{i,j}+p<\ell_1
=\lambda_{i,\delta_{i}}-\ell_2+\ell_1$, so
$$\delta_{i}+j-\lambda_{i,\delta_{i}}-\lambda_{i,j}
> (n-2-2\ell_2)+1+(\ell_2-\ell_1+p) =n-1-\ell_1-\ell_2+p,$$
again yielding that \eqref{eq:semipossible-2} is not possible.

This completes the proof of the lemma
in the case that $\bar{\lambda}_{i'}=\lambda_{i'}$ for all $i'$. We will see 
that the general case follows. The main observation is the following: if
$\bar{\lambda}_{i,j}=\bar{\lambda}_{i-1,j}+1$, and we let $j'$ be
such that $\bar{a}^{i+1}_{j}=a^{i+1}_{j'}$,
then we necessarily have 
$\lambda_{i,j'}=\lambda_{i-1,j'}+1$, and we cannot have any swaps in the
$i$th column involving the $j'$th row. Indeed, the identity
$\bar{\lambda}_{i,j}=\bar{\lambda}_{i-1,j}+1$ means that we have
$\bar{a}^{i}_j=\bar{a}^{i+1}_j$,
which means that $a^{i+1}_{j'}=a^{i}_{j''}$ for the $j''$ such that
exactly $j$ values of $a^{i}_m$ are less than $a^{i}_{j''}$.
We also have exactly $j$ values of $a^{i+1}_{m}$ less than
$a^{i+1}_{j'}$. It then follows that we must have $j''=j'$: we cannot
have $a^{i}_{j'}>a^{i}_{j''}$, since then we would have
$a^{i}_{j'}>a^{i+1}_{j'}$. But if $a^{i}_{j'}<a^{i}_{j''}$,
then $j'$ occurs among the values of $m$ with $a^{i}_m<a^{i}_{j''}$,
so there is necessarily some $m$ with $a^{i+1}_m<a^{i+1}_{j'}$ but
$a^{i}_m \geq a^{i}_{j''}$, again leading to a contradiction.
This proves the observation, noting that the fact that $j'=j''$ rules
out any swaps involving the $j'$th row. 

We then conclude that the impossibility
of \eqref{eq:impossible-1} and \eqref{eq:impossible-2} reduces to the
case that $\bar{\lambda}_{i'}=\lambda_{i'}$, since both equations can
be phrased in terms of the values of $j-\lambda_{i',j}=a^{i'+1}_j-g(i')$,
and our above observation implies that we have 
$a^{i}_{\delta_{i}}=a^{i+1}_{\delta_{i}}
=\bar{a}^{i+1}_{\bar{\delta}_{i}}=\bar{a}^{i}_{\bar{\delta}_{i}}$
(here we use $\bar{\delta}_i$ to denote the values of $\delta$ coming from
$\bar{T}$).
Next, suppose that we have some $j$ with $\lambda_{i,j}<\lambda_{i-1,j}$;
we claim that if $j'$ is such that 
$a^{i+1}_j=\bar{a}^{i+1}_{j'}$, and $j''$ is such that 
$a^{i}_j=\bar{a}^{i}_{j''}$, then we
necessarily also have that $\bar{\lambda}_{i,j'}<\bar{\lambda}_{i-1,j'}$
and $\bar{\lambda}_{i,j''}<\bar{\lambda}_{i-1,j''}$. Given this claim,
the impossibility of \eqref{eq:semipossible-1} and \eqref{eq:semipossible-2}
follows from the case that $\bar{\lambda}_{i'}=\lambda_{i'}$ for all $i'$. 
By our above observations on 
the case $\bar{\lambda}_{i,j}=\bar{\lambda}_{i-1,j}+1$, it suffices to prove
that $\bar{\lambda}_{i,j'} \neq \bar{\lambda}_{i-1,j'}$ 
and $\bar{\lambda}_{i,j''}\neq \bar{\lambda}_{i-1,j''}$,
or equivalently, 
that $\bar{a}^{i+1}_{j'} \neq \bar{a}^{i}_{j'}+1$,
and
$\bar{a}^{i+1}_{j''} \neq \bar{a}^{i}_{j''}+1$.
But in order to have 
$a^{i+1}_j=\bar{a}^{i+1}_{j'} = \bar{a}^{i}_{j'}+1$, we would need to have
$a^{i+1}_j-1$ occurring among the $a^{i}_{\bullet}$, with precisely
$j'$ strictly smaller values also occurring. But by definition we have $j'$
values strictly smaller than $a^{i+1}_j$ occurring in $a^{i+1}_{\bullet}$,
and using our observation on lack of swaps when 
$\bar{\lambda}_{i,j}=\bar{\lambda}_{i-1,j}+1$ we see that every one of
these also must yield a value of $a^{i}_{\bullet}$ strictly smaller than
$a^{i+1}_j-1$.
But we have in addition that $a^{i}_j<a^{i+1}_j-1$, so we conclude that
there are at least $j'+1$ values in $a^{i}_{\bullet}$ strictly less
than $a^{i+1}_j-1$, proving the desired inequality by contradiction.

Similarly, in order to have
$\bar{a}^{i+1}_{j''} = \bar{a}^{i}_{j''}+1=a^{i}_j+1$, we would need to 
have $a^{i}_j+1$ occurring among the $a^{i+1}_{\bullet}$, with precisely
$j''$ strictly smaller values also occurring. By definition, we have
only $j''$ values among the $a^{i}_{\bullet}$ strictly smaller than
$a^{i}_j$, and every value of 
$a^{i+1}_{\bullet}$ which is strictly smaller than $a^{i}_j+1$ must come
from one of these.
But again using our observation on the lack of swaps when 
$\bar{\lambda}_{i,j}=\bar{\lambda}_{i-1,j}+1$, we see that the value 
$a^{i}_j+1$ in $a^{i+1}_{\bullet}$ must itself come from a row in
$a^{i}_{\bullet}$ with value strictly smaller than $a^{i}_j$, so
we conclude that if $a^{i}_j+1$ occurs in $a^{i+1}_{\bullet}$, there
must be strictly fewer than $j''$ entries in $a^{i+1}_{\bullet}$ which
are strictly smaller than it. This proves the claim, and the lemma.
\end{proof}

\section{An independence criterion}\label{sec:crit}

Suppose we have a limit linear series, and fix choices of sections
$s^i_j$ matching the vanishing orders in our table. We make the following
definition:

\begin{defn}\label{defn:pot-secs} Given an unimaginative multidegee $\omega$, 
for all $(j_1,j_2)$, let $n_{(j_1,j_2)}$ be the number of places (i.e., 
collections of contiguous columns) where the $(j_1,j_2)$ row could 
potentially appear in the multidegree $\omega$. Let
$s_{(j_1,j_2),i}$ for $i=1,\dots,n_{(j_1,j_2)}$ be the induced sections
in multidegree $\omega$ with precisely the given support. Then the full
collection of $s_{(j_1,j_2),i}$ are the \textbf{potentially appearing}
sections in multidegree $\omega$, and their span in 
$\Gamma(X_0,(\sL^{\otimes 2})_\omega)$ is the \textbf{potential ambient space}.
\end{defn}

Note that in the above, we require that each $s_{(j_1,j_2),i}$ be
potentially starting in its first column of support and potentially ending
in its last column of support. Thus, there may be individual columns in
which the $(j,j')$ row satisfies the inequalities to potentially appear
in that column, but which does not occur in any of the $s_{(j_1,j_2),i}$
because it fails necessarily inequalities in other columns.

The $s_{(j_1,j_2),i}$ are each unique up to scaling given a choice of the
$s^i_j$. The $s^i_j$
are not unique, but they can differ only by multiples of $s^i_{j'}$ with
strictly higher vanishing at both points. Then if $s^i_j$ has potential
support (in the $i$th column), necessarily $s^i_{j'}$ has a connected
component of potential support consisting precisely of the $i$th column.
We conclude that the potential ambient space is independent of the 
choice of the $s^i_j$. Consequently, the dimension of the span -- and
in particular the linear independence -- of
the potentially appearing sections is likewise independent of choices.

As in \cite{o-l-t-z}, we will give an elementary independence criterion in 
given multidegrees, stated in terms of iterated dropping of sections. However,
while in \cite{o-l-t-z} we determined the image of each $s_j \otimes s_{j'}$
in multidegree $\omega$ and phrased our criterion for linear independence in
terms of dropping rows, in order for us to handle degenerate cases it will
be important to shift our attention from rows to potentially appearing 
sections. The below definition is to be applied during the iterative
procedure, so refers to ``remaining'' sections (i.e., those which have not
yet been dropped).

\begin{defn} We say that the $i$th column of $T$ is \textbf{semicritical} in
multidegree $\omega$ if it satisfies the following conditions:
\begin{itemize}
\item it has a value of $\delta_i$ (in particular, it has genus $1$); 
\item the minimal values among the potentially appearing sections remaining
in the two subcolumns of column $i$ add to at least $2d-2$;
\item if the $(j,\delta_i)$ row remains in the $i$th column for some 
$j\neq \delta_i$, then the $j$th row is not exceptional. 
\end{itemize}
If further the minimal values among the remaining potentially appearing
sections are not both one less than the values in the $(\delta_i,\delta_i)$ 
row, we say that the $i$th column is \textbf{critical}.
\end{defn}

The following is our criterion for checking that the potentially appearing
sections are linearly independent in a given multidegree.

\begin{lem}\label{lem:independence-crit} For a given limit linear series,
and given unimaginative multidegree $\omega$, suppose that we can drop all
potentially appearing sections by iterative application of the following 
rules:
\begin{ilist}
\itm if in some column $i$, there is a unique remaining potentially appearing
section supported in that column having minimal $a^i_{(j,j')}$ value, or
a unique one having minimal $b^i_{(j,j')}$, then the one achieving the
minimum may be dropped;
\itm if there are at most two remaining potentially appearing
sections with support in some genus-$1$ column $i$, and neither of them 
involves an exceptional row, then they can both be dropped;
\itm if there are $i<i'$ such that
the block of columns from $i$ to $i'$ has the following properties, then
all the remaining potentially appearing sections supported in this block
can be dropped:
\begin{itemize}
\item there are at most $3$ remaining potentially appearing sections
supported in each of the $i$th and $i'$th columns; 
\item within the block, there are at most three potentially appearing 
sections continuing from any column to the next;
\item every column strictly between $i$ and $i'$ has degree $2$;
\item both the $i$th column and the $i'$th column are semicritical,
and either $i$ is critical with no remaining potentially appearing section 
ending in the $i$th column, or $i'$ is critical with 
no remaining potentially appearing section starting in the $i'$th column.
\end{itemize}
\end{ilist}

Then the potentially appearing sections in multidegree $\omega$ are linearly
independent.
\end{lem}

\begin{proof} Suppose we had a hypothetical linear dependence among the
potentially appearing sections. We claim that in each case (i), (ii), (iii),
the coefficients of the relevant potentially appearing sections would be
forced to vanish. In case (i), this is clear: the uniqueness of the minimal 
value of
$a^i_{(j,j')}$ means that $s^i_{(j,j')}$ vanishes to strictly smaller
order at $P_i$ than any other remaining potentially appearing section,
and similarly for $b^i_{(j,j')}$. In both cases, the coefficient would have
to vanish in any linear dependence. 

In case (ii), we need to see that for a fixed column $i$, any two 
$s^i_{(j,j')}$ have to be linearly independent provided that they do not
involve any exceptional rows.
If either of them involves $\delta_i$, this is
automatic, since either the $a^i_{(j,j')}$ or $b^i_{(j,j')}$ values are
forced to be distinct.
On the other hand, if neither involves $\delta_i$, we claim that the
sections in question must have distinct zeroes on $Z_i$ away from $P_i$
and $Q_i$. Indeed, if we have $a,b,a',b'$ with $a+b=d-1=a'+b'$, 
then the unique sections $s,s'$ of our given line bundle vanishing to order 
at least $a$ at $P_i$ and $b$ at $Q_i$
(respectively, $a'$ at $P_i$ and $b'$ at $Q_i$) have $\dv s=aP_i+bQ_i+R$ and
$\dv s'=a'P_i+b'Q_i+R'$ for some $R,R'$. We see that we have a linear 
equivalence
$R-R' \sim (a'-a)P_i+(b'-b)Q_i$, and if $0 \leq a,a' \leq d$, we see that 
$R\neq R'$ because of our running generality hypothesis on $P_i,Q_i$. Thus,
tensors of different sections of this form always have zeroes in distinct
places on $Z_i$, and must be linearly independent.

For case (iii), note that the condition that the degree is $2$
on every column between $i$ and $i'$ means by Corollary \ref{cor:2-case} 
that there is at most one potentially appearing section starting and at 
most one ending in each of these columns.
Noting that the situation is fully symmetric, 
suppose without loss of generality that $i'$ is critical, with no
remaining potentially appearing sections starting in it. If $i$ or $i'$ has
fewer than three remaining potentially appearing sections, we may use
(ii) to drop these,
and then move iteratively through the rest
of the block, using that at most one potentially appearing section starts
or ends in each column to repeatedly use (i) to drop the remaining sections 
from the block. Thus, suppose that $i$ and $i'$ both have three remaining 
potentially appearing sections. Note also that if any column $i''$ has only 
one potentially appearing section spanning $i''$ and $i''+1$, then
the minimal value in the right subcolumn of $i''$ is necessarily unique, 
so we can use case (i) to drop the section in question. Moreover, there can 
be at most one other potentially appearing section supported in column $i''$
(the one ending there), so we can drop this one as well, and then we can
move iteratively left and right to drop the entire block. Thus, we may
further suppose that every column in the block has at least two potentially 
appearing sections spanning it and the next column.

Next, normalize our sections as follows: scale all sections spanning the 
$i'-1$ and $i'$ column so they agree at $Q_{i'-1}$, and then go back one
column at a time, scaling any newly appearing section so that its value
at the previous node agrees with the value of a section which has already
been fixed.
In this way, we will fix a normalization of all our sections except for
those which are supported in only one column. Although the normalization
depends on some choices, they are of a discrete nature, and can be fixed
based purely on the discrete data of the limit linear series.

Now, consider a hypothetical nonzero linear dependence involving the rows in 
our block. First, the coefficients of the linear dependence cannot vanish 
identically in the remaining potential sections of any column, since 
otherwise the condition that at most one potentially appearing section ends 
in each column would imply that there was a column with exactly one nonzero 
coefficient among its remaining potentially appearing sections.
Next, we see that the coefficients are unique up to simultaneous
scaling for the three potentially appearing sections in column $i$.
Indeed, since we have assumed that $i$ is semicritical, its three 
potentially appearing sections must be pairwise independent.

Since we have at most one new potentially appearing section in each column, 
we find that the coefficient for any new one is always uniquely determined 
by the previous ones. Since there are no new potentially appearing sections
in column $i'$, we find that even before considering this column, we have 
already uniquely determined all of the coefficients (up to simultaneous 
scaling) of all of the potentially appearing sections remaining in the 
block. Moreover, we claim that these coefficients (excluding the ones for
potentially appearing sections supported only in a single column) are 
uniquely determined up to finite indeterminacy
by the marked curves $Z_i,\dots,Z_{i'-1}$ together with the discrete data
of the limit linear series.
Indeed, there are only two ways in which nontrivial moduli can enter the
picture: if there are columns $i''$ between $i$ and $i'-1$ either having no 
$\delta_{i''}$, or having some sections $s^{i''}_j$ which are not uniquely
determined up to scalar. This becomes slightly delicate, since in both these
cases, varying the moduli could affect both the normalization we have
chosen and the linear dependence. However, we will show that in both cases,
there will in fact be only finitely many possibilities which still preserve
the linear dependence. Note that by hypothesis, neither of these nontrivial
moduli occurs in the $i$th column. Note also that we cannot have both
occurring at once, as the $s^{i''}_j$ can only fail to be determined up to
scalar if they involve an exceptional row, and since we have
assumed we have degree $2$ between $i$ and $i'$, these can only appear if
paired with the $\delta_{i''}$ row.

First consider the case that we have no $\delta_{i''}$. Then since we have
degree $2$, every potentially appearing section in column $i''$ must
extend to both the preceding and subsequent columns; in particular, there can
be at most three such sections. If there are fewer than three, they cannot
be independent, leading to an immediate contradiction. If there are three,
say $s_0^{i''},s_1^{i''},s_2^{i''}$, then they are necessarily dependent
with a unique dependence 
$c_0 s_0^{i''}+c_1 s_1^{i''}+c_2 s_2^{i''}=0$ which can be determined by
requiring that it holds at both $P_{i''}$ and $Q_{i''}$. We claim that 
for any fixed choice of $c_0,c_1,c_2$ (not all zero), there can be only
finitely many choices of the line bundle $\sL^{i''}$ such that the resulting
cancellation holds at both points. For this claim, we can renormalize
our sections so that the values of the $s_j^{i''}$ agree at $P_{i''}$,
and we just want to see that the values at $Q_{i''}$ must move 
nondegenerately in $\PP^2$ as $\sL^{i''}$ varies. But this is precisely the
content of Proposition \ref{prop:nondegen-again}.

Next, suppose that we have an exceptional row $j$ involved in column 
$i''$, necessarily paired with the $\delta_{i''}$ row. As before, a linear 
dependence in the $i''$ necessarily has to give cancellation at both 
$P_{i''}$ and $Q_{i''}$. Suppose that the $j$th
row and the $\delta_{i''}$th row have entries $a,b$ and $a',b'$ respectively,
so that $a+b=d-2$ and 
$a'+b'=d$. There are two cases: if $a=a'-1$, so that also $b=b'-1$ (and
$i''$ has a swap in it), then the moduli for the section $s^{i''}_j$ 
consists simply of adding multiples of the section 
$s^{i''}_{\delta_{i''}}$, which doesn't affect the
value at either $P_{i''}$ or $Q_{i''}$, and only affects the coefficient of 
the $(\delta_{i''},\delta_{i''})$ row, which in this case is supported purely
in the $i''$ column.\footnote{In this situation, varying $s^{i''}_j$
doesn't even change the limit linear series, but insofar as we made a choice
in our setup, we have to consider its possible effects.}
On the other hand,
if $a \neq a'-1$, observe that since the degree is $2$ in this
column, we cannot have any other sections involving $\delta_{i''}$ starting
or ending in the column,
and therefore we have no sections starting or
ending in the column. Thus, there are at most three potentially appearing
sections in column $i''$, and the other ones can't involve any exceptional 
row,
and must therefore be linearly independent. It follows that in our linear
dependence, the coefficient of $s_{(j,\delta_{i''})}$ must be nonzero.
Now, varying $s_j$ will change the relationship
between the values at $P_{i''}$ and $Q_{i''}$ (we can view the moduli for
$s_j$ as adding multiples of a section vanishing to order $a+1$ at $P_{i''}$ 
and order $b$ at $Q_{i''}$). Since this variation of moduli affects only a
single potentially appearing section, and we know it must have nonzero
coefficient in our linear dependence, there is only one choice of $s^{i''}_j$
compatible with the previously determined linear dependence, and we have
no nontrivial moduli in this case.

Finally, note that although our normalization was not determined for
potentially appearing sections supported in a single column, scaling these
does not affect the coefficients of any of the sections spanning the $i'-1$
and $i'$ column, so we have that the possible coefficients of these sections
are determined up to finitely many possibilities.
It thus suffices to show that if we vary
the gluing points on the component corresponding to the final column,
the (unique, if it exists) linear independence on the three potentially
appearing sections varies nontrivially. 

Now, necessarily the last column has the same value $a$ in all three rows in
its left subcolumn. On the right subcolumn, the criticality condition
rules out that there is a unique minimum value, although if there were the 
situation would be even simpler, since we could just drop the potentially
appearing sections in this column right away.
If $b$ is the minimal value for the right subcolumn, we similarly see that
we must have $a+b=2d-2$, or we could not have three (or even two) remaining 
potentially appearing sections.
Thus, the only two cases to consider are that
$b$ is attained twice, or in
all three rows. The last condition in the
definition of criticality implies that none of the $(a,b)$ rows are
obtained by adding the $\delta_{i'}$ row to an exceptional row. 
Now, if all three rows are $(a,b)$ rows, we can directly apply Proposition
\ref{prop:nondegenerate} to conclude that the linear dependence in the 
$i'$th column varies
nontrivially with $P_{i'},Q_{i'}$, as desired. On the other hand, if two
rows are $(a,b)$ rows, we again apply Proposition \ref{prop:nondegenerate} 
to these two rows, and since we have
normalized all three rows so that the values at $P_{i'}$ agree, we
again see that the linear dependence among the three has to vary
nontrivially with $P_{i'},Q_{i'}$, as desired.
\end{proof}

\section{The $r=6$ case}\label{sec:r6}

We now specialize to $r=6$, and suppose we have $g=21+\epsilon$ and 
$d=24+\epsilon$ for some $\epsilon \geq 0$, so that $\rho=\epsilon$.
Then our multidegree has total degree $2d=2g+6$, so it is determined by
placing $3$s in six columns, and $2$s in the rest. 

Although it turns out we will have flexibility in which multidegree to
consider, for the purposes of classifying cases, it is helpful to
introduce the following.

\begin{defn}\label{defn:default-multidegree} Given a limit linear
series, the \textbf{default multidegree} $\omega_{\df}$ is determined by
placing a $3$:
\begin{enumerate}
\item in the first column; 
\item in the first column with 
$\bar{\lambda}_i^1+\bar{\lambda}_i^2=5$;
\item in the first column with 
$\bar{\lambda}_i^1+\bar{\lambda}_i^3=7$; 
\item in the column immediately after the last column with 
$\bar{\lambda}_i^1+\bar{\lambda}_i^3=7$;
\item in the column immediately after the last column with
$\bar{\lambda}_i^2+\bar{\lambda}_i^3=9$; 
\item in the last column.
\end{enumerate}
\end{defn}

Note that $\bar{\lambda}_i^{\ell}$ can only increase in a genus-$1$
column, so the default multidegree is unimaginative.

\begin{prop}\label{prop:crossing-bd} Fix an unimaginative multidegree.
Then for a column $i$, there can be at most three rows spanning columns $i$
and $i+1$ except in the following circumstances:
\begin{ilist}
\itm $\gamma_i = 0$ and $\bar{\lambda}_i^1+\bar{\lambda}_i^3\geq 8$;
\itm $\gamma_i = 2$ and $\bar{\lambda}_i^1+\bar{\lambda}_i^3\geq 7$;
\itm $\gamma_i = 4$ and $\bar{\lambda}_i^1+\bar{\lambda}_i^3\leq 7$;
\itm $\gamma_i = 6$ and $\bar{\lambda}_i^1+\bar{\lambda}_i^3\leq 6$.
\end{ilist}

In particular, in the default multidegree there are never more than three
rows spanning a given pair of columns.
\end{prop}

\begin{proof} We will use the criterion from Proposition 
\ref{prop:when-appears}; since this only involves the values of 
$j-\lambda_{i,j}=a^{i+1}_j-g(i)$, the general case reduces immediately
to the notationally simpler situation that $\bar{\lambda}_i=\lambda_i$ for 
all $i$. We thus assume that we are in this situation.
Then, because the sequence $j-\lambda_{i,j}$ is strictly increasing
in $j$, we see that pairs $(j_1,j_2)$ satisfying the identity for appearing
in the $i$th and $(i+1)$st columns from Proposition \ref{prop:when-appears} 
must be strictly nested, so we can have at most $r/2+1=4$ of them, and
we can only have all of these
if $\lambda_{i,j}+\lambda_{i,r-j}$ is constant for all $j$, so that in 
particular
$$2\lambda_{i,r/2}
=\lambda_{i,0}+\lambda_{i,r}.$$
We also see that we have to have
$$\gamma_i=r-2\lambda_{i,r/2}=6-2\lambda_{i,r/2},$$
(so in particular $\gamma_i$ has to be even) and more generally for 
$j=0,\dots,r$ we have
$$\gamma_i=r-\lambda_{i,j}-\lambda_{i,r-j}.$$
Summing, we find that 
$$\sum_{j=0}^r \lambda_{i,j} 
=\frac{(r+1)(r-\gamma_i)}{2}=7(3-\frac{\gamma_i}{2}),$$
so $\lambda_{i,r/2}=3-\frac{\gamma_i}{2}$.

If $\gamma_i=0$ we must have 
$\lambda_{i,r/2}=3$ and we conclude that we would have to have
$\bar{\lambda}_i^1+\bar{\lambda}_i^3\geq 8$.

If $\gamma_i=2$, we need $\lambda_{i,r/2}=2$. Let $n$ be the number of 
values of $j$ with $\lambda_{i,j} \leq 0$; then we must have 
$\lambda_{i,r-j} \geq 4$ for the same $n$ values of $j$, so
$\lambda_i^1+\lambda_i^3\geq (r+1-n)+n=7$, as desired.

Similarly, if $\gamma_i=4$ then $\lambda_{i,r/2}=1$, so if we have $n$
values of $j$ with $\lambda_{i,j} \geq 3$, then we also have
$\lambda_{i,r-j} \leq -1$, so as before we find
$\lambda_i^1+\lambda_i^3\leq (r+1-n)+n = 7$.

Finally, if $\gamma_i=6$ we have
$\lambda_{i,r/2}=0$ and we conclude that we would have to have
$\bar{\lambda}_i^1+\bar{\lambda}_i^3\leq 6$, as claimed.
\end{proof}

We can now prove the following theorem, which will in particular prove
the desired maximal rank statement in all sufficiently nondegenerate
cases for all $\epsilon$ in our family of cases. It will also suffice
to prove the genus-$22$ case of our main theorem.

\begin{thm}\label{thm:basic-r6} In the default multidegree, we can
always drop all potentially appearing sections using the rules from
Lemma \ref{lem:independence-crit}, so the potentially appearing
sections are all linearly independent.
\end{thm}

\begin{proof}
In the first column, unless there is a swap with $\delta_1=1$ we will have 
at most the rows $(0,0)$, $(0,1)$ and $(0,2)$ among the potentially 
appearing sections, while if there is a swap with $\delta_1=1$ we will have 
at most the rows $(0,1)$ and $(1,1)$ potentially appearing.
In either case, these must all have distinct orders of vanishing, so can 
all be dropped.
According to Corollary \ref{cor:2-case}, we will have at most one new 
row with a potentially appearing section in each column until we get to the 
next column of degree $3$, so these can all be dropped. 

Now, suppose that $i$ is minimal such that
$\bar{\lambda}_{i}^1+\bar{\lambda}_{i}^2=5$. Then we are looking at
$\ell_1=1$ and $\ell_2=2$, so
$\gamma_{i}-1=1=5-1-\ell_1-\ell_2$, and according to Corollary 
\ref{cor:3-case} and Lemma \ref{lem:impossible}, we have no
potentially appearing sections supported entirely in the $i$th column. 
Any other new potentially appearing sections
would have to be supported in the $i$th and $(i+1)$st columns, so
by Proposition \ref{prop:crossing-bd}, we have at most three of these.
Note that if we choose $i'$ minimal so that 
$\bar{\lambda}_{i'}^1+\bar{\lambda}_{i'}^2=6$, then in this case
$\gamma_{i'}=2=6-1-\ell_1-\ell_2$, so according to Corollary 
\ref{cor:2-case} and Lemma \ref{lem:impossible}, there is no
row starting in the $i'$th column. 
We then see that the $i$th
(respectively, $i'$th) columns are critical: if $a,b$ are the minimum values 
in the subcolumns, they have to add to at least $2d-2$ or the rows would 
not be potentially starting in the $i$th column (respectively, potentially
supported in the $i'$th column). 
The last condition of semicriticality and the condition for criticality then 
follow from the second and first parts of Lemma \ref{lem:impossible}, 
respectively.
It follows that the hypotheses of Lemma \ref{lem:independence-crit} (iii)
are satisfied, so we can drop all rows occurring in this block.
We can then again handle any additional columns
before the next degree-$3$ one.

The setup being symmetric, we can also go from right
to left in the same manner, eliminating all potentially appearing sections
occurring in any
columns outside the middle two degree-$3$ columns. For these columns,
we are considering $\ell_1=1$ and $\ell_2=3$, so we have 
$\gamma_{i+1}-1=2=7-1-\ell_1-\ell_2$ and
$\gamma_{i+1}-1=3=8-1-\ell_1-\ell_2$ respectively,
and according to Corollary \ref{cor:3-case} and Lemma 
\ref{lem:impossible}, neither column has any potentially appearing section
supported entirely in it. As before, we find we must 
have a block satisfying the hypotheses of Lemma \ref{lem:independence-crit}
(iii), which we can then eliminate.
\end{proof}

If the specialization of our linear series contains the ``expected'' 
sections $s_j$ for
every $j=0,\dots,r$ in the expected multidegrees $\omega_j$ (as in Proposition
\ref{prop:simple}), then Theorem \ref{thm:basic-r6} implies that the 
images of each $s_j \otimes s_{j'}$ in the default multidegree are
linearly independent, so the multiplication map has the desired rank 
$\binom{r+2}{2}=28$. However, some linear series may have more degenerate
specializations, and the remainder of the paper will be devoted to
applying Theorem \ref{thm:basic-r6} (and variants thereof) to handle these 
situations as well.
For this, the statement in terms of potentially appearing sections
(as opposed to the separate rows considered in \cite{o-l-t-z}) is 
crucial. In interesting cases, we can have strictly more than 
$28$ potentially appearing sections. This does not contradict the fact that 
we know the multiplication map can have rank at most $28$, because these do 
not occur separately in the linked linear series coming as the specialization
of any fixed family of linear series on the smooth fibers. In most limits, 
for every $(j_1,j_2)$ we will have a unique linear combination of the 
potentially appearing sections in the $(j_1,j_2)$ row which actually
arise in the specialization. What makes the degenerate cases more 
interesting is that in these cases, we may have more than one linear 
combination occurring from a given row, 
precisely in situations where the specialization fails to contain any
potentially appearing sections from some other row -- see Example
\ref{ex:g22-degen} below.

\begin{ex}\label{ex:g22-simple} We continue with the running example of
Example \ref{ex:g22-grd} in Table \ref{table:g22-ex}.
\end{ex}

\begin{sidewaystable}
\begin{center}
\vspace{13cm}
\resizebox{\textwidth}{!}{
\begin{tabular}{llr|lr|lr|lr|lr|lr|lr|lr|lr|lr|lr|lr|lr|lr|lr|lr|lr|lr|lr|lr|lr|lr}
 & & $47$ & $3$ & $45$ & $5$ & $43$ & $7$ & $41$ & $9$ & $38$ & $12$ & $36$ 
& $14$ & $33$ & $17$ & $31$ & $19$ & $29$ & $21$ & $27$ & $23$ & $25$ & $25$ 
& $23$ & $27$ & $21$ & $29$ & $19$ & $31$ & $17$ & $33$ & $14$ & $36$ & $12$ 
& $38$ & $9$ & $41$ & $7$ & $43$ & $5$ & $45$ & $3$ & $47$ & \\
\hline
$(0,0)$ & \cellcolor[gray]{.8} $0$ & \cellcolor[gray]{.8} $50$ & $0$ & $48$ 
& $2$ & $46$ & $4$ & $44$ & $6$ & $42$ & $8$ & $40$ & $10$ & $38$ & $12$ 
& $38$ & $12$ & $36$ & $14$ & $34$ & $16$ & $32$ & $18$ & $30$ & $20$ & $28$ 
& $22$ & $26$ & $24$ & $24$ & $26$ & $24$ & $26$ & $22$ & $28$ & $20$ & $30$ 
& $18$ & $32$ & $16$ & $34$ & $14$ & $36$ & $12$ \\
$(0,1)$ & \cellcolor[gray]{.8} $1$ & \cellcolor[gray]{.8} $48$ & $2$ & $47$ 
& $3$ & $45$ & $5$ & $43$ & $7$ & $41$ & $9$ & $39$ & $11$ & $37$ & $13$ 
& $36$ & $14$ & $34$ & $16$ & $33$ & $17$ & $31$ & $19$ & $29$ & $21$ & $27$ 
& $23$ & $25$ & $25$ & $23$ & $27$ & $22$ & $28$ & $21$ & $29$ & $19$ & $31$ 
& $17$ & $33$ & $15$ & $35$ & $13$ & $37$ & $11$ \\
$(0,2)$ & \cellcolor[gray]{.8} $2$ & \cellcolor[gray]{.8} $47$ 
& \cellcolor[gray]{.8} $3$ & \cellcolor[gray]{.8} $45$ 
& \cellcolor[gray]{.8} $5$ & \cellcolor[gray]{.8} $44$ & $6$ & $42$ & $8$ 
& $40$ & $10$ & $38$ & $12$ & $36$ & $14$ & $35$ & $15$ & $32$ & $18$ & $30$ 
& $20$ & $28$ & $22$ & $27$ & $23$ & $25$ & $25$ & $23$ & $27$ & $21$ & $29$ 
& $20$ & $30$ & $18$ & $32$ & $16$ & $34$ & $15$ & $35$ & $13$ & $37$ & $11$ 
& $39$ & $9$ \\
$(1,1)$ & $2$ & $46$ & \cellcolor[gray]{.8} $4$ & \cellcolor[gray]{.8} $46$ 
& $4$ & $44$ & $6$ & $42$ & $8$ & $40$ & $10$ & $38$ & $12$ & $36$ & $14$ 
& $34$ & $16$ & $32$ & $18$ & $32$ & $18$ & $30$ & $20$ & $28$ & $22$ & $26$ 
& $24$ & $24$ & $26$ & $22$ & $28$ & $20$ & $30$ & $20$ & $30$ & $18$ & $32$ 
& $16$ & $34$ & $14$ & $36$ & $12$ & $38$ & $10$ \\
$(0,3)$ & $3$ & $46$ & $4$ & $44$ & $6$ & $42$ & \cellcolor[gray]{.8} $8$ 
& \cellcolor[gray]{.8} $41$ & \cellcolor[gray]{.8} $9$ 
& \cellcolor[gray]{.8} $39$ & $11$ & $37$ & $13$ & $35$ & $15$ & $34$ & $16$ 
& $33$ & $17$ & $31$ & $19$ & $30$ & $20$ & $28$ & $22$ & $26$ & $24$ & $24$ 
& $26$ & $22$ & $28$ & $21$ & $29$ & $19$ & $31$ & $18$ & $32$ & $16$ & $34$ 
& $14$ & $36$ & $12$ & $38$ & $10$ \\
$(1,2)$ & $3$ & $45$ & $5$ & $44$ & \cellcolor[gray]{.8} $6$ 
& \cellcolor[gray]{.8} $43$ & \cellcolor[gray]{.8} $7$ 
& \cellcolor[gray]{.8} $41$ & \cellcolor[gray]{.8} $9$ 
& \cellcolor[gray]{.8} $39$ & $11$ & $37$ & $13$ & $35$ & $15$ & $33$ & $17$ 
& $30$ & $20$ & $29$ & $21$ & $27$ & $23$ & $26$ & $24$ & $24$ & $26$ & $22$ 
& $28$ & $20$ & $30$ & $18$ & $32$ & $17$ & $33$ & $15$ & $35$ & $14$ & $36$ 
& $12$ & $38$ & $10$ & $40$ & $8$ \\
$(0,4)$ & $4$ & $45$ & $5$ & $43$ & $7$ & $41$ & $9$ & $39$ 
& \cellcolor[gray]{.8} $11$ & \cellcolor[gray]{.8} $38$ 
& \cellcolor[gray]{.8} $12$ & \cellcolor[gray]{.8} $36$ 
& \cellcolor[gray]{.8} $14$ & \cellcolor[gray]{.8} $34$ & $16$ & $33$ & $17$ 
& $31$ & $19$ & $29$ & $21$ & $27$ & $23$ & $25$ & $25$ & $24$ & $26$ & $22$ 
& $28$ & $20$ & $30$ & $19$ & $31$ & $17$ & $33$ & $15$ & $35$ & $13$ & $37$ 
& $12$ & $38$ & $10$ & $40$ & $8$ \\
$(1,3)$ & $4$ & $44$ & $6$ & $43$ & $7$ & $41$ & $9$ & $40$ 
& \cellcolor[gray]{.8} $10$ & \cellcolor[gray]{.8} $38$ 
& \cellcolor[gray]{.8} $12$ & \cellcolor[gray]{.8} $36$ 
& \cellcolor[gray]{.8} $14$ & \cellcolor[gray]{.8} $34$ & $16$ & $32$ & $18$ 
& $31$ & $19$ & $30$ & $20$ & $29$ & $21$ & $27$ & $23$ & $25$ & $25$ & $23$ 
& $27$ & $21$ & $29$ & $19$ & $31$ & $18$ & $32$ & $17$ & $33$ & $15$ & $35$ 
& $13$ & $37$ & $11$ & $39$ & $9$ \\
$(2,2)$ & $4$ & $44$ & $6$ & $42$ & $8$ & $42$ & $8$ & $40$ 
& \cellcolor[gray]{.8} $10$ & \cellcolor[gray]{.8} $38$ 
& \cellcolor[gray]{.8} $12$ & \cellcolor[gray]{.8} $36$ 
& \cellcolor[gray]{.8} $14$ & \cellcolor[gray]{.8} $34$ & $16$ & $32$ & $18$ 
& $28$ & $22$ & $26$ & $24$ & $24$ & \cellcolor[gray]{.8} $26$ 
& \cellcolor[gray]{.8} $24$ & $26$ & $22$ & $28$ & $20$ & $30$ & $18$ & $32$ 
& $16$ & $34$ & $14$ & $36$ & $12$ & $38$ & $12$ & $38$ & $10$ & $40$ & $8$ 
& $42$ & $6$ \\
$(0,5)$ & $5$ & $44$ & $6$ & $42$ & $8$ & $40$ & $10$ & $38$ & $12$ & $36$ 
& $14$ & $35$ & \cellcolor[gray]{.8} $15$ & \cellcolor[gray]{.8} $33$ 
& \cellcolor[gray]{.8} $17$ & \cellcolor[gray]{.8} $32$ & $18$ & $30$ & $20$ 
& $28$ & $22$ & $26$ & $24$ & $24$ & $26$ & $22$ & $28$ & $21$ & $29$ & $19$ 
& $31$ & $18$ & $32$ & $16$ & $34$ & $14$ & $36$ & $12$ & $38$ & $10$ & $40$ 
& $9$ & $41$ & $7$ \\
$(1,4)$ & $5$ & $43$ & $7$ & $42$ & $8$ & $40$ & $10$ & $38$ & $12$ & $37$ 
& $13$ & $35$ & \cellcolor[gray]{.8} $15$ & \cellcolor[gray]{.8} $33$ 
& \cellcolor[gray]{.8} $17$ & \cellcolor[gray]{.8} $31$ 
& \cellcolor[gray]{.8} $19$ & \cellcolor[gray]{.8} $29$ 
& \cellcolor[gray]{.8} $21$ & \cellcolor[gray]{.8} $28$ & $22$ & $26$ & $24$ 
& $24$ & $26$ & $23$ & $27$ & $21$ & $29$ & $19$ & $31$ & $17$ & $33$ & $16$ 
& $34$ & $14$ & $36$ & $12$ & $38$ & $11$ & $39$ & $9$ & $41$ & $7$ \\
$(2,3)$ & $5$ & $43$ & $7$ & $41$ & $9$ & $40$ & $10$ & $39$ & $11$ & $37$ 
& $13$ & $35$ & \cellcolor[gray]{.8} $15$ & \cellcolor[gray]{.8} $33$ 
& \cellcolor[gray]{.8} $17$ & \cellcolor[gray]{.8} $31$ 
& \cellcolor[gray]{.8} $19$ & \cellcolor[gray]{.8} $29$ 
& \cellcolor[gray]{.8} $21$ & \cellcolor[gray]{.8} $27$ 
& \cellcolor[gray]{.8} $23$ & \cellcolor[gray]{.8} $26$ & $24$ & $25$ & $25$ 
& $23$ & $27$ & $21$ & $29$ & $19$ & $31$ & $17$ & $33$ & $15$ & $35$ & $14$ 
& $36$ & $13$ & $37$ & $11$ & $39$ & $9$ & $41$ & $7$ \\
$(0,6)$ & $6$ & $43$ & $7$ & $41$ & $9$ & $39$ & $11$ & $37$ & $13$ & $35$ 
& $15$ & $33$ & $17$ & $32$ & \cellcolor[gray]{.8} $18$ 
& \cellcolor[gray]{.8} $31$ & \cellcolor[gray]{.8} $19$ 
& \cellcolor[gray]{.8} $29$ & \cellcolor[gray]{.8} $21$ 
& \cellcolor[gray]{.8} $27$ & \cellcolor[gray]{.8} $23$ 
& \cellcolor[gray]{.8} $25$ & \cellcolor[gray]{.8} $25$ 
& \cellcolor[gray]{.8} $23$ & \cellcolor[gray]{.8} $27$ 
& \cellcolor[gray]{.8} $21$ & \cellcolor[gray]{.8} $29$ 
& \cellcolor[gray]{.8} $19$ & \cellcolor[gray]{.8} $31$ 
& \cellcolor[gray]{.8} $18$ & $32$ & $17$ & $33$ & $15$ & $35$ & $13$ & $37$ 
& $11$ & $39$ & $9$ & $41$ & $7$ & $43$ & $6$ \\
$(1,5)$ & $6$ & $42$ & $8$ & $41$ & $9$ & $39$ & $11$ & $37$ & $13$ & $35$ 
& $15$ & $34$ & $16$ & $32$ & $18$ & $30$ & $20$ & $28$ 
& \cellcolor[gray]{.8} $22$ & \cellcolor[gray]{.8} $27$ 
& \cellcolor[gray]{.8} $23$ & \cellcolor[gray]{.8} $25$ 
& \cellcolor[gray]{.8} $25$ & \cellcolor[gray]{.8} $23$ 
& \cellcolor[gray]{.8} $27$ & \cellcolor[gray]{.8} $21$ 
& \cellcolor[gray]{.8} $29$ & \cellcolor[gray]{.8} $20$ & $30$ & $18$ & $32$ 
& $16$ & $34$ & $15$ & $35$ & $13$ & $37$ & $11$ & $39$ & $9$ & $41$ & $8$ 
& $42$ & $6$ \\
$(2,4)$ & $6$ & $42$ & $8$ & $40$ & $10$ & $39$ & $11$ & $37$ & $13$ & $36$ 
& $14$ & $34$ & $16$ & $32$ & $18$ & $30$ & $20$ & $27$ & $23$ & $25$ & $25$ 
& $23$ & $27$ & $22$ & \cellcolor[gray]{.8} $28$ & \cellcolor[gray]{.8} $21$ 
& \cellcolor[gray]{.8} $29$ & \cellcolor[gray]{.8} $19$ 
& \cellcolor[gray]{.8} $31$ & \cellcolor[gray]{.8} $17$ 
& \cellcolor[gray]{.8} $33$ & \cellcolor[gray]{.8} $15$ & $35$ & $13$
& $37$ & $11$ & $39$ & $10$ & $40$ & $9$ & $41$ & $7$ & $43$ & $5$ \\
$(3,3)$ & $6$ & $42$ & $8$ & $40$ & $10$ & $38$ & $12$ & $38$ & $12$ & $36$ 
& $14$ & $34$ & $16$ & $32$ & $18$ & $30$ & \cellcolor[gray]{.8} $20$ 
& \cellcolor[gray]{.8} $30$ & $20$ & $28$ & $22$ & $28$ & $22$ & $26$ & $24$ 
& $24$ & $26$ & $22$ & $28$ & $20$ & $30$ & $18$ & $32$ & $16$ & $34$ & $16$ 
& $34$ & $14$ & $36$ & $12$ & $38$ & $10$ & $40$ & $8$ \\
$(1,6)$ & $7$ & $41$ & $9$ & $40$ & $10$ & $38$ & $12$ & $36$ & $14$ & $34$ 
& $16$ & $32$ & $18$ & $31$ & $19$ & $29$ & $21$ & $27$ & $23$ & $26$ & $24$ 
& $24$ & $26$ & $22$ & $28$ & $20$ & $30$ & $18$ & \cellcolor[gray]{.8} $32$ 
& \cellcolor[gray]{.8} $17$ & \cellcolor[gray]{.8} $33$ 
& \cellcolor[gray]{.8} $15$ & $35$ & $14$ & $36$ & $12$ & $38$ & $10$ & $40$ 
& $8$ & $42$ & $6$ & $44$ & $5$ \\
$(2,5)$ & $7$ & $41$ & $9$ & $39$ & $11$ & $38$ & $12$ & $36$ & $14$ & $34$ 
& $16$ & $33$ & $17$ & $31$ & $19$ & $29$ & $21$ & $26$ & $24$ & $24$ & $26$ 
& $22$ & $28$ & $21$ & $29$ & $19$ & $31$ & $18$ & $32$ & $16$ 
& \cellcolor[gray]{.8} $34$ & \cellcolor[gray]{.8} $14$ 
& \cellcolor[gray]{.8} $36$ & \cellcolor[gray]{.8} $12$ 
& \cellcolor[gray]{.8} $38$ & \cellcolor[gray]{.8} $10$ & $40$ & $9$ & $41$ 
& $7$ & $43$ & $6$ & $44$ & $4$ \\
$(3,4)$ & $7$ & $41$ & $9$ & $39$ & $11$ & $37$ & $13$ & $36$ & $14$ & $35$ 
& $15$ & $33$ & $17$ & $31$ & $19$ & $29$ & $21$ & $28$ & $22$ & $26$ 
& \cellcolor[gray]{.8} $24$ & \cellcolor[gray]{.8} $25$ 
& \cellcolor[gray]{.8} $25$ & \cellcolor[gray]{.8} $23$ 
& \cellcolor[gray]{.8} $27$ & \cellcolor[gray]{.8} $22$ & $28$ & $20$ & $30$ 
& $18$ & $32$ & $16$ & $34$ & $14$ & $36$ & $13$ & $37$ & $11$ & $39$ & $10$ 
& $40$ & $8$ & $42$ & $6$ \\
$(2,6)$ & $8$ & $40$ & $10$ & $38$ & $12$ & $37$ & $13$ & $35$ & $15$ & $33$ 
& $17$ & $31$ & $19$ & $30$ & $20$ & $28$ & $22$ & $25$ & $25$ & $23$ & $27$ 
& $21$ & $29$ & $20$ & $30$ & $18$ & $32$ & $16$ & $34$ & $15$ & $35$ & $13$ 
& $37$ & $11$ & \cellcolor[gray]{.8} $39$ & \cellcolor[gray]{.8} $9$ 
& \cellcolor[gray]{.8} $41$ & \cellcolor[gray]{.8} $8$ & $42$ & $6$ & $44$ 
& $4$ & $46$ & $3$ \\
$(3,5)$ & $8$ & $40$ & $10$ & $38$ & $12$ & $36$ & $14$ & $35$ & $15$ & $33$ 
& $17$ & $32$ & $18$ & $30$ & $20$ & $28$ & $22$ & $27$ & $23$ & $25$ & $25$ 
& $24$ & $26$ & $22$ & $28$ & $20$ & \cellcolor[gray]{.8} $30$ 
& \cellcolor[gray]{.8} $19$ & \cellcolor[gray]{.8} $31$ 
& \cellcolor[gray]{.8} $17$ & \cellcolor[gray]{.8} $33$ 
& \cellcolor[gray]{.8} $15$ & $35$ & $13$ & $37$ & $12$ & $38$ & $10$ & $40$ 
& $8$ & $42$ & $7$ & $43$ & $5$ \\
$(4,4)$ & $8$ & $40$ & $10$ & $38$ & $12$ & $36$ & $14$ & $34$ & $16$ & $34$ 
& $16$ & $32$ & $18$ & $30$ & $20$ & $28$ & $22$ & $26$ & $24$ & $24$ & $26$ 
& $22$ & $28$ & $20$ & $30$ & $20$ & $30$ & $18$ & $32$ & $16$ 
& \cellcolor[gray]{.8} $34$ & \cellcolor[gray]{.8} $14$ 
& \cellcolor[gray]{.8} $36$ & \cellcolor[gray]{.8} $12$
& \cellcolor[gray]{.8} $38$ & \cellcolor[gray]{.8} $10$ & $40$ & $8$ & $42$ 
& $8$ & $42$ & $6$ & $44$ & $4$ \\
$(3,6)$ & $9$ & $39$ & $11$ & $37$ & $13$ & $35$ & $15$ & $34$ & $16$ & $32$ 
& $18$ & $30$ & $20$ & $29$ & $21$ & $27$ & $23$ & $26$ & $24$ & $24$ & $26$ 
& $23$ & $27$ & $21$ & $29$ & $19$ & $31$ & $17$ & $33$ & $16$ 
& \cellcolor[gray]{.8} $34$ & \cellcolor[gray]{.8} $14$ 
& \cellcolor[gray]{.8} $36$ & \cellcolor[gray]{.8} $12$ 
& \cellcolor[gray]{.8} $38$ & \cellcolor[gray]{.8} $11$ & $39$ & $9$ & $41$ 
& $7$ & $43$ & $5$ & $45$ & $4$ \\
$(4,5)$ & $9$ & $39$ & $11$ & $37$ & $13$ & $35$ & $15$ & $33$ & $17$ & $32$ 
& $18$ & $31$ & $19$ & $29$ & $21$ & $27$ & $23$ & $25$ & $25$ & $23$ & $27$ 
& $21$ & $29$ & $19$ & $31$ & $18$ & $32$ & $17$ & $33$ & $15$ & $35$ & $13$ 
& $37$ & $11$ & \cellcolor[gray]{.8} $39$ & \cellcolor[gray]{.8} $9$ 
& \cellcolor[gray]{.8} $41$ & \cellcolor[gray]{.8} $7$ 
& \cellcolor[gray]{.8} $43$ & \cellcolor[gray]{.8} $6$ & $44$ & $5$ & $45$ 
& $3$ \\
$(4,6)$ & $10$ & $38$ & $12$ & $36$ & $14$ & $34$ & $16$ & $32$ & $18$ 
& $31$ & $19$ & $29$ & $21$ & $28$ & $22$ & $26$ & $24$ & $24$ & $26$ & $22$ 
& $28$ & $20$ & $30$ & $18$ & $32$ & $17$ & $33$ & $15$ & $35$ & $14$ & $36$ 
& $12$ & $38$ & $10$ & $40$ & $8$ & $42$ & $6$ & \cellcolor[gray]{.8} $44$ 
& \cellcolor[gray]{.8} $5$ & \cellcolor[gray]{.8} $45$ 
& \cellcolor[gray]{.8} $3$ & \cellcolor[gray]{.8} $47$ 
& \cellcolor[gray]{.8} $2$ \\
$(5,5)$ & $10$ & $38$ & $12$ & $36$ & $14$ & $34$ & $16$ & $32$ & $18$ 
& $30$ & $20$ & $30$ & $20$ & $28$ & $22$ & $26$ & $24$ & $24$ & $26$ & $22$ 
& $28$ & $20$ & $30$ & $18$ & $32$ & $16$ & $34$ & $16$ & $34$ & $14$ & $36$ 
& $12$ & $38$ & $10$ & $40$ & $8$ & $42$ & $6$ & $44$ & $4$ 
& \cellcolor[gray]{.8} $46$ & \cellcolor[gray]{.8} $4$ & $46$ & $2$ \\
$(5,6)$ & $11$ & $37$ & $13$ & $35$ & $15$ & $33$ & $17$ & $31$ & $19$ 
& $29$ & $21$ & $28$ & $22$ & $27$ & $23$ & $25$ & $25$ & $23$ & $27$ & $21$ 
& $29$ & $19$ & $31$ & $17$ & $33$ & $15$ & $35$ & $14$ & $36$ & $13$ & $37$ 
& $11$ & $39$ & $9$ & $41$ & $7$ & $43$ & $5$ & $45$ & $3$ & $47$ & $2$ 
& \cellcolor[gray]{.8} $48$ & \cellcolor[gray]{.8} $1$ \\
$(6,6)$ & $12$ & $36$ & $14$ & $34$ & $16$ & $32$ & $18$ & $30$ & $20$ 
& $28$ & $22$ & $26$ & $24$ & $26$ & $24$ & $24$ & $26$ & $22$ & $28$ & $20$ 
& $30$ & $18$ & $32$ & $16$ & $34$ & $14$ & $36$ & $12$ & $38$ & $12$ & $38$ 
& $10$ & $40$ & $8$ & $42$ & $6$ & $44$ & $4$ & $46$ & $2$ & $48$ & $0$ 
& \cellcolor[gray]{.8} $50$ & \cellcolor[gray]{.8} $0$ \\
\hline
 & & $47$ & $3$ & $45$ & $5$ & $43$ & $7$ & $41$ & $9$ & $38$ & $12$ & $36$ 
& $14$ & $33$ & $17$ & $31$ & $19$ & $29$ & $21$ & $27$ & $23$ & $25$ & $25$ 
& $23$ & $27$ & $21$ & $29$ & $19$ & $31$ & $17$ & $33$ & $14$ & $36$ & $12$ 
& $38$ & $9$ & $41$ & $7$ & $43$ & $5$ & $45$ & $3$ & $47$ & \\
\end{tabular}
}
\end{center}

\caption{The above table is the $T$ obtained from the tensor square of the
limit linear series
considered in Example \ref{ex:g22-grd}, which has $r=6$, $g=22$, and $d=25$.
We have also included the $w$ corresponding to the default multidegree 
$\omega_{\df}$; for ease of reading, we place the
multidegree at both the top and bottom of the table, and include not only
the values $c_i$ for $i=2,\dots,22$, but also $2d-c_i$ in the preceding 
subcolumns. 
We have highlighted the potentially appearing sections; note
that the $(2,2)$ row contains two, while the rest all have a unique one.
These two potentially appearing sections are thus treated separately in
Theorem \ref{thm:basic-r6}; the first appears as part of a block in the
$5$th and $6$th columns which is eliminated using rule (iii) of Lemma
\ref{lem:independence-crit}, while the second occurs as the only new
potentially appearing sections appearing in the $12$th column, which
is part of another block, extending from the $7$th column to the $16$th
column, which is again eliminated using rule (iii), after all other
potentially appearing sections have been eliminated on both the left and
right. The only other block that requires rule (iii) contains the
$17$th and $18$th columns, and is eliminated after the potentially 
appearing sections appearing to the right have all been dropped. Following
the proof of Theorem \ref{thm:basic-r6}, we see that we can eliminate 
all sections outside the aforementioned three blocks going inward from both 
the left and right ends, using only iterated applications of rule (i).
\label{table:g22-ex}}
\end{sidewaystable}

Ultimately, the default multidegree used in Theorem \ref{thm:basic-r6} will 
be sufficient to handle the genus-$22$ case, and most of the genus-$23$
cases. However, for certain degenerate cases we will need to consider other 
multidegrees instead. 

We will thus want to develop the following results describing the flexibility
we have in choosing the multidegree while maintaining linear independence.

\begin{prop}\label{prop:multideg-flexible} Suppose we have an unimaginative
multidegree $\omega$ determined by placing degree $3$ in genus-$1$ columns as 
follows:
\begin{enumerate}
\item in one column which is either the first, or a column with no exceptional
rows and
satisfying
$\bar{\lambda}_i^1+\bar{\lambda}_i^2 \leq 4$ and 
$\bar{\lambda}_{i,0}  \leq 2$;
\item in one column with $\bar{\lambda}_i^1+\bar{\lambda}_i^2=5$ but
$\bar{\lambda}_{i-1}^1+\bar{\lambda}_{i-1}^2= 4$;
\item in one column between the first column with 
$\bar{\lambda}_i^1+\bar{\lambda}_i^2=6$ and the first column with
$\bar{\lambda}_i^1+\bar{\lambda}_i^3=7$ (inclusive); 
\item in one column between the column immediately after the last 
column with $\bar{\lambda}_i^1+\bar{\lambda}_i^3=7$ and the column
immediately after the last column with 
$\bar{\lambda}_i^2+\bar{\lambda}_i^3=8$ (inclusive);
\item in one column with $\bar{\lambda}_i^2+\bar{\lambda}_i^3=10$ but
$\bar{\lambda}_{i-1}^2+\bar{\lambda}_{i-1}^3=9$; 
\item in one column which is either the last, or a column with no exceptional
rows and satisfying
$\bar{\lambda}_{i-1}^2+\bar{\lambda}_{i-1}^3\geq 10$ and 
$\bar{\lambda}_{i-1,6} \geq 1$.
\end{enumerate}

Then the potentially appearing sections in multidegree $\omega$ are still
linearly independent.
\end{prop}

\begin{proof} The main new ingredient is verifying that if we place 
the first degree $3$ in a (genus-$1$) column after the first, but still
satisfying $\bar{\lambda}_i^1+\bar{\lambda}_i^2 \leq 4$ and
$\bar{\lambda}_{i,0} \leq 2$, then provided we also have no exceptional rows,
we will
in fact obtain at most two potentially appearing sections starting
in the $i$th column. Note that in this case, in particular there are no
swaps in the $i$th column. By Proposition \ref{prop:when-appears},
for the $(j,j')$ row to have a potentially appearing
section starting in the $i$th column, we will need 
$j+j'-\lambda_{i-1,j}-\lambda_{i-1,j'}>\gamma_{i-1}=0$ and
$j+j'-\lambda_{i,j}-\lambda_{i,j'}\leq \gamma_{i}=1$, or equivalently
\begin{equation}\label{eq:gamma1}
\lambda_{i-1,j}+\lambda_{i-1,j'}<
j+j' \leq 1+ \lambda_{i,j}+\lambda_{i,j'}.
\end{equation}
For this assertion, since we are assuming no swaps occur in the $i$th
column, it suffices to check the case with $\bar{\lambda}_i=\lambda_i$
for all $i$, which simplifies notation.
Now, since $\bar{\lambda}_i^1+\bar{\lambda}_i^2 \leq 4$, we must have
$\lambda_{i,j}\leq 0$ for $j \geq 4$ and 
$\lambda_{i,j}\leq 1$ for $j=2,3$. 
It follows that to satisfy
the righthand inequality above, we must have at least one of $j,j'$ 
equal to $0$ or $1$. Moreover, by Corollary \ref{cor:3-case} we have that 
for $j=0,1$, if $j \neq \delta_i$, then there is at most one value of $j'$ 
satisfying the above inequalities. 
In particular, we conclude that if 
$\delta_i\neq 0,1$, we have at most two potentially appearing sections,
as claimed. 

Now, if $\delta_i=0$, we need to see that we can have at most 
two rows of the form $(0,j')$ appearing in the $i$th column, and if
two appear, then none of the form $(1,j')$ can appear for $j'>0$. Suppose 
first $(0,0)$ is potentially starting in the $i$th column. By
\eqref{eq:gamma1} this could only happen if 
$\lambda_{i-1,0}<0$, so
$\lambda_{i,j'}=\lambda_{i-1,j'}<0$ for all $j'>0$, 
and then $(0,j')$ cannot satisfy the righthand side of \eqref{eq:gamma1}
for any $j'>0$. On the other hand, if 
we have $j''>j'>0$ such that $(0,j')$ and $(0,j'')$ both appear,
then we have
$$\lambda_{i-1,0}+\lambda_{i-1,j'}<
j' < j'' \leq 1+ \lambda_{i,0}+\lambda_{i,j''}
\leq 2+\lambda_{i-1,0}+\lambda_{i-1,j'},$$
so the only possibility is that 
$j'= 1+\lambda_{i-1,0}+\lambda_{i-1,j'}$ and
$j''=j'+1$, with 
$\lambda_{i,j''}= \lambda_{i,j'}$. It immediately
follows that 
we could not have $(0,j''')$ appearing for any $j'''\neq 0,j',j''$.
We also check that $(1,j''')$ cannot be potentially appearing for any 
$j''' >0$ in this situation. Indeed, $1+j'''$ will be too large if
$j''' \geq j'$. 
The $(1,1)$ row cannot satisfy \eqref{eq:gamma1} for parity reasons, and 
in order to have $(1,2)$ appearing we would need $j'' \geq 4$, but
we note that in this case
$1+ \lambda_{i,0}+\lambda_{i,j''}\leq 3$,
contradicting \eqref{eq:gamma1}.

Finally, consider the case that $\delta_i=1$. If the $(1,1)$ row is 
potentially starting in the $i$th column, by parity we have to have
$1=\lambda_{i,1}$, so for all $j>1$ we have
$\lambda_{i,j}=\lambda_{i-1,j} \leq 
\lambda_{i-1,1} =0$. Then we cannot have $(1,j')$ potentially 
starting for any $j'>1$,
so we have at most two rows potentially starting.
On the other hand, if we have $j''>j'>1$ potentially starting in the $i$th
column, we are just as above forced to have
$j'= \lambda_{i-1,1}+\lambda_{i-1,j'}$ and
$j''=j'+1$, with 
$\lambda_{i,j''}= \lambda_{i,j'}$, and we
claim we cannot have $(0,j''')$ potentially starting for any $j'''$. 
Indeed, if $j''' \leq j''$, then we have
$j''' - \lambda_{i-1,j'''} \leq 
j'' - \lambda_{i-1,j''}$,
so
\begin{align*}j''' & \leq j'' - \lambda_{i-1,j''}
+ \lambda_{i-1,j'''} \\ & = 
1+\lambda_{i-1,1}+\lambda_{i-1,j''}
- \lambda_{i-1,j''}
+ \lambda_{i-1,j'''} \\ & \leq
\lambda_{i-1,0} + \lambda_{i-1,j'''},
\end{align*}
violating \eqref{eq:gamma1}. But $j'' \geq 3$, so if $j''' \geq 4$
we cannot satisfy \eqref{eq:gamma1} without violating our hypothesis that
$\lambda_{i-1,0}\leq 2$. We thus conclude the desired statement
on the number of potentially appearing sections starting in column $i$.

Now, since we have assumed that our first column with degree $3$ has
no exceptional rows, the fact that it has at most two
potentially appearing sections starting in it means that we can still
eliminate sections from left to right until we reach the second column
of degree $3$, just as in the proof of Theorem \ref{thm:basic-r6}. We also
see that the second column of degree $3$ will still be critical,
with at most three potentially appearing sections starting in it.
The next step depends on the location of the third column of degree $3$.
If the first column with $\bar{\lambda}_i^1+\bar{\lambda}_i^2=6$
still has degree $2$, we will eliminate this block from left to right,
as before. On the other hand, if the first column with
$\bar{\lambda}_i^1+\bar{\lambda}_i^3=7$ has degree $2$, we do not need
to have eliminated everything from the left in order to eliminate the
central block, since the potentially supported rows in multidegree $\omega$
will be precisely the same as the potentially starting rows in $\omega_{\df}$. 
Thus, if the third column of degree $3$ is strictly between
these, we can eliminate both adjacent blocks first, and then eliminate
all potentially appearing sections one by one from both sides until we
reach this final column, which can have at most one remaining potentially
appearing section by Corollary \ref{cor:3-case}. 
However, if the third column of degree $3$ is the 
first column with $\bar{\lambda}_i^1+\bar{\lambda}_i^2=6$, we see that this
will be critical with at most three potentially appearing sections ending
in it,
and we will instead 
eliminate the central block first, and then eliminate the block between
the second and third columns of degree $3$ last.

The situation is symmetric on the right, so we see that in all cases we
will be able to eliminate all potentially appearing sections in a suitable
order.
\end{proof}

It will also be important to consider moving degree $3$ into a column with
a swap, which we analyze with the below lemma.

\begin{lem}\label{lem:3-at-swap} Suppose that in the $i$th column, 
we have a $\delta_i$, and for some $j_0$ we have
$a^i_{\delta_i}=a^i_{j_0}+1$,
but $a^{i+1}_{\delta_i}=a^{i+1}_{j_0}-1$. Suppose further that our 
multidegree assigns degree $3$ to the $i$th column, and write it as usual
as $\md_{2d}(w)$ for $w=(c_2,\dots,c_g)$.

Then we can have at most four potentially appearing
sections start in the $i$th column. Moreover, we can only have four
if either we have
\begin{equation}\label{eq:3-at-swap} 2\bar{a}^i_3=c_i+1,
\end{equation} 
or if 
$(\delta_i,\delta_i)$ is potentially starting, and one of the following 
three possibilities holds:
\begin{enumerate}
\item $a^i_{(\delta_i,\delta_i)}=c_i+1$;
\item $a^i_{(\delta_i,\delta_i)}=c_i+2$, and the $(j_0,\delta_i)$ row
is potentially starting;
\item $a^i_{(\delta_i,\delta_i)}=c_i+3$, with $a^i_{\delta_i}=\bar{a}^i_4$. 
\end{enumerate}

We can also have at most four potentially appearing sections end in the
$i$th column, with four of them ending only if either we have
\begin{equation}\label{eq:3-at-swap-2}
2\bar{a}^i_3=c_i,
\end{equation} 
or if $(\delta_i,\delta_i)$ is potentially ending, and
one the following three possibilities holds:
\begin{enumerate}
\item $a^i_{(\delta_i,\delta_i)}=c_i$, with $a^i_{\delta_i}=\bar{a}^i_3$; 
\item $a^i_{(\delta_i,\delta_i)}=c_i+1$, and the $(j_0,\delta_i)$ row
is potentially ending;
\item $a^i_{(\delta_i,\delta_i)}=c_i+2$; 
\end{enumerate}
\end{lem}

We have written the above to allow for swaps having occurred prior to the
$i$th column. If no swaps have occurred, the $j_0$ in the lemma statement is 
necessarily $\delta_i-1$, and the third exceptional case would require
$\delta_i=4$ (respectively, $\delta_i=2$) in the statement on potential
support starting (respectively, ending).

\begin{proof} In order to have a potentially appearing section start in the
$(j,j')$ row, we must have $a^i_{(j,j')}>c_i$ and 
$a^{i+1}_{(j,j')} \leq c_{i+1}=c_i+3$. It immediately follows that if
$j,j' \neq \delta_i$, then we must have $a^i_{(j,j')}=c_i+1$, and 
neither $j$ nor $j'$ equal to $j_0$. 
If $j'=\delta_i$ for $j \neq \delta_i$,
we could also have $a^i_{(j,\delta_i)}=c_i+2$, provided that $j \neq j_0$.
And $j_0$ can occur only if $a^i_{(j_0,\delta_i)}=c_i+1$, or equivalently if
$a^i_{(\delta_i,\delta_i)}=c_i+2$. Recall that Corollary \ref{cor:3-case}
says that if $(j,\delta_i)$ occurs for some $j \neq \delta_i$, then there is 
no $j' \neq \delta_i$ with $(j,j')$ also occurring. 
Next, we note that we can have at most two rows of
the form $(j,\delta_i)$ occurring. Indeed, if 
$a^i_{(\delta_i,\delta_i)} = c_i+3$, then 
$a^i_{(j,\delta_i)}=c_i+2$ only for $j=j_0$, so the $(j_0,\delta_i)$ row
does not occur, and we can have at most one additional row, having
$a^i_{(j,\delta_i)}=c_i+1$. On the other hand, if
$a^i_{(\delta_i,\delta_i)} \neq c_i+3$, then we have at most two rows, 
because they have to satisfy 
$a^i_{(j,\delta_i)}=c_i+1$ or $c_i+2$.
We also observe that we can have a row of the
form $(j,j)$ for $j \neq \delta_i$ only for a unique choice of $j$,
necessarily with $2a^i_j=c_i+1$ and $j \neq j_0$, and then we cannot have 
$(\delta_i,\delta_i)$ occurring, since $a^i_j \neq a^i_{\delta_i}-1$
for $j \neq j_0$. 

Now, if we do not have any $(j,\delta_i)$ occurring, then $j_0$ also cannot
occur, and we are left with only $5$ values of $j$ from which to choose
distinct pairs, so we can obtain at most three pairs (allowing one of them
to have repeated entries). Similarly, if we have exactly one $(j,\delta_i)$,
then necessarily $j\neq j_0$, or we would be in the 2nd exceptional case
with also $(\delta_i,\delta_i)$ occurring, so we have that for the remaining
pairs we must choose from values not equal to $\delta_i,j_0,j$, leaving
four values, and at most two pairs. We therefore see that in order to have
four rows potentially starting, two of them need to involve $\delta_i$.

Next, if we have $(j_1,\delta_i)$ and $(j_2,\delta_i)$ occurring, with 
neither $j_1,j_2$ equal to $\delta_i$ (and hence also neither equal to 
$j_0$), then any remaining rows have to be chosen as distinct pairs from 
the remaining $(r+1)-4=3$ indices, with at most one pair having repeated
value. We thus obtain at most four rows, with four occurring only if
$2a^i_j=a^i_{j_3}+a^i_{j_4}=c_i+1$ for some 
$j,j_3,j_4 \neq \delta_i,j_0,j_1,j_2$. Moreover, we see that
there must be exactly three values of $j'$ with $a^i_{j'}<a^i_j$ in 
this case: if $a^i_{\delta_i}<a^i_j$, then these are $\delta_i$, $j_0$,
and exactly one of $j_3,j_4$, with necessarily $j_1,j_2$ and the other
of $j_3,j_4$ having $a^i_{j'}>a^i_j$. If $a^i_{\delta_i}>a^i_j$, then
$a^i_{j_0}$ must also be greater than $a^i_j$, so we similarly find
exactly three values are smaller. Thus \eqref{eq:3-at-swap} must hold.

It remains to consider the case that $(\delta_i,\delta_i)$ row is potentially
starting, and the only thing left to prove is the description of case (3),
where $a^i_{(\delta_i,\delta_i)}=c_i+3$. Here, we must also have a $j$ with
$a^i_{(j,\delta_i)}=c_i+1$, and if we have two additional rows appearing,
these must come from two additional pairs nested around $a^i_{(j,\delta_i)}$ 
in value, so since $a^i_j<a^i_{j_0}<a^i_{\delta_i}$ in this case, we obtain
the desired statement.

The statement on rows ending is symmetric.
\end{proof}

We next give two background propositions which do not require that $r=6$.

\begin{prop}\label{prop:potential-support} For a fixed limit linear series,
$w$, and column $i$, if $a^i_{(j,j')}>c_i$,
then $(j,j')$
has a component of potential support strictly right of $i-1$, and if
if $a^i_{(j,j')}<c_i$, then $(j,j')$
has a component of potential support strictly left of $i$. 

Conversely, suppose further that $w$ is unimaginative. If $(j,j')$ has a 
component
of potential support strictly right of $i-1$, and if neither $j$ nor $j'$
is exceptional in any column strictly right of $i-1$, then 
$a^i_{(j,j')}>c_i$. Similarly, if $(j,j')$ has a component of potential
support strictly left of $i$, and if neither $j$ nor $j'$ is exceptional
in any column strictly left of $i$, then
$a^i_{(j,j')}<c_i$.

In particular, in the unimaginative case, if the potential support of
$(j,j')$ is disconnected, then at least one of $j,j'$ must be exceptional
somewhere.
\end{prop}

\begin{proof} The first part is straightforward, and we omit the proof.
For the second part, the point is that the unimaginative hypothesis together
with the non-exceptional hypothesis together imply that the relevant portion
of the sequence $a^{i'}_{(j,j')}-c_{i'}$ is nondecreasing in the relevant
range as $i'$ decreases, so in the first case if its positivity for some 
$i' \geq i$ implies it remains positive at $i'=i$, while in the second case
its negativity for some $i' \leq i$ implies it remains negative at $i'=i$.
\end{proof}

\begin{prop}\label{prop:rho-bound} The number of swaps in a given limit
linear series is bounded by $\rho$.

Suppose that we have $\rho$ swaps. Then the swaps
must all be minimal, and occur in genus-$1$ columns, and we cannot have
any exceptional behavior other than what is needed for the swaps.
Moreover, for any unimaginative $w$, the potential
support of the $(j,j')$ row is connected unless the sum of the number of
swaps for which the $j$th row is exceptional and the number of swaps for
which the $j'$th row is exceptional is at least $2$.
\end{prop}

\begin{proof} The first assertions follow immediately from 
Remark \ref{rem:rho-meaning}. For the last, we can have disconnected
potential support in the $(j,j')$ row only if 
the sequence $a^i_{(j,j')}-c_i$ goes from positive to negative as $i$ 
decreases, possibly over multiple columns. But we observe that if only one of
$j$ and $j'$ are exceptional at a swap, which is moreover minimal and in a 
genus-$1$ column, then the sequence $a^i_{(j,j')}-c_i$ can decrease only by 
$1$ as $i$ decreases. Thus, if this occurs only once, it cannot go from 
positive to negative, and we cannot have disconnected potential support.
\end{proof}

Using Lemma \ref{lem:3-at-swap}, we can prove the following.

\begin{cor}\label{cor:alt-multideg} Suppose that $\rho=2$ and $r=6$. Then:
\begin{enumerate}
\item if we are in the ``first $3$-cycle'' situation of Proposition 
\ref{prop:3-cycle-1},
there exists an unimaginative multidegree $\omega'$ such that the
$(j_0-1,j_0)$ row has a unique potentially appearing section in multidegree 
$\omega'$, whose support does not contain $i_0$ or $i_1$, and such that
all the potentially appearing sections are linearly independent.
\item if we are in the ``second $3$-cycle'' situation of Proposition 
\ref{prop:3-cycle-2},
suppose that in the default multidegree $\omega_{\df}$, we have the inequalities
$$2a^{i_0}_{j_0-1}\leq c_{i_0}-1, \text{ and }
2a^{i_1+1}_{j_0-1}\geq c_{i_1+1}+1,$$
with exactly one of the two inequalities satisfied with equality. 
Then there exists an unimaginative multidegree $\omega'$ such that the
$(j_0-1,j_0-1)$ row does not have potentially appearing sections both left 
of $i_0$ and right of $i_1$ in multidegree $\omega'$, and such that
all the potentially appearing sections are linearly independent.
\end{enumerate}
\end{cor}

\begin{proof} (1) In this situation, the $j_0$ and $j_0-1$ rows each
have only one column adding to $d-2$, 
so the potential support of $(j_0-1,j_0)$ can be disconnected only if
$a^{i_0}_{j_0-1}+a^{i_0}_{j_0}=c_{i_0}-1$ and
$a^{i_1+1}_{j_0-1}+a^{i_1+1}_{j_0}=c_{i_1+1}+1$, we have degree $2$
in every column from $i_0$ to $i_1$, inclusive, and we do not have
$\delta_i$ equal to $j_0-1$ or $j_0$ for any $i$ between $i_0$ and $i_1$.
It then follows in particular that $a^{i_0}_{(j_0+1,j_0+1)}\geq c_{i_0}+2$ and
$a^{i_1+1}_{(j_0+1,j_0+1)}\leq c_{i_1+1}-2$, or equivalently, 
$a^{i_1}_{(j_0+1,j_0+1)}\leq c_{i_1}$.

Consider the default multidegree $\omega_{\df}$. If 
the $(j_0-1,j_0)$ row has connected potential support, we are done: since
$\omega_{\df}$ has degree $2$ in both the $i_0$ and $i_1$ columns, the 
$(j_0-1,j_0)$ row cannot have any support in either of these. On the other 
hand, if the $(j_0-1,j_0)$ row has disconnected potential support, then 
we will use Lemma \ref{lem:3-at-swap} 
to verify that we can 
move one degree $3$ into either the $i_0$ or $i_1$ column
while maintaining the independence conclusion of Theorem 
\ref{thm:basic-r6}. We will then obtain the desired statement:
certainly, the $(j_0-1,j_0)$ row will have connected potential 
support. If the $3$ was moved to the $i_0$ column, then the $(j_0-1,j_0)$
still cannot have any potential support at $i_1$. If the $3$ was moved
from the right, we still have 
$a^{i_0}_{j_0-1}+a^{i_0}_{j_0}=c_{i_0}-1$, ruling out potential support
at $i_0$, but if it was moved from the left, then this will decrease
$c_{i_0}$ by $1$, and we will then have 
$a^{i}_{j_0-1}+a^{i}_{j_0}=c_{i}$ for
$i_0 \leq i \leq i_1$, meaning that any potential support at $i_0$ would
have to continue right to $i_1$, but we will still have
$a^{i_1+1}_{j_0-1}+a^{i_1+1}_{j_0}=c_{i_1}+1$, so there cannot be any
potential support at $i_1$. A similar analysis holds if we moved the 
$3$ to $i_1$, proving the desired result.

To prove that we can always move a $3$ as desired, we first make some 
general observations
regarding when we will be able to move degree $3$ from the left or right
onto $i_0$ or $i_1$.
Recall that $\delta_{i_0}=\delta_{i_1}=j_0+1$.
Since
$a^{i_1}_{(j_0+1,j_0+1)}\leq c_{i_1}$, moving a degree $3$ to
$i_1$ from the right will always lead to at most $3$ rows starting in
the $i_1$ column, 
unless $2\bar{a}^{i_1}_3 = c_{i_1}+1$, or equivalently,
\begin{equation}\label{eq:bad-gamma-2}
5- \gamma_{i_1-1}=2\bar{\lambda}_{i_1-1,3}.
\end{equation}
In addition,
$a^{i_1}_{(j_0-1,j_0+1)}< c_{i_1}$, so the $(j_0-1,j_0+1)$ row will not
be among the appearing rows.

We next consider what happens if we move a degree $3$ to $i_0$ from the 
left. This will decrease $c_{i_0}$ by $1$, so we have to rule out that
in multidegree $\omega_{\df}$ we have $2\bar{a}^{i_0}_3 = c_{i_0}$, or 
equivalently,
\begin{equation}\label{eq:bad-gamma-3}
6- \gamma_{i_0-1}=2\bar{\lambda}_{i_0-1,3}.
\end{equation}
Additionally, if $a^{i_0}_{(j_0+1,j_0+1)}\geq c_{i_0}+3$ in $\omega_{\df}$, 
then 
after moving the degree $3$ to $i_0$, none of the other exceptional cases 
of Lemma \ref{lem:3-at-swap} can occur, so as long as we do not have
\eqref{eq:bad-gamma-3}, 
we will have at most three rows with potential support starting at 
$i_0$. The only other possibility is that
$a^{i_0}_{(j_0+1,j_0+1)}= c_{i_0}+2$, which is equivalent to
$2j_0-\gamma_{i_0-1}=2\lambda_{i_0-1,j_0+1}$; moreover,
after moving a $3$ from the left to $i_0$ we will
have $a^{i_0}_{(j_0+1,j_0+1)}= c_{i_0}+3$, so we could potentially be only
in the third exceptional case in Lemma \ref{lem:3-at-swap}. Thus, the only
case for concern is that $j_0+1=4$, so we simply need to check that in
cases where we wish to move a $3$ from the left, we never have
\begin{equation}\label{eq:bad-gamma}
6-\gamma_{i_0-1}=2\lambda_{i_0-1,4}.
\end{equation}
Finally, in either case after the move we will have 
$a^{i_0}_{(j_0,j_0+1)}=a^{i_0}_{(j_0+1,j_0+1)}-1 \geq c_{i_0}+2$, so the
$(j_0,j_0+1)$ row cannot be among the rows starting at $i_0$.

We now describe how to modify our default multidegree, depending on the
location of $i_0$ and $i_1$. 
If we have $\gamma_{i_0}=\gamma_{i_1}=1$, 
then we will move the next $3$ from the right to column $i_1$, 
and we will obtain at most three rows with potential support starting in
$i_1$: by the above observation, it suffices to rule out
\eqref{eq:bad-gamma-2}, but we have $5-\gamma_{i_1-1}=4$.
To have equality we would need
$\bar{\lambda}_{i_1-1,3}=2$, which would imply
$\bar{\lambda}^1_{i_1-1}+\bar{\lambda}^2_{i_1-1} \geq 8$, in which case we 
would not have
had $\gamma_{i_1}=1$ in $\omega_{\df}$.

Next, suppose $\gamma_{i_0}=\gamma_{i_1}=2$, and
we have $\bar{\lambda}_{i_1}^1+\bar{\lambda}_{i_1}^2< 6$.
In this case, we will move the $3$ to $i_0$ from the left, and 
$\gamma_{i_0-1}=2$ in $\omega_{\df}$, so if either 
\eqref{eq:bad-gamma} or \eqref{eq:bad-gamma-3} is satisfied, we must have
$\lambda_{i_0-1,3}\geq 2$. But this would force
$$\bar{\lambda}_{i_1}^1+\bar{\lambda}_{i_1}^2\geq 
\bar{\lambda}_{i_0-1}^1+\bar{\lambda}_{i_0-1}^2 \geq 8,$$ 
contradicting the hypothesis for the case in question.
We again conclude that we have at most $3$ rows starting, and again the
$(j_0,j_0+1)$ row is not among them.

On the other hand, if $\gamma_{i_0}=\gamma_{i_1}=2$, and
we have $\bar{\lambda}_{i_1}^1+\bar{\lambda}_{i_1}^2 \geq 6$, then we
will move a $3$ to $i_1$ from the right, and 
\eqref{eq:bad-gamma-2} is not satisfied for parity reasons, so we will have 
at most three new rows starting.
Finally, if $\gamma_{i_0}=\gamma_{i_1}=3$, neither \eqref{eq:bad-gamma} nor 
\eqref{eq:bad-gamma-3} can be 
satisfied for parity reasons, so we can move a $3$ from the left to $i_0$, 
and have at most three starting rows.

The remaining cases are treated symmetrically, with rows starting replaced
by rows ending.
In each case, we see that the basic structure of the proof of Theorem
\ref{thm:basic-r6} is preserved by our change of multidegree,
so our linear independence is likewise preserved, yielding the desired 
statement.

(2) Suppose that in multidegree $\omega_{\df}$, we have
$$2a^{i_0}_{j_0-1}=c_{i_0}-1, \text{ but }2a^{i_1+1}_{j_0-1}> c_{i_1+1}+1.$$
We will show that we can always move a $3$ from the left to a genus-$1$
column on or right of $i_0$, while preserving linear independence.
This will eliminate potential support in the $(j_0-1,j_0-1)$ row left of
$i_0$, as desired. 
Now, in this situation, we necessarily have
$$2(j_0-1)+1-\gamma_{i_0-1}=2\lambda_{i_0-1,j_0-1},$$
so in particular $\gamma_{i_0-1}$ must be odd.

\textbf{Case $\gamma_{i_0-1}=1$}. We have 
$j_0-1=\lambda_{i_0-1,j_0-1}$.
Because $\bar{\lambda}^1_{i_0-1}+\bar{\lambda}^2_{i_0-1}<5$ necessarily,
this forces $j_0-1=1$, so we have
$\bar{\lambda}_{i_0,1}=\lambda_{i_0-1,1} =1$.

First, if $i_1$ is the genus-$1$ column immediately following $i_0$, 
we observe that if we move the first $3$ to $i_0$, considering only
the inequalities at $i_0$, there can be at most three rows with potential
support starting at $i_0$: $(1,2)$, $(2,2)$ and $(0,j)$ for a unique $j>2$.
But in this case the actual potential support of $(1,2)$ is connected and 
supported strictly to the right of $i_1$.
Thus, there are in fact at most two rows with potential support starting
at $i_0$, and neither of them involves the exceptional row (specifically, 
$j=1$),
so even after moving the first $3$ to $i_0$ we will be able to eliminate
potentially appearing sections from left to right as before.

Next, suppose that $i_1$ is not the genus-$1$ column immediately following 
$i_0$, and denote this column by $i$. Suppose also that there is no
degree-$3$ column between $i_0$ and $i_1$, so that in particular
$\bar{\lambda}_i^1+\bar{\lambda}_i^2 \leq 4$. We observe that we must
also have $\bar{\lambda}_{i,0}= 1$, since we have 
$\bar{\lambda}_{i_0,1}=
\bar{\lambda}_{i_0,2}=1$, and we must have 
$\bar{\lambda}_{i_1-1,2}=
\bar{\lambda}_{i_1-1,3}\geq 1$, so the only way we can avoid having
a column of degree $3$ before $i_1$ is if also 
$\bar{\lambda}_{i_1-1,0}= 1$.
We can then apply Proposition \ref{prop:multideg-flexible} to move the first
$3$ to column $i$, and we will still obtain linear independence. 

Finally, if we have a column of degree $3$ between $i_0$ and $i_1$,
say in column $i$, so that $\bar{\lambda}_i^1+\bar{\lambda}_i^2 = 5$, then 
we claim
that if we move the first $3$ from the left to $i_0$, we will have
at most two potentially appearing sections ending in column $i$, and
at most two potentially appearing sections supported in 
the first column with $\bar{\lambda}_{i'}^1+\bar{\lambda}_{i'}^2 = 6$.
This will prove the desired statement, since we can then eliminate the
potentially appearing sections starting from $i'$ and moving both left
and right from there.
For checking the possible inequalities in column $i'$, moving the
$3$ from the left to $i_0$ won't affect anything, so the argument for
Theorem \ref{thm:basic-r6} implies \textit{a 
priori} that there are at most three rows satisfying the inequalities at
$i'$ for potentially appearing sections to be supported there. We will
check that there is always one such row which satisfies the inequalities
at $i'$, but does not in fact have potential support there. 
Because we have a $3$ between $i_0$ and $i_1$, we must have 
$2a^{i_1+1}_{j_0-1}= c_{i_1+1}+2$.
If $i'<i_1$, the row in question is $(1,1)=(j_0-1,j_0-1)$: indeed, in this 
situation we will have
$a^{i''}_{(j_0-1,j_0-1)}=c_{i''}$ for all $i''$ with $i<i'' \leq i_1$,
so $(j_0-1,j_0-1)$ does satisfy the necessary inequalities at $i'$, but 
its actual potential support (after moving the $3$ to $i_0$) is strictly 
to the right of $i_1$.
On the other hand, if $i'>i_1$, the row in question will be 
$(1,2)=(j_0-1,j_0)$: we have 
$a^{i_1+1}_{j_0-1}+a^{i_1+1}_{j_0} \leq c_{i_1+1}$, so 
$a^{i_1}_{j_0-1}+a^{i_1}_{j_0} < c_{i_1}$, and because the potential
support is connected, it must be strictly left of $i_1$. However,
we claim that we must have $a^{i_1+1}_{j_0-1}+a^{i_1+1}_{j_0} = c_{i_1+1}$, 
and that this must extend through the column $i'$, so that the inequalities
for potential support are satisfied at $i'$. Indeed, the only way this
could fail is if $\delta_{i''}=j_0$ for some $i''$ with $i_0<i'' < i'$.
But we know that 
$\lambda_{i_0-1,j_0-1}=\lambda_{i_0-1,j_0}=1$,
so if $\delta_{i''}=j_0$ anywhere after $i_0$, it increases
$\bar{\lambda}_{i''}^1+\bar{\lambda}_{i''}^2$ to at least $5$. 
Thus, this could only happen
for $i''<i'$ if $i''=i$, which then forces us to have
$\bar{\lambda}_{i_0}^1=3$ and $\bar{\lambda}_{i_0}^2=1$. However, in this
case, because we cannot have a gap between the $j_0-1$ and $j_0+1$ column
at $i_1$, this would force us to also increase 
$\bar{\lambda}_{i''}^1$ to $4$ before $i_1$, which violates our hypothesis
that $i'>i_1$. Thus, in either situation we have shown that the column
$i'$ has at most two potentially appearing sections supported on it, and
it remains to check that the column $i$ has at most two potentially 
appearing sections ending in it. But we either have 
$\bar{\lambda}_{i-1}^1=4$ and 
$\bar{\lambda}_{i-1}^2=0$ or
$\bar{\lambda}_{i-1}^1=3$ and 
$\bar{\lambda}_{i-1}^2=1$,
and one can calculate directly that because we cannot have $\delta_i=0$ or $4$
in the second case, $\delta_i=1$ in either case (recalling that by column
$i$ we have had a swap between rows $1$ and $2$), or $\delta_i=3$ in the
first case, the only rows with potential support ending in column $i$ are 
$(1,2)$ and $(0,j)$ for a unique value of $j$,
yielding the desired statement.

\textbf{Case $\gamma_{i_0-1}=3$}. We can have either 
$j_0-1=2$ and $\lambda_{i_0-1,j_0-1}=1$, 
or $j_0-1=3$ and $\lambda_{i_0-1,j_0-1}=2$. 
First, suppose that
$(j_0-1,j_0)$ has potential support strictly to the right of $i_1$, or 
equivalently, that
there are no columns between $i_0$ and $i_1$ having degree $3$, or with 
$\delta_i=j_0-1$ or $j_0$.
In this case, if we move a $3$ from the left
to $i_0$, by Lemma \ref{lem:3-at-swap} at most four rows satisfy the 
inequalities at $i_0$ to have
potentially appearing sections starting in $i_0$, and we see that these 
include $(j_0-1,j_0)$. 
But $(j_0-1,j_0)$ does not actually have potential support at $i_0$, so
in this case we have at most three rows starting at $i_0$, and none of
them involve the exceptional row (specifically, $j_0-1$), so we
can eliminate this central block just as in Theorem \ref{thm:basic-r6},
and we conclude we still have linear independence. 

Now, the possibility
that we have $\delta_i=j_0-1$ in between $i_0$ and $i_1$ is ruled out
by the inequality $2a^{i_1+1}_{j_0-1}> c_{i_1+1}+1$. 
If there is a column with $\delta_i=j_0$, but no column having
degree $3$ between $i_0$ and $i_1$, we will move the third degree-$3$
from the left to $i_1$, and the $(j_0+1,j_0+1)=(\delta_{i_1},\delta_{i_1})$
row is supported strictly to the right of $i_1$.
In addition
\eqref{eq:3-at-swap} is ruled out by parity reasons, so by Lemma 
\ref{lem:3-at-swap} we have at most three rows starting at $i_1$, and we
also see that $(j_0-1,j_0+1)$ is not among them, as it will have potential
support strictly to the right of $i_1$. Thus, no row involving $j_0-1$
(the exceptional row) has potential support starting at $i_1$,
and in this case we can eliminate all potentially appearing sections just
as in Theorem \ref{thm:basic-r6}.

Next, suppose there is some column with degree $3$ between $i_0$ and 
$i_1$, but no column with $\delta_i=j_0$. In this case, we will
move the $4$th $3$ to the first column $i$ with 
$\bar{\lambda}_{i}^2+\bar{\lambda}_{i}^3 = 9$, and the $3$rd $3$ to $i_0$.
If 
$\bar{\lambda}_{i_1}^2+\bar{\lambda}_{i_1}^3 < 9$, then according to
Proposition \ref{prop:multideg-flexible}, moving the $4$th $3$ doesn't
disrupt linear independence,
and then we are in exactly the same situation as the first case
considered above, with
$(j_0-1,j_0)$ having potential support strictly to the right of $i_1$.
On the other hand, if 
$\bar{\lambda}_{i_1}^2+\bar{\lambda}_{i_1}^3 = 9$,
we will still maintain linear independence, but for different reasons:
we claim that will have at most three rows ending in the $i$th column, no 
row ending in the first column $i'$ with
$\bar{\lambda}_{i'}^2+\bar{\lambda}_{i'}^3 = 8$, and only two rows
ending in the first column $i''$ with
$\bar{\lambda}_{i''}^2+\bar{\lambda}_{i''}^3 = 10$. Thus, we will be able
to eliminate potentially appearing sections from the right, treating the
columns from $i'$ to $i$ as a block to which to apply Lemma
\ref{lem:independence-crit} (3),
and we will in this way eliminate
all potentially appearing sections supported on either side of $i_0$.
This leaves at most one potentially appearing section, which can then
be eliminated. Thus, it suffices to prove the above claim. By the 
argument for Proposition \ref{prop:multideg-flexible}, we have no 
potentially appearing section supported only in the $i$th column,
and at most three continuing from the previous column, so there are
at most three ending in the $i$th column, as claimed. The fact that there
are no rows ending in the $i'$th column is immediate from Corollary
\ref{cor:2-case} and Lemma \ref{lem:impossible}. Finally, we know from
the proof of Theorem \ref{thm:basic-r6} that there at most three rows
satisfying the inequalities in column $i''$ to have potential support 
ending there. Moreover, we see that $(j_0-1,j_0)$ is necessarily one of
them. Indeed, since we have one column with degree $3$ and none with
$\delta_i=j_0-1$ or $j_0$ between $i_0$ and $i_1$, we see that we
necessarily have $a^{i_1+1}_{(j_0-1,j_0)}=c_{i_1+1}$ even after changing
the multidegree. But after $i_1$, any column with $\delta_i=j_0-1$ or $j_0$
will increase $\bar{\lambda}_{i}^2+\bar{\lambda}_{i}^3$,
so this cannot occur strictly between $i_1$ and $i''$, and we conclude that 
$a^{i''}_{(j_0-1,j_0)}=c_{i''}$ as well. Since column $i''$ has degree $3$,
this means that $(j_0-1,j_0)$ satisfies the inequalities to have potential
support ending at $i''$.
But again using that the $4$th $3$ is still left of $i_1$, the 
actual potential support of $(j_0-1,j_0)$ is
contained to the left of $i_1$, so we conclude that column $i''$ has at
most two rows with potential support ending there, completing the proof
of the claim.

It remains to analyze the possibility that we have a column of degree $3$
and a column with $\delta_i=j_0$ in between $i_0$ and $i_1$.
Recall that we have either 
$j_0-1=2$ and $\lambda_{i_0-1,j_0-1}=1$, 
or $j_0-1=3$ and $\lambda_{i_0-1,j_0-1}=2$. 
We first claim that in the latter
case, we cannot have $\delta_i=j_0$ in between $i_0$ and $i_1$
without forcing there to be two columns of degree $3$ in between, or
equivalently, forcing 
$\bar{\lambda}_{i_1}^2+\bar{\lambda}_{i_1}^3 \geq 10$. Indeed, since we
cannot have a gap between $a^{i_0}_{j_0-1}$ and $a^{i_0}_{j_0}$ for the
swap, we must have 
$\bar{\lambda}_{i_0}^2 \geq 5$, and then for the same reason at $i_1$
we must have 
$\bar{\lambda}_{i_1}^2 \geq 6$. But having some $\delta_i=j_0$ also
requires
$\bar{\lambda}_{i}^3 \geq 4$, so we conclude that we would necessarily
have 
$\bar{\lambda}_{i_1}^2+\bar{\lambda}_{i_1}^3 \geq 10$, as claimed.
Thus, it suffices to treat the situation that 
$\lambda_{i_0-1,j_0-1}=1$. In this situation, we 
have 
$\bar{\lambda}_{i_1}^2+\bar{\lambda}_{i_1}^3 \leq 6$,
and we will
move the third $3$ to column $i_0$ and the fourth $3$ to column $i_1$.
We claim that we will have at most two rows with potentially appearing
sections ending in $i_1$, and neither involves the exceptional row 
(specifically, $j_0-1$, which is $2$). Thus, we will be able to eliminate
all potentially appearing sections from the left and from the right of $i_0$,
and finally eliminate the at most one potentially appearing section supported
only at $i_0$. To verify the claim, we see that we necessarily have
$$5 \leq \bar{\lambda}_{i_1}^1 \leq 7, \bar{\lambda}_{i_1}^2 =3, \text{ and }
1 \leq \bar{\lambda}_{i_1}^3 \leq 3.$$
We compute that the only rows satisfying the inequalities to potentially end
at $i_1$ are $(3,4)$, $(0,6)$, $(1,5)$, $(1,6)$ and $(3,5)$, 
but by the
uniqueness part of Corollary \ref{cor:3-case}, we see that the only way we
can have three of these occurring at once is if we have $(3,4)$, $(0,6)$ and
$(1,5)$. However, we also have that $(3,4)$ can only end
if $\bar{\lambda}_{i_1}^3 \leq 2$, $(0,6)$ can only end if 
$\bar{\lambda}_{i_1}^1 \leq 6$, and $(1,5)$ can only end if one of the
preceding two inequalities is strict. But together these imply that
$\bar{\lambda}_{i_1}^1+\bar{\lambda}_{i_1}^3 \leq 7$, meaning that we 
cannot have all the rows ending at $i_1$ under our hypothesis that
the $4$th $3$ comes before $i_1$. 

This concludes the case $\gamma_{i_0-1}=3$.

\textbf{Case $\gamma_{i_0-1}=5$}. We necessarily have
$j_0-1=\lambda_{i_0-1,j_0-1}+2$, and since
$\bar{\lambda}_{i_0-1}^2+\bar{\lambda}_{i_0-1}^3 \geq 10$, we find that
$j_0-1=4$ is the only possibility.
But then if we move the fifth $3$ to $i_0$, even if we obtain two rows
involving $\delta_{i_0}=5$ with potential support ending at $i_0$,
we can have at most one more (necessarily of the form $(j,6)$ for some $j$).
Moreover, the $(4,5)$ row is not one of these, as it will have potential
support starting, not ending, at $i_0$. We can therefore still eliminate the 
block spanning from the first column with 
$\bar{\lambda}_{i}^2+\bar{\lambda}_{i}^3 = 9$ to column $i_0$ just
as before.

The case that
$2a^{i_1+1}_{j_0-1}= c_{i_1+1}+1$ but $2a^{i_0}_{j_0-1}<c_{i_0}-1$ is
handled completely symmetrically, completing the proof. 
\end{proof}

\section{Proofs in the degenerate case}\label{sec:proofs}

To conclude the proof of our main theorem, we show that there are always
multidegrees such that on the one hand, the potentially appearing sections 
are still linearly independent, and on the other hand, tensors coming
from any exact linked linear series must generate at least 
$\binom{r+2}{2}=28$ linearly independent combinations of the potentially 
appearing sections. The key point is that even though there are cases where 
some row may not have any potentially appearing section occuring in our 
linked linear series in the chosen multidegree, in those cases we have to 
have more than one combination of sections from some other row. In fact, 
the arguments of this section are independent of $r$.

So far, in \S \ref{sec:degen-link} we proved statements on existence of mixed
sections, while in the following sections, we proved statements on linear
independence of potentially appearing sections. These threads are related
by the following.

\begin{lem}\label{lem:mixed-in-span} Suppose $s,s'$ are mixed sections
of multidegrees $\md_d(w)$ and $\md_d(w')$, and let $\md_{2d}(w'')$ be 
another multidegree.
 Then $f_{w+w',w''} (s \otimes s')$ lies in the potential
ambient space in multidegree $\md(w'')$.
\end{lem}

\begin{proof} By definition of mixed sections as sums, it suffices to
treat the case that $s$ is obtained purely from gluing together $s^i_j$
for fixed $j$, and $s'$ is obtained from gluing together $s^i_{j'}$ for
fixed $j'$. But in this case the result is clear, since 
$f_{w+w',w''} (s \otimes s')$ must be a combination of potentially appearing
sections from the $(j,j')$ row.
\end{proof}

The following lemma is convenient for cutting down the number of 
possibilities to consider.

\begin{lem}\label{lem:unmixed} Let $s,s'$ be mixed sections
of multidegrees $\md(w)$ and $\md(w')$ and types $(\vec{S},\vec{j})$ and 
$(\vec{S}',\vec{j}')$ respectively. Suppose that for some $i$
with $1<i<N$, we have $\ell_1 \neq \ell_2$ and $\ell_1' \neq \ell_2'$ such
that $i \in S_{\ell_1} \cap S_{\ell_2} \cap S_{\ell_1'}' \cap S_{\ell_2'}'$.
Then for any unimaginative $w''$, 
the map $f_{w+w',w''}$ vanishes identically on $Z_i$.

If further either $\{j_{\ell_1},j_{\ell_2}\}=\{j_{\ell_1'}',j_{\ell_2'}'\}$
or 
$\{j_{\ell_1},j_{\ell_2}\}\cap\{j_{\ell_1'}',j_{\ell_2'}'\}=\emptyset$,
then the same conclusion holds when $i=1$ or $i=N$. 
\end{lem}

\begin{proof} First consider the case $1<i<N$, and write $w=(c_2,\dots,c_N)$
and $w'=(c_2',\dots,c_N')$. The 
hypotheses mean that $w$ allows for support of both
$s^i_j$ and $s^i_{j'}$ for some distinct $j,j'$, so we need to have
$a^i_j,a^i_{j'} \geq c_i$ and $b^i_j,b^i_{j'} \geq d-c_{i+1}$.
Without loss of generality, suppose $a^i_j<a^i_{j'}$. We must have either
$a^i_j+b^i_j<d$ or $a^i_{j'}+b^i_{j'}<d$. Then if 
$b^i_j>b^i_{j'}$ we have either
$c_{i+1} \geq d-b^i_{j'}>d-b^i_j > a^i_j \geq c_i$ or
$c_{i+1} \geq d-b^i_{j'}>a^i_{j'}>a^i_j \geq c_i$, so in either
case we have $c_{i+1} \geq c_i+2$. On the other hand, if $b^i_j<b^i_{j'}$
we have  
$c_{i+1} \geq d-b^i_j>d-b^i_{j'} \geq a^i_{j'} > a^i_j \geq c_i$, so
we again have $c_{i+1} \geq c_i+2$. The same argument holds
for $w'$, so we conclude that $c_{i+1}+c_{i+1}' \geq c_i+c_i'+4$, which
implies that $f_{w+w',w''}$ vanishes on $Z_i$, since
if we write $w''=(c_2'',\dots,c_N'')$, the unimaginative hypothesis means
that $c_{i+1}'' \leq c_i''+3$.

Next, if $i=1$, the unimaginative hypothesis means that $c_2$ is equal
to $2$ or $3$, and it follows (see the proof of Theorem \ref{thm:basic-r6})
that only the rows $(0,0)$, $(0,1)$, $(1,1)$ and $(0,2)$ can have potential 
support in the column, with not both $(0,0)$ and $(1,1)$ occurring.
If $f_{w+w',w''}$ is nonzero on $Z_1$, then
$f_{w+w',w''}(s \otimes s')$ must have $(j_{\ell_u},j_{\ell_v'}')$ parts
with potential support at $i=1$ for $u=1,2$ and $v=1,2$, and this isn't
possible if either $\{j_{\ell_1},j_{\ell_2}\}=\{j_{\ell_1'}',j_{\ell_2'}'\}$
or $\{j_{\ell_1},j_{\ell_2}\}\cap\{j_{\ell_1'}',j_{\ell_2'}'\}=\emptyset$.
The case $i=N$ is symmetric.
\end{proof}

We treat the case of a single swap as follows.

\begin{prop}\label{prop:single-swap-indep} Suppose a limit linear series
contains precisely one swap, occurring between the rows $j_0,j_0-1$ in
column $i_0$. 
Then for any multidegree $\omega$, with notation as in 
Proposition \ref{prop:single-swap},
the images in multidegree $\omega$ of the tensors of pairs of
the $s_j$ for $j \neq j_0$, and 
$s_{j_0}',s_{j_0}''$
contain $\binom{r+2}{2}$ independent linear combinations of the
potentially appearing sections.
\end{prop}

Note that in the proposition statement, we are not asserting that the
actual global sections in multidegree $\omega$ are linearly independent,
merely that the relevant vectors of coefficients (expressing the sections
in question as combinations of the potentially appearing sections) are 
linearly independent.

\begin{proof} For any row $(j,j')$ with
neither $j,j'$ equal to $j_0$, since we have $s_{j}$ and $s_{j'}$
in our linked linear series, we obtain a nonzero contribution from an
$s_{(j,j'),i}$. In particular, considering all $j,j' \neq j_0-1,j_0$ 
we obtain $\binom{r}{2}$ combinations of the potentially appearing
sections, necessarily independent because they are supported in distinct 
rows. Next, consider 
$j \neq j_0,j_0-1$. Then we have the three global sections 
$s_j \otimes s_{j_0}'$, $s_j \otimes s_{j_0}''$ and $s_j \otimes s_{j_0-1}$,
each of which has nonzero image in multidegree $\omega$. We claim that these
three images must contain at least two distinct linear combinations of the
$s_{(j,j_0),i}$ and $s_{(j,j_0-1),i}$.
If $s_j \otimes s_{j_0}'$ has support
in any columns greater than or equal to $i_0$, this necessarily includes 
a nonzero
combination of the $s_{(j,j_0),i}$, which is distinct from the image
of $s_j \otimes s_{j_0-1}$, and we are done. The same holds if 
$s_j \otimes s_{j_0}''$ has support in any columns less than or equal
to $i_0$. The final case is that $s_j \otimes s_{j_0}'$ has support
only in columns strictly less than $i_0$, and $s_j \otimes s_{j_0}''$
has support only in columns strictly greater than $i_0$. In this case,
both may be linear combinations of the $s_{(j,j_0-1),i}$, but since 
their support is disjoint, they must be two distinct combinations, as 
desired.

Thus, we have produced $\binom{r}{2}+2(r-1)=\binom{r+2}{2}-3$ independent
combinations of potentially appearing sections, supported among the
rows $(j,j')$ with $j \neq j_0-1,j_0$.
Finally, we consider the tensors of $s_{j_0-1}, s_{j_0}', s_{j_0}''$,
and claim we obtain three distinct linear combinations, necessarily
supported among the rows $(j_0-1,j_0-1), (j_0-1,j_0), (j_0,j_0)$. 
Consider the images of $s_{j_0}' \otimes s_{j_0}'$, 
$s_{j_0}' \otimes s_{j_0}''$, and $s_{j_0}'' \otimes s_{j_0}''$.
If any of their images contain any portion of the $(j_0,j_0)$ row,
then considering $s_{j_0-1} \otimes s_{j_0-1}, s_{j_0-1} \otimes s_{j_0}',
s_{j_0-1} \otimes s_{j_0}'$, the same argument as above shows we obtain
two distinct combinations of type $(j_0-1,j_0-1)$ and/or $(j_0-1,j_0)$,
so we are done. But the only alternative is that the first
three tensors come from the $(j_0-1,j_0-1)$, $(j_0-1,j_0)$ and
$(j_0-1,j_0-1)$ rows respectively, with the first and last having
disjoint support. Thus, in this case these three are all linearly
independent, and we again obtain the desired conclusion.
\end{proof}

We are now ready to prove the genus $22$ case of Theorem \ref{thm:main};
in fact, we will prove a more general statement for $\rho=1$ cases of
the strong maximal rank conjecture.

The main point is that if we have a smoothing family $\pi:X \to B$ as in
Situation \ref{sit:smoothing},
and a generic linear series $(\sL_{\eta},V_{\eta})$, which
after base change and blowup we may assume is rational on the generic
fiber, we can apply the linked linear series construction both to 
$(\sL_{\eta},V_{\eta})$ and to $(W_{\eta},\sL_{\eta}^{\otimes 2})$,
where $W_{\eta}$ is the image of the multiplication map \eqref{eq:mult-map}.
Then we will have that for any multidegree $\omega$ of total degree
$2d$, and any multidegrees $\omega',\omega''$ of total degree $d$, and any 
sections $s' \in V_{\omega'}$ and $s'' \in V_{\omega''}$, then necessarily 
$f_{\omega'+\omega'',\omega} (s' \otimes s'')$ lies in $W_\omega$; see the 
discussion following
Situation 4.10 of \cite{o-l-t-z}. Thus, in order to give a lower bound on
the rank of \eqref{eq:mult-map}, we can choose many different 
$\omega',\omega''$
and $s',s''$, and show that they span a certain-dimensional subspace of
$(\sL^{\otimes 2})_{\omega}$.

\begin{thm}\label{thm:rho-1} Fix $g,r,d$ with $r \geq 3$ and $\rho = 1$.
In characteristic $0$, suppose that for every curve $X_0$ of genus $g$ as 
in Situation 
\ref{sit:chain-2}, and every refined limit $\fg^r_d$ on $X_0$, there is a 
multidegree $\omega$ such that the potentially appearing sections in 
multidegree $\omega$ are linearly independent.

Then the strong maximal rank conjecture holds for $(g,r,d)$, and more
specifically, if we define $\cD_{(g,r,d)} \subseteq \cM_g$ to be the
set of curves which have a $\fg^r_d$ for which \eqref{eq:mult-map} is not
injective, then the closure in $\overline{\cM}_g$ of $\cD_{g,r,d}$ does 
not contain a general chain of genus-$1$ curves.
\end{thm}

\begin{proof} According to the above discussion together with Theorem
\ref{thm:refined-specialize} and Proposition \ref{prop:exact-limit},
we need to show that an arbitrary exact linked linear series 
on $X_0$ lying over a refined limit linear series admits some multidegree
$\omega$ such that the combined images 
$f_{\omega'+\omega'',\omega} (s' \otimes s'')$ span
an $\binom{r+2}{2}$-dimensional space. For the $\omega$ in the statement, it
then suffices to show that these sections give $\binom{r+2}{2}$ independent
combinations of the potentially appearing sections.
In this case, since $\rho=1$, we can have at most one swap; see Remark
\ref{rem:rho-meaning}. If we
have no swaps, we obtained the desired independence directly from the
independence of the potentially appearing sections, 
using Proposition \ref{prop:simple}. On the other hand, if we have a single 
swap, the desired result follows from Proposition 
\ref{prop:single-swap-indep}.

Because we have proved the statement for all $X_0$ at once, we conclude
the stronger assertion on the closure of $\cD_{(g,r,d)}$ (see the proof
of the last part of Corollary 4.11 of \cite{o-l-t-z} for details of a 
similar argument).
\end{proof}

In particular, the genus-$22$ case of Theorem \ref{thm:main} follows
immediately from Theorem \ref{thm:rho-1} together with Theorem 
\ref{thm:basic-r6}.

\begin{ex}\label{ex:g22-degen} We continue with the running genus-$22$ 
example of Examples \ref{ex:g22-grd} and \ref{ex:g22-simple}. Observe that
in the default multidegree, the (unique) potentially appearing section in
row $(2,3)$ extends from the $7$th column to the $11$th column. This
means that if $s_3'$ and $s_3''$ have the smallest possible portions coming
from the $j=3$ row, so that $s_3'$ only has nonzero $s^i_3$ parts for
$i \geq 8$ and $s_3''$ for $i \leq 10$, then the potentially appearing
section for the $(2,3)$ row cannot come from either $s_2 \otimes s_3'$ 
or $s_2 \otimes s_3''$. This means that these sections (or more precisely,
their images in multidegree $\omega_{\df}$) are forced to yield potentially
appearing sections from the $(2,2)$ row, with $s_2 \otimes s_3'$ necessarily
yielding the one supported from columns $5$ through $7$, and 
$s_2 \otimes s_3''$ necessarily yielding the one supported in column $12$.
Thus, we explicitly see the lack of a $(2,3)$ section being offset by
the inclusion of two independent $(2,2)$ sections.
\end{ex}

We now move on to the $\rho=2$ case, as needed for the genus-$23$ case
of Theorem \ref{thm:main}. Propositions \ref{prop:simple} and
\ref{prop:single-swap-indep} will still suffice to handle the cases that
we have fewer than two swaps, so what remains is to analyze the four cases
with two swaps, which we treat one by one. In all four cases, we will
have swaps occurring in distinct columns $i_0<i_1$, and we will find
convenient to introduce shorthand notation as follows: we will write
for instance
$$s_{j_0}' \otimes s_{j_0+1}''
 = (j_0-1,j_0+1)_{L}+ (j_0-1,j_0)_{R}+ (j_0,j_0+1)$$
to indicate that the image of $s_{j_0}' \otimes s_{j_0+1}''$ in the
relevant multidegree is a combination of potentially appearing sections
from the $(j_0-1,j_0+1)$, $(j_0-1,j_0)$ and $(j_0,j_0+1)$ rows, where
the first must be supported strictly left of $i_0$, and the second strictly
right of $i_1$, and the third has no restrictions on its support. We will
also use subscripts $C$ to denote support strictly between $i_0$ and $i_1$,
$LC$ to denote support strictly left of $i_1$, and $CR$ to denote support
strictly right of $i_0$.

The first case to address is the following.

\begin{prop}\label{prop:repeat-swap-indep} Suppose that we are in the
``repeated swap'' case described in Proposition \ref{prop:repeat-swap}, 
so that our limit linear series
contains precisely two swaps, and these both occur in the same pair of rows, 
say $j_0,j_0-1$.  Then for any unimaginative multidegree $\omega$,
with notation as in Proposition \ref{prop:repeat-swap},
the images in multidegree $\omega$ of the tensors of pairs of
the $s_j$ for $j \neq j_0,j_0-1$, and 
$s_{j_0-1}', s_{j_0-1}'',s_{j_0}',s_{j_0}''$
contain $\binom{r+2}{2}$ independent linear combinations of the
potentially appearing sections.
\end{prop}

\begin{proof} Just as in the proof of Proposition 
\ref{prop:single-swap-indep}, for $j,j' \neq j_0,j_0-1$, our linked linear 
series contains $s_j$ and
$s_{j'}$, so the image of $s_j \otimes s_{j'}$ always gives a potentially 
appearing section from row $(j,j')$.

Now consider $j \neq j_0,j_0-1$; we claim that
$s_j \otimes s_{j_0-1}', s_j \otimes s_{j_0-1}'', 
s_j \otimes s_{j_0}', s_j \otimes s_{j_0}''$
cannot all coincide, and hence have a two-dimensional span. 
Indeed, if $s_j \otimes s_{j_0-1}''$ coincides with $s_j \otimes s_{j_0}'$,
they must be of the form $(j,j_0-1)_L+(j,j_0)_R$.
But the former cannot occur in $s_j \otimes s_{j_0}''$, and the latter
cannot occur in $s_j \otimes s_{j_0-1}'$, so we obtain the desired
independence for these sections.

It remains to show that we have at least three independent sections among
all tensors of the $s_{j_0-1}',s_{j_0-1}'',s_{j_0}',s_{j_0}''$. 
We first consider the four tensor squares; according to Lemma 
\ref{lem:unmixed},
these can only contain types $(j_0-1,j_0-1)$ and $(j_0,j_0)$,
with no type $(j_0-1,j_0)$ appearing. Now, the possible $(j_0,j_0)$
parts of $s_{j_0-1}'^{\otimes 2}$ and $s_{j_0-1}''^{\otimes 2}$ are
disjoint, so we conclude that either these two are distinct, or they
are of pure type $(j_0-1,j_0-1)$. Similarly, 
the sections $s_{j_0}'^{\otimes 2}$ and $s_{j_0}''^{\otimes 2}$ are either
distinct or of pure type $(j_0,j_0)$. Thus, it suffices to show that we
cannot have all of our tensors in the span of a single pair of sections,
each of pure type $(j_0-1,j_0-1)$ or $(j_0,j_0)$. Now, 
$s_{j_0}' \otimes s_{j_0}''$ cannot have a $(j_0-1,j_0-1)$ part, and
$s_{j_0-1}' \otimes s_{j_0-1}''$ cannot have a $(j_0,j_0)$ part, so the
only possibility to consider is that one of our sections is purely of type
$(j_0-1,j_0-1)$, and the other is purely of type $(j_0,j_0)$. 

If the $(j_0-1,j_0-1)$ part occurs in $s_{j_0-1}'' \otimes s_{j_0}'$, 
it must be supported strictly to the left of $i_0$. Then 
$s_{j_0-1}'' \otimes s_{j_0}''$ cannot have a $(j_0-1,j_0-1)$ part,
so must be of type $(j_0,j_0)$, and the support must be strictly to the
right of $i_1$.  On the other hand, if
the $(j_0,j_0)$ part occurs in $s_{j_0-1}'' \otimes s_{j_0}'$,
it must again be supported strictly to the right of $i_1$, and then
$s_{j_0-1}' \otimes s_{j_0}'$ cannot have a $(j_0,j_0)$ part, so must be of
type $(j_0-1,j_0-1)$, again supported to the left of $i_0$. 
But in either case, $s_{j_0-1}' \otimes s_{j_0}''$ cannot be a 
linear combination of these two sections, as desired. 
\end{proof}

We now start imposing that $\rho=2$ and that $X_0$ is left-weighted.
Roughly speaking, the first gives us control over potential support
of sections, while the second ensures that the actual support occurs
where we want it to. 

\begin{prop}\label{prop:rho-2} Suppose that $\rho=2$, and $\omega$ is an
unimaginative multidegree. Then:

\begin{enumerate}
\item In the ``disjoint swap'' situation treated in Proposition 
\ref{prop:disjoint-swap}, 
suppose without loss of generality that $i_0<i_1$. Then
we necessarily have that $i_0$ and $i_1$ both have genus $1$, the two swaps
are minimal, and no rows are exceptional except row $j_0-1$ at $i_0$ and
row $j_1-1$ at $i_1$. In multidegree $\omega$, the potential support of every
$(j,j')$ is connected except possibly for $(j_0-1,j_0-1)$, 
$(j_1-1,j_1-1)$, and $(j_0-1,j_1-1)$. Moreover,
if $(j_0-1,j_1-1)$ has disconnected potential support in
multidegree $\omega$, the potential support must be made up of two components,
one contained strictly to the right of $i_1$, and one contained strictly
to the left of $i_0$, and 
the potential support of $(j_0-1,j_1)$ is contained strictly right of $i_1-1$, 
and
the potential support of $(j_0,j_1-1)$ is 
contained strictly left of $i_0+1$. Finally, if 
the potential support of $(j_0-1,j_1)$ is contained strictly left of $i_0$,
then $(j_0-1,j_1-1)$ must also have a component of potential support
contained strictly left of $i_0$, and if 
the potential support of $(j_0,j_1-1)$ is contained strictly right of $i_1$,
then $(j_0-1,j_1-1)$ must also have a component of potential support
contained strictly right of $i_1$.

\item In the ``first $3$-cycle'' situation described in Proposition 
\ref{prop:3-cycle-1}, 
we necessarily have that $i_0$ and $i_1$ both have genus $1$, the two swaps
are minimal, and no rows are exceptional except row $j_0$ at $i_0$ and
row $j_0-1$ at $i_1$. In multidegree $\omega$, the potential support of every
In multidegree $\omega$, the potential support of every
$(j,j')$ is connected except possibly for $(j_0-1,j_0-1)$, 
$(j_0-1,j_0)$, and $(j_0,j_0)$. Moreover, if for
some $j$, the potential support of $(j,j_0)$ has a component strictly to 
the left of $i_0$, then the potential support of $(j,j_0-1)$ is entirely
contained strictly to the left of $i_0$, and if 
the potential support of $(j,j_0-1)$ has a component strictly to 
the right of $i_1$, then the potential support of $(j,j_0)$ is entirely
contained strictly to the right of $i_1$.

Finally, if $(j_0-1,j_0)$ 
has potential support contained entirely strictly to the left of $i_1$, 
then the potential support of 
$(j_0-1,j_0+1)$ cannot be contained to the right of $i_1$; if it has 
potential support contained
entirely strictly to the right of $i_0$, then 
the potential support of $(j_0,j_0+1)$ cannot be contained to 
the left of $i_0$; and
if it has potential support contained entirely strictly between $i_0$ and 
$i_1$, then $(j_0-1,j_0-1)$ has potential support contained entirely
strictly to the left of $i_1$, and $(j_0,j_0)$ has potential support 
contained entirely strictly to the right of $i_0$.
\end{enumerate}
\end{prop}

\begin{proof} We write as usual $\omega=\md(w)$ with
$w=(c_2,\dots,c_g)$.

(1) The first assertions follow from Proposition 
\ref{prop:rho-bound}, and following the proof we see further 
in order for $(j_0-1,j_1-1)$ to
have disconnected support, the support must be split between strictly 
right of $i_1$ and strictly left of $i_0$, as claimed.
Next, if the potential support of $(j_0-1,j_1-1)$ has a component lying
strictly right of $i_1$, then we have
$$a^{i_1}_{(j_0-1,j_1)}=a^{i_1+1}_{(j_0-1,j_1)}-1
=a^{i_1+1}_{(j_0-1,j_1-1)}-2>c_{i_1+1}-2\geq c_{i_1},$$
and this implies (using our previous connectedness statement)
that the potential support of $(j_0-1,j_1)$ is supported strictly to
the right of $i_1-1$, as desired. The corresponding statement on support 
left of $i_0$ and $i_0+1$ follows similarly.
Finally, if the potential support of $(j_0-1,j_1)$ is contained strictly 
left of $i_0$, then we have
$$a^{i_0}_{(j_0-1,j_1-1)}< a^{i_0}_{(j_0-1,j_1)} <c_{i_0},$$
so $(j_0-1,j_1-1)$ must also have a component of potential support
strictly left of $i_0$, as desired. The last statement on support strictly
right of $i_1$ follows similarly.

(2) 
Most of the argument is similar to (1). 
For the support of $(j,j_0)$ to have a component strictly to the
left of $i_0$ we must have $a^{i_0}_{(j,j_0)} \leq c_{i_0}-1$, and then
$a^{i_0}_{(j,j_0-1)}<c_{i_0}-1$, so arguing as in (1) we conclude that (even 
if $j=j_0-1$ or $j_0$) the support of $(j,j_0-1)$ is connected and
strictly to the left of $i_0$. The statement on support to the right of
$i_1$ is proved in exactly the same way.
For the last assertion, note that the $(j_0-1,j_0-1)$ row has no support at 
$i_1$, and the $(j_0,j_0)$ has no support at $i_0$, since both sum to
$2d-4$ in the relevant columns.
\end{proof}

\begin{prop}\label{prop:left-weighted} Suppose that $X_0$ is left-weighted,
and that the rows $j,j'$ have no exceptional behavior in any genus-$0$
columns. Then the image of $s_j \otimes s_{j'}$ in any unimaginative 
multidegree $\omega$ is equal to the leftmost potentially appearing section 
in the $(j,j')$ row.
\end{prop}

\begin{proof} The lack of exceptional behavior away from genus-$1$ components
means that the $a^i_{(j,j')}$ are constant on the genus-$0$ components. The
idea is then that the left-weighting means that the leftmost negative value
of $a^i_{(j,j')}-c_i$ is repeated so many times that it must lead to a strict 
minimum of the partial sums. Compare the proof of
Proposition \ref{prop:mixed-controlled}, where in \eqref{eq:mixed-controlled}
we now replace $d$ by $2d$ due to having passed to the tensor square.
\end{proof}

\begin{prop}\label{prop:disjoint-swaps-indep} 
Suppose that we are in the
``disjoint swap'' case described in Proposition \ref{prop:disjoint-swap},
so that our limit linear series
contains precisely two swaps, and these occur in disjoint pairs of rows, 
say $j_0,j_0-1$ and $j_1,j_1-1$. 
Suppose further that $\rho=2$, and that $X_0$ is left-weighted.
Then for any unimaginative multidegree $\omega$, with notation as in 
Proposition \ref{prop:disjoint-swap}, if we suppose that we have chosen
$s_{j_0}'$ and $s_{j_1}'$ as allowed by Proposition 
\ref{prop:mixed-controlled}, then
the images in multidegree $\omega$ of the tensors of pairs of
the $s_j$ for $j \neq j_0,j_1$, and $s_{j_0}',s_{j_0}'',s_{j_1}',s_{j_1}''$
contain $\binom{r+2}{2}$ independent linear combinations of the
potentially appearing sections.
\end{prop}

\begin{proof} 
Without loss of generality, assume that $i_0<i_1$. 
First note that by Proposition \ref{prop:rho-2} (1), the hypothesis that 
$\rho=2$ means that in order to have
two swaps, they both must occur at genus-$1$ components. Then by
Proposition \ref{prop:mixed-controlled}, we may assume that 
$s_{j_1}'$ is controlled, and that the $j_1$-part of
$s_{j_1}'$ does not contain any genus-$1$ component left of $i_1$.
We also have that every $(j,j')$ has
connected potential support unless $j,j' \in \{j_0-1,j_1-1\}$.

Now, if we have $j,j' \neq j_0,j_0-1,j_1,j_1-1$, then we know
that $f_{w_j+w_j',w} (s_j \otimes s_{j'})$ is nonzero and composed of
$s^i_{(j,j')}$. Now, suppose $j \neq j_0,j_0-1,j_1,j_1-1$. Then the
same argument as in Proposition \ref{prop:single-swap-indep}
also shows that if we consider the images in multidegree $\omega$
of $s_j \otimes s_{j_0-1}$, $s_j \otimes s_{j_0}'$, and 
$s_j \otimes s_{j_0}''$, we either obtain one section of type $(j,j_0-1)$
and one with a contribution of type $(j,j_0)$, or two sections of type 
$(j,j_0-1)$, but having disjoint support. The same holds with $j_1$ in
place of $j_0$. Together, these produce 
$\binom{r-2}{2}+4(r-3)=\binom{r+2}{2}-10$ linearly independent combinations.
It thus suffices to show that we have $10$ linearly independent combinations
coming from tensor products of pairs of the sections 
$s_{j_0-1},s_{j_0}',s_{j_0}'',s_{j_1-1},s_{j_1}',s_{j_1}''$. 
Just as in the proof of Proposition \ref{prop:single-swap-indep},
tensor products of the first three sections yield three
independent combinations, with contributions contained among the types 
$(j_0-1,j_0-1)$, $(j_0-1,j_0)$, and $(j_0,j_0)$. Tensor products of the
last three sections likewise yield three combinations, with $j_1$
replacing $j_0$ in the types. 

It remains to consider the tensors with types contained among $(j_0-1,j_1-1)$,
$(j_0-1,j_1)$, $(j_0,j_1-1)$ and $(j_0,j_1)$.
First suppose that $(j_0-1,j_1-1)$ has connected potential support
in multidegree $\omega$. Then just as in the single-swap case, at least one of
$s_{j_0-1} \otimes s_{j_1}', s_{j_0-1} \otimes s_{j_1}''$ must involve
a $(j_0-1,j_1)$ part, and at least one of $s_{j_0}' \otimes s_{j_1-1},
s_{j_0}'' \otimes s_{j_1-1}$ must involve a $(j_0,j_1-1)$ part. Since
$s_{j_0-1} \otimes s_{j_1-1}$ is pure of type $(j_0-1,j_1-1)$, and all
of these have unique potential support, we find that the span of these
sections contains the (unique) pure types of each of $(j_0-1,j_1-1)$,
$(j_0,j_1-1)$ and $(j_0-1,j_1)$. Thus, if we have anything with a
nonzero part of type $(j_0,j_1)$, this gives a fourth independent
combination. On the other hand, if nothing has a $(j_0,j_1)$ part, then
we must have the following:
$$s_{j_0}' \otimes s_{j_1}''=(j_0-1,j_1)_L+(j_0,j_1-1)_R,$$
$$s_{j_0}' \otimes s_{j_1}'=(j_0-1,j_1-1)_L+ (j_0,j_1-1)_{LC}
+(j_0-1,j_1)_L, \text{ and}$$
$$s_{j_0}'' \otimes s_{j_1}''=(j_0-1,j_1)_{CR}+(j_0-1,j_1-1)_R+
(j_0,j_1-1)_R.$$

First consider the possibility that the 
$(j_0-1,j_1)_L$ part of $s_{j_0}' \otimes s_{j_1}''$ is nonzero. Then
by Proposition \ref{prop:rho-2} (1), we have that $(j_0-1,j_1-1)$ has
support strictly left of $i_0$ too, which in turn means that 
$(j_0,j_1-1)$ can't have support strictly right of $i_1$. But this leaves 
no possibility for $s_{j_0}'' \otimes s_{j_1}''$. On the other hand, if
the $(j_0,j_1-1)_R$ part of $s_{j_0}' \otimes s_{j_1}''$ is nonzero, we
have that $(j_0-1,j_1-1)$ must have support strictly right of $i_1$, and
hence that $(j_0-1,j_1)$ can't have support strictly left of $i_0$, 
leaving no possibility for $s_{j_0}' \otimes s_{j_1}'$. We conclude that
it is not possible for these tensors not to have some $(j_0,j_1)$ part,
giving the desired four independent combinations when $(j_0-1,j_1-1)$
has connected potential support.

It remains to treat the case that $(j_0-1,j_1-1)$ has disconnected 
potential support in multidegree $\omega$. Then Proposition \ref{prop:rho-2}
tells us that this potential support has two parts, contained strictly 
left of $i_0$ and right of $i_1$ respectively. Moreover, it says that
the potential support of $(j_0-1,j_1)$ is contained strictly right of 
$i_1-1$ and the potential support of $(j_0,j_1-1)$ is 
contained strictly left of $i_0+1$. This forces 
$s_{j_0}' \otimes s_{j_1}''$ to be of pure $(j_0,j_1)$ type. Now,
we observe that two of the sections $s_{j_0-1} \otimes s_{j_1-1}, 
s_{j_0-1} \otimes s_{j_1}', s_{j_0-1} \otimes s_{j_1}''$ must be
independent, either involving a $(j_0-1,j_1)$ part and a $(j_0-1,j_1-1)$
part, or two $(j_0-1,j_1-1)$ parts. Similarly,
$s_{j_0}' \otimes s_{j_1-1}$ and $s_{j_0}'' \otimes s_{j_1-1}$ must either
involve a $(j_0,j_1-1)$ part or two $(j_0-1,j_1-1)$ parts. We see that the 
only way to avoid having four independent combinations would be if these 
five tensors are all of pure type $(j_0-1,j_1-1)$, necessarily achieving 
support independently both on the left and right. But we note that
because the potential support of $(j_0,j_1-1)$ is contained strictly 
left of $i_0+1$, and because (in the disconnected support case) we must
have $a^i_{(j_0-1,j_1-1)}=c_i$ for $i_0<i \leq i_1$, the only way that 
$s_{j_0}''\otimes s_{j_0-1}$ can fail to have a $(j_0,j_1-1)$ part 
is if $s_{j_0}''$ is not controlled, and more specifically if 
its $j_0$ portion does not extend more than halfway to the next genus-$1$ 
component after $i_0$.
On the other hand, $s_{j_1}'$ is controlled and has $j_1$ part not 
containing any genus-$1$ component smaller than $i_1$, so we conclude that
in this situation
its $j_1$ part is disjoint from the $j_0$ part of $s_{j_0}''$, and then
$s_{j_0}'' \otimes s_{j_1}'=(j_0,j_1-1)+(j_0-1,j_1)$, and gives a fourth
independent combination.
This completes the proof of the proposition.
\end{proof}

We now move on to consider the two remaining cases, both involving
a pair of swaps in overlapping columns. 

\begin{prop}\label{prop:3-cycle-1-indep}
Suppose that we are in the ``first $3$-cycle'' situation described in
Proposition \ref{prop:3-cycle-1}, so that our limit linear series 
contains precisely two swaps, with
one swap between the $j_0$th and $(j_0+1)$st rows occurring in the $i_0$th 
column, and a second swap between the $(j_0-1)$st and $(j_0+1)$st rows in 
the $i_1$st column for some $i_1>i_0$. Suppose further that $\rho=2$, and
that we have an unimaginative multidegree $\omega$ such that the
$(j_0-1,j_0)$ row has a unique potentially appearing section in multidegree 
$\omega$, whose support does not contain $i_0$ or $i_1$. Then with notation 
as in Proposition \ref{prop:3-cycle-1},
the images in multidegree $\omega$ of the tensors of pairs of
the $s_j$ for $j \neq j_0+1$, and $s_{j_0+1}',s_{j_0+1}'',s_{j_0+1}'''$
contain $\binom{r+2}{2}$ independent linear combinations of the
potentially appearing sections.
\end{prop}

\begin{proof} 
We first show that for 
$j \neq j_0-1,j_0,j_0+1$, the sections 
$$s_j \otimes s_{j_0-1}, s_j \otimes s_{j_0}, s_j \otimes s_{j_0+1}',
s_j \otimes s_{j_0+1}'', s_j \otimes s_{j_0+1}'''$$
must yield at least three independent combinations.
But the first two tensors yield $(j,j_0-1)$ and 
$(j,j_0)$ parts, so if any of the last three have any $(j,j_0+1)$ part, 
we obtain the desired independence. On the other hand, if not we find that
$$s_j \otimes s_{j_0+1}'=(j,j_0-1)_{LC}+(j,j_0)_L;$$
$$s_j \otimes s_{j_0+1}''=(j,j_0-1)_R+(j,j_0)_{CR};$$
$$s_j \otimes s_{j_0+1}'''=(j,j_0)_L+(j,j_0-1)_R.$$
If the $(j,j_0)_L$ part of the last tensor is nonzero, then 
by Proposition \ref{prop:rho-2}, the potential support of both the 
$(j,j_0-1)$ and $(j,j_0)$ rows are connected and contained strictly to the
left of $i_0$, leaving no possibility for the second tensor.  
But if the $(j,j_0-1)_R$ part of the last tensor is nonzero, then we similarly
have that the potential support of both the 
$(j,j_0-1)$ and $(j,j_0)$ rows are contained strictly to the
right of $i_1$, leaving no possibility for the first tensor.  
Thus, we reach a contradiction, and conclude that we must obtain a
$(j,j_0+1)$ part, giving the desired three independent combinations.

Next, we consider the $15$ tensors arising from
$$s_{j_0-1},s_{j_0},s_{j_0+1}',s_{j_0+1}'',s_{j_0+1}''';$$
we need to show that these yield $6$ independent linear combinations.
By hypothesis, we have that the potential
support of the $(j_0-1,j_0)$ row is connected and does not contain $i_0$
or $i_1$, so we organize cases according to its support.
First suppose that the support of the $(j_0-1,j_0)$ row is entirely to the
left of $i_0$; then according to Proposition \ref{prop:rho-2}, the
same holds for the $(j_0-1,j_0-1)$ row, and the $(j_0-1,j_0+1)$ row
cannot have its support to the right of $i_1$. We then see that 
$s_{j_0-1} \otimes s_{j_0+1}''$ cannot have any $(j_0-1,j_0-1)$ or
$(j_0-1,j_0)$ parts, so must be of $(j_0-1,j_0+1)$ type. Similarly,
$s_{j_0+1}'' \otimes s_{j_0+1}'''$ cannot have any $(j_0-1,j_0-1)$,
$(j_0-1,j_0)$, or $(j_0-1,j_0+1)$ parts, so it must contain $(j_0,j_0+1)$
or $(j_0+1,j_0+1)$ parts. In addition, the pair
$s_{j_0} \otimes s_{j_0+1}''$ and $s_{j_0} \otimes s_{j_0+1}'''$ 
must contain either a $(j_0,j_0+1)$ part, or two distinct $(j_0,j_0)$ parts,
supported left and right of $i_0$, respectively. Given that we always
have $(j_0-1,j_0-1)$, $(j_0-1,j_0)$ and $(j_0,j_0)$ parts, the only
way we could fail to have produced six independent combinations is if
$s_{j_0+1}'' \otimes s_{j_0+1}'''$ has type $(j_0,j_0+1)$, and we have
only one $(j_0,j_0)$ part. But then considering 
$s_{j_0+1}''^{\otimes 2}$ and $s_{j_0+1}'''^{\otimes 2}$ and using
Lemma \ref{lem:unmixed}, we find
that we must produce a $(j_0+1,j_0+1)$ part or two distinct
$(j_0,j_0)$ parts, so we necessarily obtain the sixth combination.

Similarly, if the potential support of the $(j_0-1,j_0)$ row is entirely to 
the right of $i_1$, then Proposition \ref{prop:rho-2} tells us that
the same holds for $(j_0,j_0)$, and that the potential support of the
$(j_0,j_0+1)$ row cannot be to the left of $i_0$. Then
$s_{j_0} \otimes s_{j_0+1}'$ must be of $(j_0,j_0+1)$ type, and
$s_{j_0+1}' \otimes s_{j_0+1}'''$ must have $(j_0-1,j_0+1)$
or $(j_0+1,j_0+1)$ parts. The pair
$s_{j_0-1} \otimes s_{j_0+1}'$ and $s_{j_0-1} \otimes s_{j_0+1}'''$ 
must contain either a $(j_0-1,j_0+1)$ part, or two distinct $(j_0-1,j_0-1)$ 
parts, and in either case the tensors
$s_{j_0+1}'^{\otimes 2}$ and $s_{j_0+1}'''^{\otimes 2}$ (together with
the usual tensors of $s_{j_0-1}$ and $s_{j_0}$) must complete the
six independent combinations.

Finally, if the potential support of the $(j_0-1,j_0)$ row is between
the $i_0$ and $i_1$ columns, then by Proposition \ref{prop:rho-2}, we know
that the potential support of $(j_0-1,j_0-1)$ is left of $i_1$ and the
potential support of $(j_0,j_0)$ is right of $i_0$. We then see that the
tensors $s_{j_0-1} \otimes s_{j_0+1}'''$, $s_{j_0} \otimes s_{j_0+1}'''$, 
and $s_{j_0+1}'''^{\otimes 2}$ must be pure of types $(j_0-1,j_0+1)$,
$(j_0,j_0+1)$, and $(j_0+1,j_0+1)$ respectively, yielding the desired
six combinations.
\end{proof}

\begin{prop}\label{prop:3-cycle-2-indep}
Suppose that we are in the ``second $3$-cycle'' situation described in
Proposition \ref{prop:3-cycle-2}, so that our limit linear series contains 
precisely two swaps, with
one swap between the $(j_0-1)$st and $j_0$th rows occurring in the $i_0$th 
column, and a second swap between the $(j_0-1)$st and $(j_0+1)$st rows in 
the $i_1$st column for some $i_1>i_0$. Suppose further that $\rho=2$,
that $X_0$ is left-weighted, and that we have an unimaginative 
$w=(c_2,\dots,c_N)$ satisfying one of the following three conditions:
\begin{enumerate}
\itm the $(j_0-1,j_0-1)$ row does not have potentially appearing sections 
both left of $i_0$ and right of $i_1$; or
\itm
$2a^{i_0}_{j_0-1}= c_{i_0}-1, \text{ and } 2a^{i_1+1}_{j_0-1}= c_{i_1+1}+1$; or
\itm
$2a^{i_0}_{j_0-1}= c_{i_0}-2, \text{ and } 2a^{i_1+1}_{j_0-1}= c_{i_1+1}+2$,
and $w$ has degree $2$ in both $i_0$ and $i_1$.
\end{enumerate}
Then 
with notation as in Proposition \ref{prop:3-cycle-2},
the images in multidegree $\md(w)$ of the tensors of pairs of
the $s_j$ for $j \neq j_0,j_0+1$, and 
$s_{j_0}',s_{j_0}'',s_{j_0+1}',s_{j_0+1}'',s'''$
contain $\binom{r+2}{2}$ independent linear combinations of the
potentially appearing sections.
\end{prop}

\begin{proof} First suppose $j \neq j_0-1,j_0,j_0+1$; we show that we can
always obtain three linearly independent combinations of potential
appearing sections from the rows $(j,j_0-1)$, $(j,j_0)$ and $(j,j_0+1)$.
$s_j \otimes s_{j_0-1}$ always yields a pure $(j,j_0-1)$ part.
If $S_2'=S_4''=\{1,\dots,N\}$, then
$s_j \otimes s_{j_0}'$ has a nonzero $(j,j_0)$ part and no $(j,j_0+1)$ part,
while 
$s_j \otimes s_{j_0+1}''$ has a nonzero $(j,j_0+1)$ part, so we get the
desired three combinations. Otherwise, 
we have 
\begin{align*}
s_j \otimes s_{j_0}' & = (j,j_0-1)_L+(j,j_0) \\
s_j \otimes s_{j_0}'' & = (j,j_0)+(j,j_0+1)_{R'} + (j,j_0-1)_{CR} \\
s_j \otimes s_{j_0+1}' & = (j,j_0-1)_{LC}+(j,j_0)_{L'} + (j,j_0+1) \\
s_j \otimes s_{j_0+1}'' & = (j,j_0+1)+(j,j_0-1)_R \\
s_j \otimes s''' & = (j,j_0)+(j,j_0-1)_C + (j,j_0+1),
\end{align*}
where $R'$ and $L'$ denote possible support at and right of $i_1$ and at and 
left of $i_0$, respectively, and if $s_j \otimes s_{j_0}''$ has a nonzero
$(j,j_0+1)$ part with support containing $i_1$, its $(j,j_0)$ part must
be nonzero, and similarly for the $(j,j_0)$ and $(j,j_0+1)$ parts of 
$s_j \otimes s_{j_0+1}'$.
Now, suppose that $(j,j_0-1)$ has connected potential support which is
not contained strictly right of
$i_0$. Then $(j,j_0+1)$ cannot have any potential support strictly right of 
$i_1$ without also forcing $(j,j_0-1)$ to have potential support strictly 
right of $i_1$, 
so the $(j,j_0)$ part of $s_j \otimes s_{j_0}''$ must be nonzero.
But then adding
$s_j \otimes s_{j_0+1}'' = (j,j_0+1)$ and $s_j \otimes s_{j_0-1}$ yields 
three independent sections.
Similarly, if $(j,j_0-1)$ has connected potential support not contained 
strictly left of
$i_1$, then $(j,j_0)$ cannot have potential support strictly left of $i_0$, 
so $s_j \otimes s_{j_0+1}'$ has nonzero $(j,j_0+1)$ part, and adding
$s_j \otimes s_{j_0}' = (j,j_0)$ and $s_j \otimes s_{j_0-1}$ yields 
the desired combinations. For connected potential support,
the only remaining possibility is that $(j,j_0-1)$ has potential
support strictly between $i_0$ and $i_1$, in which case 
$s_j \otimes s_{j_0}' = (j,j_0)$ and 
$s_j \otimes s_{j_0+1}'' = (j,j_0+1)$.

Finally, since $\rho=2$, the only remaining possibility is that 
$(j,j_0-1)$ has potential support both left of $i_0$
and right of $i_1$, and in this case we must have  
$a^{i_0}_{(j,j_0-1)}=c_{i_0}-1$ and 
$a^{i_1+1}_{(j,j_0-1)}=c_{i_1+1}+1$. 
Then $(j,j_0+1)$ cannot have potential support strictly right of $i_1$, and 
$(j,j_0)$ cannot have potential support strictly left of $i_0$, so as 
above we find that if the $(j,j_0+1)$ part of $s_j \otimes s_{j_0}''$ is 
nonzero (necessarily with support at $i_1$),
then the $(j,j_0)$ part must also be nonzero, and if the
$(j,j_0)$ part of $s_j \otimes s_{j_0+1}'$ is nonzero, then the $(j,j_0+1)$
part must also be nonzero.
Now, we have $s_j \otimes s_{j_0}'$ and $s_j \otimes s_{j_0+1}''$ linearly
independent always, and the only way they could fail to be independent from
$s_j \otimes s'''$ is if either $s_j \otimes s_{j_0}'=(j,j_0)$ or
$s_j \otimes s_{j_0+1}''=(j,j_0+1)$, while the only way they could fail to
be independent from $s_j \otimes s_{j_0-1}$ if is either 
$s_j \otimes s_{j_0}'=(j,j_0-1)_L$ or
$s_j \otimes s_{j_0+1}''=(j,j_0-1)_R$. If 
$s_j \otimes s_{j_0}'=(j,j_0)$ and $s_j \otimes s_{j_0+1}''=(j,j_0-1)_R$, 
we see that $s_j \otimes s_{j_0+1}'$ necessarily gives a third independent
combination, while if
$s_j \otimes s_{j_0}'=(j,j_0-1)_L$ and $s_j \otimes s_{j_0+1}''=(j,j_0+1)$, 
we see that $s_j \otimes s_{j_0}''$ necessarily gives a third independent
combination.

It remains to show that we can get six independent combinations from
the rows $(j_0-1,j_0-1)$, $(j_0-1,j_0)$, $(j_0-1,j_0+1)$, $(j_0,j_0)$,
$(j_0,j_0+1)$, and $(j_0+1,j_0+1)$.
If $S_2'=S_4''=\{1,\dots,N\}$, then we immediately get that the six 
tensors coming from $s_{j_0-1},s_{j_0}',s_{j_0+1}''$ are linearly
independent, as desired.
Otherwise, we will make use of the mixed section $s'''$ to handle certain
cases. For reference, we write out the form of all the relevent tensors
of $s_{j_0-1},s_{j_0}',s_{j_0}'',s_{j_0+1}',s_{j_0+1}''$. Note that we
are making use of Lemma \ref{lem:unmixed} in the case of self-tensors.
\begin{align*}
s_{j_0-1} \otimes s_{j_0}' 
& = (j_0-1,j_0-1)_L+(j_0-1,j_0) 
\displaybreak[1] \\
s_{j_0-1} \otimes s_{j_0}'' 
& = (j_0-1,j_0)+(j_0-1,j_0+1)_{R'}+(j_0-1,j_0-1)_{CR} 
\displaybreak[1] \\
s_{j_0-1} \otimes s_{j_0+1}' 
& = (j_0-1,j_0-1)_{LC}+(j_0-1,j_0)_{L'}+(j_0-1,j_0+1) 
\displaybreak[1] \\
s_{j_0-1} \otimes s_{j_0+1}'' 
& = (j_0-1,j_0+1)+(j_0-1,j_0-1)_R 
\displaybreak[1] \\
s_{j_0}' \otimes s_{j_0}' 
& = (j_0-1,j_0-1)_L+(j_0,j_0) 
\displaybreak[1] \\
s_{j_0}'' \otimes s_{j_0}'' 
& = (j_0,j_0)+ (j_0+1,j_0+1)_{R}+(j_0-1,j_0-1)_{CR} 
\displaybreak[1] \\
s_{j_0}' \otimes s_{j_0}'' 
& = (j_0-1,j_0)+ (j_0,j_0)+(j_0,j_0+1)_{R'} 
\displaybreak[1] \\
s_{j_0+1}' \otimes s_{j_0+1}'
& = (j_0-1,j_0-1)_{LC}+ (j_0,j_0)_L+(j_0+1,j_0+1) 
\displaybreak[1] \\
s_{j_0+1}'' \otimes s_{j_0+1}''
& = (j_0+1,j_0+1)+(j_0-1,j_0-1)_R 
\displaybreak[1] \\
s_{j_0+1}' \otimes s_{j_0+1}'' 
& = (j_0-1,j_0+1)+ (j_0,j_0+1)_{L'}+(j_0+1,j_0+1) 
\displaybreak[1] \\
s_{j_0}' \otimes s_{j_0+1}'
& = (j_0-1,j_0-1)_{L}+ (j_0-1,j_0)_{LC}+(j_0-1,j_0+1)_L \\
& \quad +(j_0,j_0)_{L'} + (j_0,j_0+1) 
\displaybreak[1] \\
s_{j_0}' \otimes s_{j_0+1}''
& = (j_0-1,j_0+1)_{L}+ (j_0-1,j_0)_{R}+ (j_0,j_0+1) 
\displaybreak[1] \\
s_{j_0}'' \otimes s_{j_0+1}''
& = (j_0-1,j_0-1)_{R}+ (j_0-1,j_0)_{R}+(j_0-1,j_0+1)_{CR} \\
& \quad +(j_0,j_0+1) + (j_0+1,j_0+1)_{R'}. 
\end{align*}

As above, we separate out cases by the potential support of the
$(j_0-1,j_0-1)$ row. Note that because the entries sum to $2d-4$ in
both the $i_0$ and $i_1$ columns, the $(j_0-1,j_0-1)$ row cannot have any
potential support in either of these columns in any unimaginative 
multidegree.
First suppose the potential support is strictly 
to the left of $i_0$. In this case none of the relevant rows can 
have potential support extending right of $i_1$,
so we get  
$s_{j_0-1} \otimes s_{j_0+1}'' = (j_0-1,j_0+1)$,
$s_{j_0}'' \otimes s_{j_0}'' = (j_0,j_0)$, and
$s_{j_0+1}'' \otimes s_{j_0+1}'' = (j_0+1,j_0+1)$, and the $(j_0-1,j_0)$
part of $s_{j_0-1} \otimes s_{j_0}''$ must be nonzero.
We also have
$s_{j_0}' \otimes s_{j_0+1}'' = (j_0-1,j_0+1)_{L}+ (j_0,j_0+1)$ and
$s_{j_0}'' \otimes s_{j_0+1}''= (j_0-1,j_0+1)_{CR}+(j_0,j_0+1)
+(j_0+1,j_0+1)_{R'}$, where again the latter has to have nonzero $(j_0,j_0+1)$
part unless it is equal to $(j_0-1,j_0+1)_{CR}$, so
these must either yield a nonzero $(j_0,j_0+1)$ part, or two independent
$(j_0-1,j_0+1)$ parts (which won't happen when $\rho=2$),
and in either case together with $s_{j_0-1} \otimes s_{j_0-1}$ we get the 
desired six independent combinations.

Similarly, if the potential support of the $(j_0-1,j_0-1)$ row is
strictly to the right of $i_1$, 
we will have
$s_{j_0-1} \otimes s_{j_0}' = (j_0-1,j_0)$,
$s_{j_0}' \otimes s_{j_0}' = (j_0,j_0)$, 
$s_{j_0+1}' \otimes s_{j_0+1}' = (j_0+1,j_0+1)$, 
with $s_{j_0-1} \otimes s_{j_0+1}'$ having a nonzero $(j_0-1,j_0+1)$ part,
and
$s_{j_0}' \otimes s_{j_0+1}'' = (j_0-1,j_0)_{R}+ (j_0,j_0+1)$ and
$s_{j_0}' \otimes s_{j_0+1}' = (j_0-1,j_0)_{LC}+ (j_0,j_0+1)+(j_0,j_0)_{L'}$,
and we again obtain six independent combinations in the same manner.

If the potential support of the $(j_0-1,j_0-1)$ row is strictly between
$i_0$ and $i_1$, then none of the relevant rows can have support either
left of $i_0$ or right of $i_1$, and we get
$s_{j_0-1} \otimes s_{j_0}' = (j_0-1,j_0)$,
$s_{j_0-1} \otimes s_{j_0+1}'' = (j_0-1,j_0+1)$,
$s_{j_0}' \otimes s_{j_0}' = (j_0,j_0)$, 
$s_{j_0+1}'' \otimes s_{j_0+1}'' = (j_0+1,j_0+1)$, 
and
$s_{j_0}' \otimes s_{j_0+1}'' = (j_0,j_0+1)$.

If the $(j_0-1,j_0-1)$ row has disconnected potential support 
to the left of $i_0$ and strictly between $i_0$ and $i_1$, then
once again none of the relevant rows can have potential support extending 
right of $i_1$, and because $\rho=2$ we must have 
$a^{i_0}_{(j_0-1,j_0-1)}=c_{i_0}-1$, so none of the other relevant rows can
have their potential support contained strictly left of $i_0$, either. 
Moreover, the
$(j_0-1,j_0)$ row must have potential support containing $i_0$, 
so $s_{j_0}' \otimes s_{j_0}''$ cannot have any $(j_0-1,j_0)$ part, and
its $(j_0,j_0)$ part must be nonzero.
We then find that
$s_{j_0-1} \otimes s_{j_0+1}'' = (j_0-1,j_0+1)$,
$s_{j_0+1}'' \otimes s_{j_0+1}'' = (j_0+1,j_0+1)$, and
$s_{j_0}' \otimes s_{j_0+1}'' = (j_0,j_0+1)$. 
If the $(j_0-1,j_0)$ part of $s_{j_0-1} \otimes s_{j_0}''$ is nonzero, then
these together with $s_{j_0-1} \otimes s_{j_0-1}$ give six independent
combinations. Otherwise, we must have 
$s_{j_0-1} \otimes s_{j_0}'' = (j_0-1,j_0-1)_{C}$, and we see that
$s_{j_0-1} \otimes s_{j_0}' = (j_0-1,j_0-1)_L+(j_0-1,j_0)$ gives a
sixth independent combination.

The situation is nearly the same if the 
$(j_0-1,j_0-1)$ row has disconnected potential support 
to the right of $i_1$ and strictly between $i_0$ and $i_1$. Here we
instead obtain that $(j_0-1,j_0+1)$ must have potential support containing
$i_1$, and thus that
$s_{j_0-1} \otimes s_{j_0}' = (j_0-1,j_0)$,
$s_{j_0}' \otimes s_{j_0}' = (j_0,j_0)$, and
$s_{j_0}' \otimes s_{j_0+1}'' = (j_0,j_0+1)$, 
with $s_{j_0+1}' \otimes s_{j_0+1}''$ having nonzero $(j_0+1,j_0+1)$ part.
Then $s_{j_0-1} \otimes s_{j_0+1}'$ either has a nonzero $(j_0-1,j_0+1)$
part, or is equal to $(j_0-1,j_0-1)_C$, and in either case we obtain a
sixth combination, from $s_{j_0-1} \otimes s_{j_0-1}$ or 
$s_{j_0-1} \otimes s_{j_0+1}'' = (j_0-1,j_0+1)+(j_0-1,j_0-1)_{R}$
respectively.

If $(j_0-1,j_0-1)$ has three components of potential support, necessarily
left of $i_0$, strictly between $i_0$ and $i_1$, and right of $i_1$, then
none of the relevant rows other than $(j_0-1,j_0-1)$ can have potential
support contained strictly left of $i_0$ or strictly right of $i_1$, and
we also know that the potential support of the $(j_0-1,j_0)$ (respectively,
$(j_0-1,j_0+1)$) row contains $i_0$ (respectively, $i_1$).
 We then have that
$s_{j_0}' \otimes s_{j_0+1}'' = (j_0,j_0+1)$, and that 
$s_{j_0}' \otimes s_{j_0}''$ and $s_{j_0+1}' \otimes s_{j_0+1}''$ 
have nonzero $(j_0,j_0)$ and $(j_0+1,j_0+1)$ parts, respectively.
We also have 
$s_{j_0-1} \otimes s_{j_0}' = (j_0-1,j_0-1)_L+(j_0-1,j_0)$,
$s_{j_0-1} \otimes s_{j_0+1}'' = (j_0-1,j_0-1)_R+(j_0-1,j_0+1)$, and
$s_{j_0-1} \otimes s''' = (j_0-1,j_0)+(j_0-1,j_0-1)_C + (j_0-1,j_0+1)$.
To have a dependence between these, we need (at least one of) 
$s_{j_0-1} \otimes s_{j_0}' =(j_0-1,j_0)$ or
$s_{j_0-1} \otimes s_{j_0+1}'' = (j_0-1,j_0+1)$.
On the other hand, to have a dependence between the first five and
$s_{j_0-1} \otimes s_{j_0-1}$, we need 
$s_{j_0-1} \otimes s_{j_0}' =(j_0-1,j_0-1)_L$ or
$s_{j_0-1} \otimes s_{j_0+1}'' = (j_0-1,j_0-1)_R$.
If 
$s_{j_0-1} \otimes s_{j_0}' =(j_0-1,j_0)$ and
$s_{j_0-1} \otimes s_{j_0+1}'' = (j_0-1,j_0-1)_R$, we see that
$s_{j_0-1} \otimes s_{j_0+1}'$ must have a nonzero $(j_0-1,j_0-1)_{LC}$
or $(j_0-1,j_0+1)$ part, and thus gives a sixth independent combination.
On the other hand, if 
$s_{j_0-1} \otimes s_{j_0}' =(j_0-1,j_0-1)_L$ and
$s_{j_0-1} \otimes s_{j_0+1}'' = (j_0-1,j_0+1)$, we see that 
$s_{j_0-1} \otimes s_{j_0}''$ must have a nonzero $(j_0-1,j_0-1)_{CR}$ or
$(j_0-1,j_0)$ part, and again gives a sixth independent combination.

It remains to analyze the case that $(j_0-1,j_0-1)$ has two components
of potential support, one left of $i_0$, and the other right of $i_1$.
By hypothesis,
we only have to address the case that
$a^{i_0}_{(j_0-1,j_0-1)}=c_{i_0}-2$ and
$a^{i_1+1}_{(j_0-1,j_0-1)}=c_{i_1+1}+2$, and 
that we have degree $2$ in both $i_0$ and $i_1$.
In this situation, the $(j_0-1,j_0)$
row has potential support strictly left of $i_0$, but none of the other 
relevant rows do, and the $(j_0,j_0)$ row must have support containing
$i_0$ and extending left to at least the previous genus-$1$ component.
Similarly, the $(j_0-1,j_0+1)$ row has potential support 
strictly right of $i_1$, but none of the other relevant rows do, and the
$(j_0+1,j_0+1)$ row has support containing $i_1$ and extending to the 
right to at least the next genus-$1$ component.
We also see that the potential support of $(j_0,j_0+1)$ must be contained 
between $i_0$ and $i_1$ inclusive, and cannot be equal solely to $i_0$
or to $i_1$.
In particular, $s_{j_0}' \otimes s_{j_0+1}''$ cannot have a
$(j_0-1,j_0)$ or $(j_0-1,j_0+1)$ part, so must be equal to $(j_0,j_0+1)$.

Now, $s_{j_0-1} \otimes s_{j_0-1}=(j_0-1,j_0-1)_L$ because
$X_0$ is left-weighted, and we begin by considering the case that
no tensor has a $(j_0-1,j_0-1)_R$ part. Then we must
have 
$s_{j_0-1} \otimes s_{j_0+1}'' = (j_0-1,j_0+1)$,
$s_{j_0+1}'' \otimes s_{j_0+1}'' = (j_0+1,j_0+1)$, 
$s_{j_0}'' \otimes s_{j_0}''=(j_0,j_0)$,
and we also see that
$s_{j_0-1} \otimes s_{j_0}''$ must be $(j_0-1,j_0)$, because it could
only have a $(j_0-1,j_0+1)_{R'}$ part if the $j_0$ part of $s_{j_0}''$ 
extends through $i_1$, and in this case the fact that $X_0$ is left-weighted
gives us that $s_{j_0-1} \otimes s_{j_0}''=(j_0-1,j_0)$ regardless.
Thus, we obtain the desired six independent combinations in this case.

On the other hand, if any tensor has a $(j_0-1,j_0-1)_R$ part, we need
to produce only three more independent combinations, and we 
consider the four tensors
$s_{j_0}' \otimes s_{j_0}'' = (j_0-1,j_0)+(j_0,j_0)$,
$s_{j_0+1}' \otimes s_{j_0+1}'' = (j_0-1,j_0+1)+(j_0+1,j_0+1)$,
$s_{j_0-1} \otimes s''' = (j_0-1,j_0)+(j_0-1,j_0+1)$, and
$s''' \otimes s''' = (j_0,j_0)+(j_0+1,j_0+1)$.
These must have at least a three-dimensional span unless 
they collapse into equal pairs, and there are two possibilities for this:
either $s_{j_0}' \otimes s_{j_0}'' =s_{j_0-1} \otimes s''' = (j_0-1,j_0)$
and
$s_{j_0+1}' \otimes s_{j_0+1}'' =s''' \otimes s''' = (j_0+1,j_0+1)$,
or
$s_{j_0}' \otimes s_{j_0}'' = s''' \otimes s''' = (j_0,j_0)$ and
$s_{j_0+1}' \otimes s_{j_0+1}'' = s_{j_0-1} \otimes s''' = (j_0-1,j_0+1)$.
Moreover, Proposition \ref{prop:mixed-controlled} implies that the 
$j_0$-part of $s_{j_0}'$ doesn't
contain any genus-$1$ components left of $i_0$. Then we necessarily have
$s_{j_0}' \otimes s_{j_0}'' = (j_0-1,j_0)$, so only the first possibility
above can occur. Now, in general we have
$s_{j_0}'' \otimes s''' = (j_0,j_0)+(j_0-1,j_0+1)_{CR} + (j_0+1,j_0+1)_{R'}
+(j_0-1,j_0-1)_C +(j_0-1,j_0)_{CR}+(j_0,j_0+1)$, which in our case
simplifies to 
$s_{j_0}'' \otimes s''' = (j_0,j_0)+(j_0-1,j_0+1)_{CR} +(j_0,j_0+1)
+(j_0+1,j_0+1)_{R'}$.

If this has nonzero $(j_0,j_0)$ or $(j_0-1,j_0+1)$ term, we have our
sixth independent combination. On the other hand, if the $(j_0+1,j_0+1)$ term
is nonzero, the $(j_0,j_0+1)$ term must also be. 
Because the potential support of $(j_0,j_0+1)$ must end no later than $i_1$
and cannot be supported solely at $i_1$,
if the $(j_0,j_0+1)$ term of 
$s_{j_0}'' \otimes s'''$ is nonzero, this means that the $j_0$ part of
$s_{j_0}''$ must extend to cover all of $(j_0,j_0+1)$ (note that the proof
of Lemma \ref{lem:mixed-in-span} indicates that a $(j_0,j_0+1)$ part has to 
come from either a $j_0$ part of $s_{j_0}''$ and a $(j_0+1)$ part of $s'''$
or vice versa, but not some mixture of the two). But we know that
this contains at least one genus-$1$ component strictly right of $i_0$,
so since the support of $(j_0,j_0)$ ends at $i_0$, and $X_0$ is 
left-weighted, we conclude that we would have to have 
$s_{j_0}'' \otimes s_{j_0}''=(j_0,j_0)$ in this case. Thus, in all cases
we obtain the desired six independent combinations. 
\end{proof}

We can now prove the genus-$23$ case of our main theorem. As with the
genus-$22$ case, we phrase the result more generally to apply to other
$\rho=2$ cases in the future.

\begin{thm}\label{thm:rho-2} Fix $g,r,d$ with $r \geq 3$ and $\rho = 2$.
In characteristic $0$, suppose that for every left-weighted $X_0$ of genus 
$g$ as in Situation \ref{sit:chain-2}, and every
refined limit $\fg^r_d$ on $X_0$, there is an unimaginative 
$w=(c_2,\dots,c_N)$ such that the potentially appearing sections in 
multidegree $\md(w)$ are
linearly independent, and satisfying the following additional conditions:
\begin{ilist}
\itm if the limit $\fg^r_d$ falls into the ``first $3$-cycle'' situation 
described in Proposition \ref{prop:3-cycle-1}, we require that the
$(j_0-1,j_0)$ row has a unique potentially appearing section in multidegree 
$\md(w)$, whose support does not contain $i_0$ or $i_1$;
\itm if the limit $\fg^r_d$ falls into the ``second $3$-cycle'' situation 
described in Proposition \ref{prop:3-cycle-2}, we require that one of the
following three conditions is satisfied:
\begin{enumerate}
\item the $(j_0-1,j_0-1)$ row does not have potentially appearing sections 
both left of $i_0$ and right of $i_1$; or
\item
$2a^{i_0}_{j_0-1}= c_{i_0}-1, \text{ and } 2a^{i_1+1}_{j_0-1}= c_{i_1+1}+1$; or
\item
$2a^{i_0}_{j_0-1}= c_{i_0}-2, \text{ and } 2a^{i_1+1}_{j_0-1}= c_{i_1+1}+2$,
and $w$ has degree $2$ in both $i_0$ and $i_1$.
\end{enumerate}
\end{ilist}

Then the strong maximal rank conjecture holds for $(g,r,d)$, and more
specifically, a general curve of genus $g$ does not
have any $\fg^r_d$ for which \eqref{eq:mult-map} is not
injective.
\end{thm}

\begin{proof} The proof is essentially the same as that of Theorem
\ref{thm:rho-1}, still using
Propositions \ref{prop:simple} and \ref{prop:single-swap-indep} to treat
the cases that our refined limit linear series has no swaps or one swap,
respectively,
and adding Propositions \ref{prop:repeat-swap-indep}, 
\ref{prop:disjoint-swaps-indep}, \ref{prop:3-cycle-1-indep} and 
\ref{prop:3-cycle-2-indep} to address the cases with two swaps.
Using Remark \ref{rem:rho-meaning}, these are the only possibilities, since 
for $\rho=2$ we cannot have swaps involved more than two rows in a single
column. The only other
difference is that because we assume $X_0$ is left-weighted, we are forced
to consider only special directions of approach to $X_0$ in 
$\overline{\cM}_g$. Recalling that being left-weighted 
is preserved under the insertions of genus-$0$ chains which occur when we 
base change and then blow up to resolve
the resulting singularities, we do however conclude that for suitable
smoothing families, the generic fiber cannot carry a $\fg^r_d$ for which
\eqref{eq:mult-map} is not injective, as desired.
\end{proof}

Putting Theorem \ref{thm:rho-2} 
together with Theorem \ref{thm:basic-r6} and Corollary \ref{cor:alt-multideg},
we immediately conclude the genus-$23$ case of Theorem \ref{thm:main}.

\begin{rem}\label{rem:generalize} In our arguments for the $g=23$ case, we
used the $\rho=2$ hypothesis in two distinct ways: first, to limit the
number of swaps occurring to two, but then also to control the behavior
of the rest of the limit linear series when two swaps did occur, for 
instance limiting the number of possibilities for rows having disconnected 
potential support. This may appear discouraging from the point of view of 
generalizing to cases with higher $\rho$, but as $\rho$ increases, one 
also obtains more flexibility in choosing multidegrees while still 
maintaining linear independence of the potentially appearing sections.
Indeed, we are taking advantage of this phenomenon already in the $\rho=2$ 
case with Corollary \ref{cor:alt-multideg}.
\end{rem}

\bibliographystyle{amsalpha}
\bibliography{gen}

\end{document}